\documentclass[11pt,reqno,a4paper]{amsart}
\usepackage[T1]{fontenc}
\usepackage[utf8]{inputenc}
\usepackage{lmodern}
\usepackage{amssymb,mathtools,mathrsfs}
\usepackage{enumitem,booktabs,array,longtable}
\usepackage[margin=26mm]{geometry}
\usepackage{microtype}
\usepackage{aliascnt}
\usepackage[colorlinks=true,linkcolor=blue,citecolor=blue,urlcolor=blue,
  bookmarksnumbered=true,pdfstartview=FitH]{hyperref}

\theoremstyle{plain}
\newtheorem{theorem}{Theorem}[section]
\newtheorem{lemma}[theorem]{Lemma}
\newtheorem{proposition}[theorem]{Proposition}
\newtheorem{corollary}[theorem]{Corollary}
\newtheorem{conjecture}[theorem]{Conjecture}

\theoremstyle{definition}
\newtheorem{definition}[theorem]{Definition}

\theoremstyle{remark}
\newtheorem{remark}[theorem]{Remark}
\numberwithin{equation}{section}
\hypersetup{pdftitle={Weitzenbock Curvature Estimates in Arbitrary Degree and Sectional-Scalar Curvature Pinching},
 pdfauthor={Shimeng Zhang},
 pdfsubject={Curvature estimates, differential forms, Betti numbers, and endpoint rigidity},
 pdfkeywords={Weitzenbock, sectional curvature, scalar curvature, Betti numbers, rigidity}}
\allowdisplaybreaks[1]
\title[Weitzenb\"ock estimates in arbitrary degree]
 {Weitzenb\"ock Curvature Estimates in Arbitrary Degree\\
 and Sectional--Scalar Curvature Pinching}
\author[S.~Zhang]{Shimeng Zhang}
\address{School of Mathematical Sciences, Laboratory of Mathematics
and Complex Systems, Beijing Normal University, Beijing 100875,
P.~R.~China}
\email{202531130043@mail.bnu.edu.cn}
\date{Revised November 23, 2026}

\begin{document}
\begin{abstract}
For algebraic curvature tensors $E$ with $\sec(E)\ge0$, we give an explicit coefficient $C$ in $q_p(E)\ge-C\operatorname{Scal}(E)\operatorname{Id}$. When the dimension $n$ and degree $p$ are even and $2\le p\le n-2$, the optimal coefficient is $p(n-p)/(3n)$, and we classify the nonzero equality pairs in terms of orthogonal complex structures and powers of their K\"ahler forms. For $n=5,6$ and $p=3$ the optimal coefficient is $1/3$. Under strongly nonnegative curvature we determine the optimal coefficient for every $n$ and $p$. On closed manifolds these estimates give explicit sectional--scalar pinching thresholds for cohomology vanishing. At the weak threshold harmonic forms in the corresponding degrees are parallel. For nonflat metrics the degree-three threshold forces $b_3=0$ when $n\ge5$. At the sharp threshold in even dimension $n$, a nonflat manifold with a positive interior even Betti number is isometric, up to scale, to Fubini--Study $\mathbb{CP}^{n/2}$ or its quotient by a free antiholomorphic isometric involution.
\end{abstract}
\maketitle
\tableofcontents

\section{Introduction}
\label{sec:introduction}

How much sectional--scalar curvature pinching forces the cohomology of a closed Riemannian manifold to vanish in a given degree? At a point $x$ of a Riemannian manifold $(M^n,g)$, write
\begin{equation}\label{geom:normalization}
    K_{\min}(x)=\min_{\sigma\in\mathrm{Gr}_2(T_xM)} K_g(\sigma), \qquad S_0(x)=\frac{\operatorname{Scal}_g(x)}{n(n-1)}.
\end{equation}
Thus $S_0(x)$ is the average of the sectional curvatures at $x$. The coefficient $(n-1)/(n+2)$ in Yau's sphere pinching problem is attained by the Fubini--Study metric in every even dimension $n\ge4$. Ge \cite[Theorem~1.2]{Ge} proved that the strict inequality with this coefficient forces $H^2(M;\mathbb R)$ to vanish in even dimensions, and obtained a smaller coefficient in odd dimensions. We study the same question for $p$-forms in degrees $1\le p\le n-1$. When $n$ and $p$ are both even, the Fubini--Study coefficient $(n-1)/(n+2)$ is sharp. Whenever at least one of them is odd, a strictly smaller coefficient suffices.

The algebraic input is an estimate on the Weitzenb\"ock curvature endomorphism $q_p(E)$ on $\Lambda^p V^*$, where $E$ is an algebraic curvature tensor on a Euclidean space $V$ of dimension $n\ge4$. With our normalization,
\[
    q_p(E)=\bigl(\operatorname{Ric}(E)^\sharp\bigr)^{[p]}-2\mathcal R_E^{(p)}, \qquad q_p(I)=p(n-p)\operatorname{Id},
\]
where $I$ has constant sectional curvature one. On the cone of tensors with nonnegative sectional curvature, Ge's sharp two-form estimate is
\begin{equation}\label{intro:ge}
    q_2(E)\ge-\frac{2(\ell-1)}{3\ell}\operatorname{Scal}(E)\operatorname{Id}, \qquad \ell=\lfloor n/2\rfloor.
\end{equation}
Ge applies the Ni--Wilking four-frame inequality \cite{NiWilking} in an orthogonal block frame for a two-form \cite[Theorem~4.1]{Ge}. In higher degrees the contractions of a form need not share such a frame. We keep the Ricci term produced by contraction and average over hyperplanes. The resulting identities give an induction in the degree and dimension. For the representation-theoretic background on these curvature terms see \cite{BettiolMendes}.

For $2\le p\le n-2$, put $c=n-p$ and set
\begin{equation}\label{intro:kappa}
    \kappa_{n,p}=
    \begin{cases}
        c/3, & c\text{ even},\\[2pt]
        \displaystyle\frac c3-\frac{n+2}{3\bigl((p-1)(c+3)+1\bigr)}, & c\text{ odd}.
    \end{cases}
\end{equation}
When $c=3$ and $p\ge3$, a dimension lift gives the improved value, which we use throughout,
\begin{equation}\label{intro:t}
    \kappa_{n,n-3}=t_n=\frac{19n-82}{21n-78}\qquad(n\ge6).
\end{equation}
Pairing the estimates in complementary degrees gives the coefficient
\begin{equation}\label{intro:C}
    C_{n,p}=\frac{\kappa_{n,p}\kappa_{n,n-p}}{\kappa_{n,p}+\kappa_{n,n-p}}, \qquad 2\le p\le n-2 .
\end{equation}
Set $C_{n,p}=0$ for $p=0,1,n-1,n$. We use $\widehat C_{n,p}=C_{n,p}$ with the single complementary pair of exceptions
\begin{equation}\label{intro:Chat}
    \widehat C_{7,3}=\widehat C_{7,4}=\frac{5168}{11253}.
\end{equation}

\begin{theorem}[Estimates for nonnegative sectional curvature]\label{intro:main}
    Let $E$ be an algebraic curvature tensor with nonnegative sectional curvature on an $n$-dimensional Euclidean space $V$, $n\ge4$. Then
    \[
        q_p(E)\ge-\widehat C_{n,p}\operatorname{Scal}(E)\operatorname{Id}_{\Lambda^p V^*} \qquad (0\le p\le n).
    \]
    If $n$ and $p$ are even and $2\le p\le n-2$, the coefficient is sharp and equals $p(n-p)/(3n)$. In these degrees, equality on a unit form $\omega$ for $E\ne0$ holds precisely for the pairs
    \[
        E=tE_J, \qquad \omega=\pm\frac{\Omega_J^{p/2}}{(p/2)!\sqrt{\binom{n/2}{p/2}}}, \qquad t>0,
    \]
    where $J$ is an orthogonal complex structure on $V$ and $\Omega_J$ is its K\"ahler form. In degree three the sharp coefficient is zero for $n=4$ and $1/3$ for $n=5,6$.
\end{theorem}

Here $E_J$ is the model tensor in \eqref{eq:EJ}. For a complex structure on all of $V$, $I+E_J$ is the Fubini--Study curvature tensor of holomorphic sectional curvature four.

To pass from this estimate to cohomology, split the Riemannian curvature tensor at each point as $\operatorname{Rm}=K_{\min} I+E$. The remainder has nonnegative sectional curvature and scalar curvature $n(n-1)(S_0-K_{\min})$. The threshold at which Theorem \ref{intro:main} makes $q_p(\operatorname{Rm})$ nonnegative is therefore
\begin{equation}\label{intro:beta}
    \widehat\beta_n^{(p)}=\frac{n(n-1)\widehat C_{n,p}}{p(n-p)+n(n-1)\widehat C_{n,p}}, \qquad 1\le p\le n-1 .
\end{equation}
For these degrees we have
\begin{equation}\label{geom:even-threshold}
    \begin{aligned}
        \widehat\beta_n^{(p)}&=\frac{n-1}{n+2} && \text{if both $n$ and $p$ are even},\\[2pt]
        \widehat\beta_n^{(p)}&<\frac{n-1}{n+2} && \text{if at least one of $n$ and $p$ is odd}.
    \end{aligned}
\end{equation}
Thus whenever $n$ or $p$ is odd, the following vanishing theorem requires \emph{weaker pinching} than $K_{\min}>\frac{n-1}{n+2}S_0$. In degree two it gives Ge's coefficient in \cite[Theorem~1.2]{Ge}.

\begin{theorem}[Cohomology vanishing]\label{intro:geometry}
    Let $(M^n,g)$ be a closed, connected Riemannian manifold of dimension $n\ge4$, and let $1\le p\le n-1$. Suppose
    \[
        K_{\min}(x)\ge\widehat\beta_n^{(p)}S_0(x) \qquad \text{for every }x\in M,
    \]
    with strict inequality at some point. Then $H^p(M;\mathbb R)=H^{n-p}(M;\mathbb R)=0$. When $n$ and $p$ are both even, the coefficient in this criterion is sharp.
\end{theorem}

For example, on $\mathbb{CP}^{m}$, $n=2m\ge4$, normalize the Fubini--Study metric by $K_{\min}=1$. Then $S_0=(n+2)/(n-1)$, and the powers $\Omega^a$ of its K\"ahler form generate its even-degree cohomology. Thus the same metric realizes the sharp endpoint in each interior even degree and satisfies the strict hypothesis of Theorem \ref{intro:geometry} in every odd degree.

At the weak endpoint the Bochner formula still forces harmonic forms to be parallel. This part of Ge's weak endpoint theorem \cite[Theorem~1.3]{Ge} extends to every degree.

\begin{theorem}[Parallel forms at the weak endpoint]\label{intro:parallel}
    Let $(M^n,g)$ be a closed, connected Riemannian manifold of dimension $n\ge4$, and let $1\le p\le n-1$. If $K_{\min}\ge\widehat\beta_n^{(p)}S_0$ everywhere on $M$, then all harmonic forms of degrees $p$ and $n-p$ are parallel. For $k\in\{p,n-p\}$ and every $x\in M$,
    \[
        b_k(M)=\dim\bigl(\Lambda^k T_x^*M\bigr)^{\operatorname{Hol}_x(g)}\le\binom nk ,
    \]
    where $\operatorname{Hol}_x(g)$ denotes the full holonomy group.
\end{theorem}

No orientability hypothesis is needed in these two theorems. For a closed flat manifold $M=\mathbb R^n/\Gamma$ the formula takes the form $b_k(M)=\dim(\Lambda^k(\mathbb R^n)^*)^F$, where $F$ is the finite linear holonomy of $\Gamma$. In particular $b_3(M)>0$ precisely when $F$ fixes a nonzero three-form, and flat tori attain the binomial bound in every degree.

In even dimensions Ge's rigidity theorem \cite[Theorem~1.3]{Ge} identifies the nonflat endpoint with $b_2>0$ as Fubini--Study space. For a parallel form of degree $2a$, the relevant question on the universal cover is which deck transformations preserve $\Omega^a$.

\begin{theorem}[Even-degree rigidity]\label{geom:even-rigidity}
    Let $(M^n,g)$ be a closed, connected Riemannian manifold of dimension $n=2m\ge4$. Suppose $1\le a\le m-1$, $b_{2a}(M)>0$, and $K_{\min}\ge\frac{n-1}{n+2}S_0$ everywhere. If $g$ is nonflat, then $b_{2a}(M)=1$ and, up to scaling, $(M,g)$ is isometric to Fubini--Study $\mathbb{CP}^m$ or its quotient by a fixed-point-free antiholomorphic isometric involution. The quotient occurs only for odd $m$ and even $a$.
\end{theorem}

For $m=2s+1$ the involution $\sigma$ is unique up to conjugacy by an isometry. In homogeneous coordinates it is induced by
\[
    (z_1,z_2,\cdots,z_{2s+1},z_{2s+2})\longmapsto(-\bar z_2,\bar z_1,\cdots,-\bar z_{2s+2},\bar z_{2s+1}).
\]
Since $\sigma^*\Omega=-\Omega$, the quotient is nonorientable, its odd Betti numbers vanish, and $b_{2r}$ equals one for even $r$ and zero for odd $r$ ($0\le r\le m$). For instance, $\mathbb{CP}^3/\langle\sigma\rangle$ has $b_4=1$ and $b_2=0$. This explains why the two complementary degrees in Theorem \ref{intro:parallel} have separate holonomy formulas.

A parallel three-form $\omega$ on a manifold with irreducible holonomy satisfies $\langle\iota_X\omega,\iota_Y\omega\rangle=(3|\omega|^2/n)g(X,Y)$. For instance, the fundamental three-form of $G_2$ satisfies this identity. To rule out such a form at a nonflat endpoint we first examine which holonomy groups a positive sectional--scalar pinching coefficient allows.

\begin{theorem}[Holonomy under sectional--scalar pinching]\label{geom:three-holonomy}
    Let $(M^n,g)$ be a closed, connected, nonflat Riemannian manifold of dimension $n\ge4$. Suppose $K_{\min}\ge\beta S_0$ everywhere for some constant $0<\beta<1$. The universal Riemannian cover has no nonzero parallel three-form. If $n$ is odd, the restricted holonomy group is $\mathrm{SO}(n)$. In any dimension $n\ge4$, the additional hypothesis $q_3(\operatorname{Rm})\ge0$ everywhere implies $H^3(M;\mathbb R)=H^{n-3}(M;\mathbb R)=0$.
\end{theorem}

The explicit pinching coefficient for three-forms gives the curvature endomorphism bound in this theorem. This yields a degree-three analogue of Ge's weak endpoint theorem \cite[Theorem~1.3]{Ge}.

\begin{theorem}[The degree-three weak endpoint]\label{geom:three-rigidity}
    Let $(M^n,g)$ be a closed, connected Riemannian manifold of dimension $n\ge5$. If $K_{\min}\ge\widehat\beta_n^{(3)}S_0$ everywhere, then either $g$ is flat or $H^3(M;\mathbb R)=H^{n-3}(M;\mathbb R)=0$.
\end{theorem}

In odd dimensions, $\mathrm{SO}(n)$ fixes no nonzero forms of degrees $1,\cdots,n-1$. In particular, for a nonflat metric and $2\le p\le n-2$, the weak pinching bound in degree $p$ already forces $H^p(M;\mathbb R)=H^{n-p}(M;\mathbb R)=0$ by Theorems \ref{intro:parallel} and \ref{geom:three-holonomy}. Taking the largest of the degree-dependent thresholds therefore gives a rational homology sphere theorem. In dimension five it extends Ge's strict pinching result \cite[Corollary~1.4]{Ge} to the nonflat weak endpoint.

\begin{corollary}[Rational homology spheres in odd dimensions]\label{geom:odd-rational-sphere}
    Let $(M^n,g)$ be a closed, connected Riemannian manifold of odd dimension $n=2m+1\ge5$, and set
    \[
        B_n=\max_{2\le p\le m}\widehat\beta_n^{(p)} .
    \]
    If $g$ is nonflat and $K_{\min}\ge B_n S_0$ everywhere, then $H_*(M;\mathbb Q)\cong H_*(S^n;\mathbb Q)$. In dimension three, a closed, connected Riemannian manifold with $K_{\min}>0$ everywhere is a rational homology three-sphere.
\end{corollary}

By \eqref{geom:even-threshold}, $B_n<(n-1)/(n+2)$ for every odd $n\ge5$. In dimension three the strict statement has $B_3=0$. The product $S^2\times S^1$ with its standard product metric has $K_{\min}=0$ and $b_1=1$, so a weak inequality with this coefficient would include nonflat manifolds with nonzero first cohomology.

We close the introduction by recording what is known about the optimal \emph{algebraic} coefficient. Let $\gamma_{n,p}^*$ be the least $C\ge0$ for which $q_p(E)\ge-C\operatorname{Scal}(E)\operatorname{Id}$ whenever $\sec(E)\ge0$. Put $u=\lfloor p/2\rfloor$, $v=\lfloor(n-p)/2\rfloor$, and $D_{n,p}=2uv/(3(u+v))$, interpreted as zero if $u+v=0$. Complex-projective model tensors extended by zero give $D_{n,p}\le\gamma_{n,p}^*\le\widehat C_{n,p}$. On the cone of \emph{strongly nonnegative} curvature tensors, for which $\mathcal R_E+\eta$ is positive semidefinite for some four-form $\eta$, the sharp coefficient is $D_{n,p}$ in every degree. In dimensions at least five this cone is strictly smaller than the cone of tensors with nonnegative sectional curvature \cite{BKM}.

Section \ref{sec:conventions} fixes the curvature conventions. Sections \ref{alg:section}--\ref{sec:seven-dimensions} prove the algebraic bounds, including the seven-dimensional refinement. Section \ref{strong:section} treats strongly nonnegative curvature. Section \ref{eqcase:section} classifies the algebraic equality cases and completes the proof of Theorem~\ref{intro:main}. Section \ref{geom:section} proves the geometric results above. Section \ref{sec:open} discusses the remaining optimal constants, and Appendix \ref{sec:block-family} gives a sharp estimate for a further family of three-forms.

All arguments in this paper were developed independently of the concurrent work \cite{LMP}. After the completion of this paper, we noticed that Corollary \ref{geom:odd-rational-sphere} incidentally provides an answer to \cite[Remark~5.6]{LMP}.

\medskip
\noindent\textbf{Acknowledgment.} I thank my advisor, Jian Ge, for suggesting that I study this problem.

\smallskip
\noindent\textbf{AI disclosure.} I developed the proof of Lemma \ref{six:normal-form} and the methods in Section \ref{sec:seven-dimensions} and Appendix \ref{sec:block-family} through discussions with OpenAI's ChatGPT. I take responsibility for the mathematical arguments and the text.

\section{Preliminaries}\label{sec:conventions}
Throughout, $V$ is a real Euclidean space of dimension $n$.

\subsection{Algebraic curvature operators and the Weitzenb\"ock endomorphism}
An algebraic curvature tensor $E$ is skew in either pair, symmetric under pair interchange, and satisfies the first Bianchi identity $E_{ijkl}+E_{jkil}+E_{kijl}=0$. See \cite[Definition~2.1]{Ge}.

Our geometric convention is
\[
	R^{\nabla}(X, Y)Z = \nabla_{X}\nabla_{Y}Z-\nabla_{Y}\nabla_{X}Z-\nabla_{[X, Y]}Z, \qquad \operatorname{Rm}(X, Y, Z, W) = \langle R^\nabla(X,Y)W,Z\rangle.
\]
Thus in an orthonormal basis $\{e_i\}$, $K_{ij}=E_{ijij}$ for $i\ne j$, $\operatorname{Ric}(E)_{ij}=\sum_k E_{ikjk}$, and $S=\operatorname{Scal}(E)=2\sum_{i<j}K_{ij}$. Define $\mathrm{Ric}(E)^{\sharp}\in \mathrm{End}(V)$ by $\langle \mathrm{Ric}(E)^{\sharp}X, Y \rangle = \mathrm{Ric}(E)(X,Y)$.

Recall that the inner product on $\Lambda^2V^*$ is given by $\langle e^{ij},e^{kl}\rangle=\langle e^i\wedge e^j,e^k\wedge e^l\rangle=\delta^{ik}\delta^{jl}-\delta^{il}\delta^{jk}$, hence if $\omega=\sum_{i<j}\omega_{ij}e^i\wedge e^j$ and $\eta=\sum_{i<j}\eta_{ij}e^i\wedge e^j$, then $\langle \omega,\eta\rangle=\sum_{i<j}\omega_{ij}\eta_{ij}$. The self-adjoint curvature operator $\mathcal R_E:\Lambda^2V^*\to\Lambda^2V^*$ is determined by $\langle\mathcal R_E e^{ij},e^{kl}\rangle=E_{ijkl}$. Thus,
\[
    \mathcal R_E(e^{ij})=\sum_{k<l}E_{ijkl}e^{kl}=\frac12\sum_{k,l}E_{ijkl}e^{kl}.
\]
The first Bianchi identity gives
\[
    \mathcal R_E(e^{ij}) = \sum_{k<l}E_{ijkl}e^{kl} = \sum_{k<l}(E_{ikjl}-E_{iljk})e^{kl} = \sum_{k,l}E_{ikjl}e^{kl}.
\]
Let $\omega=\sum_{k<l}\omega_{kl}e^{k}\wedge e^{l}=\frac12\sum_{k,l}\omega_{kl}e^{k}\wedge e^{l}$ with $\omega_{kl}=-\omega_{lk}$, then
\[
    (\mathcal{R}_{E}\omega)_{ij} = \langle \mathcal{R}_{E}\omega, e^{ij} \rangle = \Bigl\langle \frac{1}{2} \sum_{k,l}\omega_{kl}\mathcal{R}_{E}(e^{kl}), e^{ij} \Bigr\rangle = \frac{1}{2} \sum_{k,l}\omega_{kl}E_{ijkl}.
\]
Similarly, the first Bianchi identity gives
\[
    (\mathcal{R}_{E}\omega)_{ij} = \frac12\sum_{k,l}E_{ijkl}\omega_{kl} = \sum_{k<l}E_{ijkl}\omega_{kl} = \sum_{k<l}(E_{ikjl}-E_{iljk})\omega_{kl} = \sum_{k,l}E_{ikjl}\omega_{kl}.
\]

\begin{definition}\label{def:ExtCurv}
    The natural extension of $\mathcal R_E$ to $p$-forms, denoted by $\mathcal R_E^{(p)}$, is defined by
    \[
        \mathcal R_E^{(p)}\omega := \sum_{i<j}\mathcal R_E(e^{ij})\wedge \iota_{j}\iota_{i}\omega = \sum_{i<j,\ a<b}E_{ijab}\,\varepsilon_a\varepsilon_b\iota_j\iota_i\omega = \sum_{i<j,\ a<b}E_{ijab}\,\varepsilon_i\varepsilon_j\iota_b\iota_a\omega,
    \]
    where $\omega\in \Lambda^p V^*$, $\iota_{i}:=\iota_{e_i}$ is the interior multiplication, and $\varepsilon_i:=e^i\wedge$.
\end{definition}

The form $\mathcal R_E^{(p)}\omega$ is again a $p$-form. For $p=2$, we have $\mathcal R_E^{(2)}=\mathcal R_E$. Since $\varepsilon_i\varepsilon_j=-\varepsilon_j\varepsilon_i$ and $\iota_i\iota_j=-\iota_j\iota_i$,
\[
    \mathcal R_E^{(p)}=\frac14\sum_{i,j,a,b}E_{ijab}\,\varepsilon_i\varepsilon_j\iota_b\iota_a=\frac14\sum_{i,j,a,b}(E_{iajb}-E_{ibja})\varepsilon_i\varepsilon_j\iota_b\iota_a = \frac12\sum_{i,j,a,b}E_{iajb}\,\varepsilon_i\varepsilon_j\iota_b\iota_a.
\]
Note that $\mathcal R_E^{(p)}$ is self-adjoint. Indeed, the adjoint relations $\varepsilon_i^*=\iota_i$ and $\iota_i^*=\varepsilon_i$ give $(\varepsilon_i\varepsilon_j\iota_b\iota_a)^* = \iota_a^*\iota_b^*\varepsilon_j^*\varepsilon_i^* = \varepsilon_a\varepsilon_b\iota_j\iota_i$. Hence $(\mathcal R_E^{(p)})^* = \mathcal R_E^{(p)}$.

\begin{definition}
    Let $B\in \mathrm{End}(V)$. For $0\le p\le\dim V$, define $B^{[p]}:\Lambda^pV^*\to\Lambda^pV^*$ by
    \[
        \bigl(B^{[p]}\omega\bigr)(X_1,\cdots,X_p) = \sum_{i=1}^p \omega(X_1,\cdots,BX_i,\cdots,X_p).
    \]
    In particular, $B^{[0]}=0$.
\end{definition}

\begin{lemma}
    If $\{e_i\}$ is an orthonormal basis with dual basis $\{e^i\}$ and $Be_i=\sum_j B_i^j e_j$, then:
    \begin{enumerate}
        \item $B^{[p]}$ maps $p$-forms to $p$-forms, and $B^{[p]}=\sum_{i,j}B_i^j\,\varepsilon_i\iota_j$.
        \item If $B$ is self-adjoint, then $B^{[p]}$ is also self-adjoint.
        \item The operator $B^{[p]}$ is a derivation: for $\omega\in\Lambda^pV^*$ and $\eta\in\Lambda^qV^*$,
        \[
            B^{[p+q]}(\omega\wedge\eta)=(B^{[p]}\omega)\wedge\eta+\omega\wedge(B^{[q]}\eta).
        \]
        \item If $B$ is self-adjoint and positive semidefinite, then $B^{[p]}$ is positive semidefinite. If $B$ is positive definite, then $B^{[p]}$ is positive definite for $1\le p\le\dim V$.
    \end{enumerate}
\end{lemma}

\begin{proof}
    The definition is alternating in $X_1,\cdots,X_p$, since each term is obtained by inserting $B$ into one argument of the alternating form $\omega$. Thus $B^{[p]}\omega$ is again a $p$-form. Evaluating on basis vectors gives $(B^{[p]}\omega)_{i_1\cdots i_p} = \sum_{a=1}^p B_{i_a}^{m}\omega_{i_1\cdots i_{a-1}\,m\,i_{a+1}\cdots i_p}$. On the other hand,
    \[
        \bigl(\varepsilon_i\iota_j\omega\bigr)(e_{i_1},\cdots,e_{i_p}) = \sum_{a=1}^p\delta_{i_a}^i \, \omega(e_{i_1},\cdots,e_{i_{a-1}},e_j,e_{i_{a+1}},\cdots,e_{i_p}).
    \]
    Therefore $(\sum_{i,j}B_i^j\varepsilon_i\iota_j\omega)_{i_1\cdots i_p} = \sum_{a=1}^p B_{i_a}^j \, \omega_{i_1\cdots i_{a-1}\,j\,i_{a+1}\cdots i_p}$, which are precisely the components of $B^{[p]}\omega$. This proves the operator formula.

    If $B$ is self-adjoint, then $B_i^j=B_j^i$ and $(B^{[p]})^* = \sum_{i,j}B_i^j\,\iota_j^*\varepsilon_i^* = \sum_{i,j}B_j^i\,\varepsilon_j\iota_i = B^{[p]}$. The derivation rule follows from that $B^{[p]}=\sum_{i,j}B_i^j\varepsilon_i\iota_j$ and that $\varepsilon_i\iota_j$ is a derivation:
    \begin{equation}\label{eq:Leibniz}
        \varepsilon_i\iota_j(\omega\wedge\eta) = (\varepsilon_i\iota_j\omega)\wedge\eta + \omega\wedge(\varepsilon_i\iota_j\eta).
    \end{equation}

    Finally, suppose \(B\) is self-adjoint and positive semidefinite. Choose an orthonormal eigenbasis with $Be_i=\lambda_i e_i$ and $\lambda_i\ge0$. Then $B^{[p]}=\sum_{i,j}B_i^j\varepsilon_i\iota_j=\sum_i\lambda_i\varepsilon_i\iota_i$. Let $i_1<\cdots<i_p$. By \eqref{eq:Leibniz},
    \[
        \varepsilon_i\iota_i\bigl(e^{i_1}\wedge\cdots\wedge e^{i_p}\bigr)=
        \begin{cases}
            e^{i_1}\wedge\cdots\wedge e^{i_p}, & i\in\{i_1,\cdots,i_p\},\\[2mm]
            0, & i\notin\{i_1,\cdots,i_p\}.
        \end{cases}
    \]
    Therefore $B^{[p]}\bigl(e^{i_1}\wedge\cdots\wedge e^{i_p}\bigr)=(\lambda_{i_1}+\cdots+\lambda_{i_p})e^{i_1}\wedge\cdots\wedge e^{i_p}$. Hence $B^{[p]}$ is diagonal with respect to the orthonormal basis $\{e^{i_1}\wedge\cdots\wedge e^{i_p} \mid i_1<\cdots<i_p\}$, with eigenvalues $\lambda_{i_1}+\cdots+\lambda_{i_p}\ge 0$. If $B$ is positive definite, all $\lambda_i>0$, so every such eigenvalue is positive for $p\ge1$.
\end{proof}

\begin{definition}\label{def:WeitCurv}
	For an algebraic curvature tensor $E$, define its \emph{Weitzenb\"{o}ck curvature endomorphism} acting on $p$-forms by
	\begin{equation}\label{eq:WeitCurv1}
		q_{p}(E):= \left( \mathrm{Ric}(E)^{\sharp} \right)^{[p]}- 2 \mathcal{R}_{E}^{(p)}.
	\end{equation}
\end{definition}

Both terms are self-adjoint, so $q_p(E)$ is self-adjoint. All operator inequalities below are inequalities of quadratic forms on $\Lambda^pV^*$. Also, $E\mapsto q_p(E)$ is linear.

The Weitzenb\"ock endomorphism has the representation formula
\begin{equation}\label{eq:Casimir-formula}
    q_p(E)=-\sum_{i<j,\,k<l}E_{ijkl}\rho_{ij}\rho_{kl}, \qquad \rho_{ij}=\varepsilon_j\iota_i-\varepsilon_i\iota_j.
\end{equation}
Indeed, expansion and the antisymmetry of $E$ give
\[
    \sum_{i<j,\,k<l}E_{ijkl}\rho_{ij}\rho_{kl} = \frac14\sum_{i,j,k,l}E_{ijkl}\rho_{ij}\rho_{kl} = \sum_{i,j,k,l}E_{ijkl}\varepsilon_j\iota_i\varepsilon_l\iota_k = \sum_{i,j,k}E_{ijki}\varepsilon_j\iota_k-\sum_{i,j,k,l}E_{ijkl}\varepsilon_j\varepsilon_l\iota_i\iota_k.
\]
The first term on the right is $-(\operatorname{Ric}(E)^\sharp)^{[p]}$ and the second term is $2\mathcal R_E^{(p)}$.

The Weitzenb\"ock formula is
\begin{equation}\label{eq:WZBKFormula}
	\Delta_{\mathrm{H}}\omega = \nabla^*\nabla \omega + q_p(\mathrm{Rm})\omega,
\end{equation}
where $\Delta_{\mathrm{H}}$ is the Hodge Laplacian and $\nabla^{*}\nabla$ is the nonnegative rough Laplacian.

We briefly check that this algebraic definition has the stated geometric sign. At a point $x$, choose an orthonormal frame $\{e_i\}$ that is synchronous, i.e., at $x$, $\nabla_ie_j=0$ for all $i$ and $j$. Then $d=\sum_i\varepsilon_i\nabla_i$, $\delta=-\sum_j\iota_j\nabla_j$ and $\nabla_i\varepsilon_j=\varepsilon_j\nabla_i$, $\nabla_i\iota_j=\iota_j\nabla_i$. The rough Laplacian $\nabla^*\nabla$ in a local orthonormal frame is $\nabla^*\nabla = -\sum_i\bigl(\nabla_i\nabla_i-\nabla_{\nabla_i e_i}\bigr)$. At $x$, $\nabla^*\nabla=-\sum_i\nabla_i\nabla_i$. These together with $\iota_j\varepsilon_i+\varepsilon_i\iota_j=\delta_{ij}\operatorname{Id}$ imply
\[
    \Delta_{\mathrm H} = d\delta+\delta d = \nabla^*\nabla - \sum_{i,j}\varepsilon_i\iota_j R^{\Lambda}(e_i,e_j),
\]
where $R^{\Lambda}$ is the curvature of the induced connection on forms. We claim that the chosen convention for $\mathrm{Rm}$ gives
\[
    R^{\Lambda}(e_i,e_j)=\sum_{k,l}\mathrm{Rm}_{ijkl}\varepsilon_k\iota_l.
\]
Indeed, $R^{\Lambda}(e_i,e_j)$ and $\varepsilon_k\iota_l$ are both derivations. Therefore, it suffices to verify their action on $1$-forms, since their action on all $p$-forms is uniquely determined by the Leibniz rule. For any $1$-form $\omega$ and vector field $Z$, we have the standard duality relation $(R^\Lambda(X,Y)\omega)(Z) = -\omega\bigl(R^\nabla(X,Y)Z\bigr)$, which gives the claimed formula. Substituting these gives
\[
    -\sum_{i,j}\varepsilon_i\iota_jR^{\Lambda}(e_i,e_j) = -\sum_{i,j,k,l}\mathrm{Rm}_{ijkl}\,\varepsilon_i\iota_j\varepsilon_k\iota_l = -\sum_{i,j,k,l}\mathrm{Rm}_{ijkl}\,\varepsilon_i(\delta_{jk}\mathrm{Id}-\varepsilon_k\iota_j)\iota_l.
\]
The right is $-\sum_{i,j,l}\mathrm{Rm}_{ijjl}\,\varepsilon_i\iota_l+\sum_{i,j,k,l}\mathrm{Rm}_{ijkl}\,\varepsilon_i\varepsilon_k\iota_j\iota_l$, where the first term is $\sum_{i,l}\mathrm{Ric}_{il}\,\varepsilon_i\iota_l = (\mathrm{Ric}^\sharp)^{[p]}$. For the second term, using $\iota_a\iota_b=-\iota_b\iota_a$ and the Bianchi identity,
\begin{align*}
    \sum_{i,j,k,l}\mathrm{Rm}_{ijkl}\,\varepsilon_i\varepsilon_k\iota_j\iota_l &= \sum_{i,a,j,b}\mathrm{Rm}_{iajb}\,\varepsilon_i\varepsilon_j\iota_a\iota_b = \frac12\sum_{i,a,j,b}(\mathrm{Rm}_{iajb}-\mathrm{Rm}_{ibja})\varepsilon_i\varepsilon_j\iota_a\iota_b\\
    &=\frac12\sum_{i,j,a,b}\mathrm{Rm}_{ijab}\,\varepsilon_i\varepsilon_j\iota_a\iota_b =-2\sum_{i<j,\ a<b}\mathrm{Rm}_{ijab}\,\varepsilon_i\varepsilon_j\iota_b\iota_a = -2\mathcal R_{\mathrm{Rm}}^{(p)}.
\end{align*}
Thus all factors agree with \eqref{eq:WZBKFormula}.

\begin{lemma}\label{lem:HodgeDual}
	Let $*:\Lambda^p\to\Lambda^{n-p}$ be the Hodge star operator, then $*\,q_p=q_{n-p}\,*$.
\end{lemma}

\begin{proof}
	For $\alpha\in\Lambda^mV^*$, we have the standard identities
    \[
        *\,\varepsilon_i\,*^{-1}\alpha=(-1)^{n-m}\iota_i\alpha, \qquad *\,\iota_i\,*^{-1}\alpha=(-1)^{n-m-1}\varepsilon_i\alpha.
    \]
    We first compute $*\,(\mathrm{Ric}(E)^{\sharp})^{[p]}\,*^{-1}=\sum_{i,j}\mathrm{Ric}_i^j\,(*\varepsilon_i*^{-1})(*\iota_j*^{-1})$ on $\Lambda^{n-p}$. Note the order of composition: first $*\iota_j*^{-1}$ acts, mapping $\Lambda^{n-p}$ to $\Lambda^{n-p+1}$, in which case $*\iota_j*^{-1}=(-1)^{p-1}\varepsilon_j$; then $*\varepsilon_i*^{-1}$ acts on $\Lambda^{n-p+1}$, in which case $*\varepsilon_i*^{-1}=(-1)^{n-(n-p+1)}\iota_i=(-1)^{p-1}\iota_i$. Thus the total coefficient is $(-1)^{(p-1)+(p-1)}=1$. Therefore
    \begin{equation}\label{eq:RicHodgeDual}
        *\,(\mathrm{Ric}(E)^{\sharp})^{[p]}\,*^{-1} = \sum_{i,j}\mathrm{Ric}_i^j\,\iota_i\varepsilon_j = \operatorname{Scal}(E)\mathrm{Id}_{\Lambda^{n-p}V^*} - (\mathrm{Ric}(E)^{\sharp})^{[n-p]},
    \end{equation}
    where the last equality uses $\iota_i\varepsilon_j=\delta_{ij}\mathrm{Id}-\varepsilon_j\iota_i$.

    We next compute $*\,\mathcal R_E^{(p)}\,*^{-1}=\sum_{i<j,\ a<b}E_{ijab} (*\varepsilon_i*^{-1})(*\varepsilon_j*^{-1})(*\iota_b*^{-1})(*\iota_a*^{-1})$. The total coefficient is $(-1)^{(p-1)+(p-2)+(p-2)+(p-1)}=1$. Therefore $*\,\mathcal R_E^{(p)}\,*^{-1} = \sum_{i<j,\ a<b}E_{ijab}\,\iota_i\iota_j\varepsilon_b\varepsilon_a$.

    We claim that
    \begin{equation}\label{eq:CurvHodgeDual}
        \mathcal R_E^{(n-p)}-*\,\mathcal R_E^{(p)}\,*^{-1} = (\mathrm{Ric}(E)^{\sharp})^{[n-p]} - \frac12\operatorname{Scal}(E)\mathrm{Id}_{\Lambda^{n-p}V^*}.
    \end{equation}
    The left is
    \begin{align*}
        &\quad\sum_{i<j,\ a<b}E_{ijab} \bigl(\varepsilon_a\varepsilon_b\iota_j\iota_i-\iota_i\iota_j\varepsilon_b\varepsilon_a\bigr) = \frac14\sum_{i,j,a,b}E_{ijab} \bigl(\varepsilon_a\varepsilon_b\iota_j\iota_i-\iota_j\iota_i\varepsilon_a\varepsilon_b\bigr)\\
        &= \frac14\sum_{i,j,a,b}E_{ijab} \bigl(\delta_{ai}\varepsilon_b\iota_j - \delta_{bi}\varepsilon_a\iota_j - \delta_{bj}\iota_i\varepsilon_a + \delta_{aj}\iota_i\varepsilon_b\bigr) = \frac12\sum_{i,j}\mathrm{Ric}_{ij}\bigl(2\varepsilon_j\iota_i-\delta_{ij}\mathrm{Id}\bigr),
    \end{align*}
    where the last two equalities use $\iota_i\varepsilon_j+\varepsilon_j\iota_i=\delta_{ij}\mathrm{Id}$. This proves the claim. Combining \eqref{eq:WeitCurv1}, \eqref{eq:RicHodgeDual} and \eqref{eq:CurvHodgeDual} yields $*\,q_p(E)\,*^{-1}=q_{n-p}(E)$.
\end{proof}

For the tensor $I_{ijkl}=\delta_{ik}\delta_{jl}-\delta_{il}\delta_{jk}$ which corresponds to constant curvature one, we have $\mathrm{Ric}(I)^{\sharp} = (n-1){\rm Id}_V$ and $\mathcal R_{I}={\rm Id}_{\Lambda^2V^*}$. Moreover,
\[
    \mathcal R_I^{(p)}=\frac{p(p-1)}2 {\rm Id}_{\Lambda^pV^*}, \qquad \left(\mathrm{Ric}(I)^\sharp\right)^{[p]}=p(n-1){\rm Id}_{\Lambda^pV^*}, \qquad q_p(I)=p(n-p){\rm Id}_{\Lambda^pV^*}.
\]
Indeed, by Definition \ref{def:ExtCurv},
\[
    \mathcal R_I^{(p)} = \sum_{i<j}\mathcal R_I(e^{ij})\wedge \iota_{j}\iota_{i} = \sum_{i<j} e^{ij}\wedge \iota_{j}\iota_{i} = \sum_{i<j}\varepsilon_i\varepsilon_j\iota_{j}\iota_{i} = \frac{p(p-1)}2 {\rm Id}_{\Lambda^pV^*}.
\]
Also, $(\mathrm{Ric}(I)^\sharp)^{[p]} =((n-1){\rm Id}_V)^{[p]} = p(n-1){\rm Id}_{\Lambda^pV^*}$. Then \eqref{eq:WeitCurv1} gives the formula for $q_p(I)$.

\subsection{The skew-endomorphism model: contractions and its full operator}

\begin{definition}\label{def:J}
	For a skew-adjoint endomorphism $A\in \mathrm{End}(V)$, define $E_A$ by
	\begin{equation}\label{eq:EJ}
        E_A(X,Y,Z,W) = \left\langle AX, Z \right\rangle\left\langle AY,W \right\rangle - \left\langle AX, W \right\rangle\left\langle AY,Z \right\rangle+2\left\langle AX, Y \right\rangle\left\langle AZ, W \right\rangle.
	\end{equation}
\end{definition}

The notation $A$ does not require $A$ to be a complex structure. Put $\Omega_A(X,Y)=\langle AX,Y\rangle$.

\begin{lemma}\label{lem:model-operator}
    The tensor $E_A$ is an algebraic curvature tensor. For orthonormal $X,Y$ it satisfies
    \begin{equation}\label{eq:model-contractions}
        K_{E_A}(X,Y)=3\Omega_A(X,Y)^2, \qquad \operatorname{Ric}(E_A)^\sharp=-3A^2, \qquad \operatorname{Scal}(E_A)=-3\operatorname{tr}(A^2).
    \end{equation}
    Let $D_A=A^{[p]}$, $L_A\alpha=\Omega_A\wedge\alpha$ and $\Lambda_A=L_A^*$. Then
    \begin{equation}\label{eq:model-operator}
        q_p(E_A)=-2(A^2)^{[p]}-D_A^2-4L_A\Lambda_A.
    \end{equation}
\end{lemma}

\begin{proof}
    The claim that $E_A$ is an algebraic curvature tensor and \eqref{eq:model-contractions} follow from direct computations. For \eqref{eq:model-operator} put $\Omega_{ia}=\langle Ae_i,e_a\rangle$. Then $D_A=\sum_{i,a}A_i^a\varepsilon_i\iota_a=\sum_{i,a}\Omega_{ia}\varepsilon_i\iota_a$. Using $\iota_a\varepsilon_j=\delta_{aj}\operatorname{Id}-\varepsilon_j\iota_a$, we have
    \[
        D_A^2 = \sum_{i,a,b}\Omega_{ia}\Omega_{ab}\varepsilon_i\iota_b-\sum_{i,j,a,b}\Omega_{ia}\Omega_{jb}\varepsilon_i\varepsilon_j\iota_a\iota_b = (A^2)^{[p]}+\sum_{i,j,a,b}\Omega_{ia}\Omega_{jb}\varepsilon_i\varepsilon_j\iota_b\iota_a.
    \]
    The first two products in $E_A$ form the pair-symmetric tensor $E^0_{ijab}:=\Omega_{ia}\Omega_{jb}-\Omega_{ib}\Omega_{ja}$.  Using the same pairwise operator formula, which is meaningful before imposing Bianchi, gives
    \[
        \mathcal R_{E^0}^{(p)} = \frac14\sum_{i,j,a,b}E^0_{ijab}\varepsilon_i\varepsilon_j\iota_b\iota_a = \frac12\sum_{i,j,a,b}\Omega_{ia}\Omega_{jb}\varepsilon_i\varepsilon_j\iota_b\iota_a = \frac12\bigl(D_A^2-(A^2)^{[p]}\bigr).
    \]
    The last term $2\Omega_A\otimes\Omega_A$ induces $\mathcal R_{2\Omega_A\otimes\Omega_A}^{(p)} = 2L_A\Lambda_A$ since $L_A=\sum_{i<j}(\Omega_A)_{ij}\varepsilon_i\varepsilon_j$ and $\Lambda_A=\sum_{a<b}(\Omega_A)_{ab}\iota_b\iota_a$. Consequently $\mathcal R_{E_A}^{(p)} = \frac12\bigl(D_A^2-(A^2)^{[p]}\bigr)+2L_A\Lambda_A$. Combine this identity with $\operatorname{Ric}(E_A)^\sharp=-3A^2$ and Definition \ref{def:WeitCurv} to obtain \eqref{eq:model-operator}.
\end{proof}

In particular, suppose that $W\subset V$ is a $J$-invariant subspace and that $J^{2}=-P_{W}$, where $P_{W}$ is the orthogonal projection onto $W$. Equivalently, $J|_W$ is an orthogonal complex structure and $J|_{W^\perp}=0$. The tensor $E_J$ is the zero extension of its restriction to $W$. In this case
\[
	\operatorname{Ric}(E_{J})(X, Y) = 3 \left\langle P_{W}X, Y \right\rangle,\quad {\rm Scal}(E_{J}) = 3 \dim W.
\]

The following lemma shows that, when $J$ is an orthogonal complex structure and $J|_{W^\perp}=0$, the power $\Omega^r$ of the K\"ahler form $\Omega$ is an eigenform of $q_{2r}(E_J)$, with eigenvalue $-4r(m-r)$.

\begin{lemma}\label{lem:model-powers}
    Suppose $W$ has dimension $2m$, $J|_W$ is an orthogonal complex structure, and $J=0$ on $W^\perp$. Set $\Omega=\Omega_J=\sum_{a=1}^m e^{2a-1,2a}$. If $0\le r\le m$, then
    \[
        |\Omega^r|^2=(r!)^2\binom mr, \qquad q_{2r}(E_J)\Omega^r=-4r(m-r)\Omega^r.
    \]
    If $\zeta$ is a unit covector in $W^\perp$, the same eigenvalue holds for $\Omega^r\wedge\zeta$ in degree $2r+1$, i.e.,
    \[
        q_{2r+1}(E_J)(\Omega^r\wedge\zeta)=-4r(m-r)(\Omega^r\wedge\zeta).
    \]
\end{lemma}

\begin{proof}
    Writing $\theta_a=e^{2a-1,2a}$ gives $\Omega^r=r!\sum_{|I|=r}\theta_I$. Since $\theta_I$ are orthonormal, we have $|\Omega^r|^2=(r!)^2\binom mr$. Since $J^{[1]}e^{2a-1}=-e^{2a}$ and $J^{[1]}e^{2a}=e^{2a-1}$, the Leibniz rule gives $J^{[2]}\theta_a=0$ and $J^{[2r]}\Omega^r=0$. Let $L:=\Omega\wedge=\sum_a\varepsilon_{2a-1}\varepsilon_{2a}$ and $\Lambda:=L^*=\sum_a\iota_{2a}\iota_{2a-1}$. Then $\Lambda\Omega^r = r!\sum_{|I|=r}\sum_{a\in I}\theta_{I\setminus\{a\}}$. Every subset $J$ of size $r-1$ is contained in $m-(r-1)$ subsets $I$ of size $r$, thus
    \[
        \Lambda\Omega^r = r!(m-r+1)\sum_{|J|=r-1}\theta_J = r(m-r+1)\Omega^{r-1}, \qquad L\Lambda\Omega^r=r(m-r+1)\Omega^r.
    \]
    On $\Lambda^{2r}W^*$, we have $-2(J^2)^{[2r]}=-2(-\operatorname{Id})^{[2r]}=4r\operatorname{Id}$. By \eqref{eq:model-operator} and $J^{[2r]}\Omega^r=0$, we have $q_{2r}(E_J)\Omega^r=-4r(m-r)\Omega^r$. Finally, since $B^{[p]}$ is a derivation for any $B\in{\rm End}(V)$,
    \[
        (J^2)^{[2r+1]}(\Omega^r\wedge\zeta)=((J^2)^{[2r]}\Omega^r)\wedge\zeta, \qquad (J^{[2r+1]})^2(\Omega^r\wedge\zeta)=((J^{[2r]})^2\Omega^r)\wedge\zeta.
    \]
    Since $\iota_i$ is an antiderivation, we have $\Lambda(\Omega^r\wedge\zeta)=(\Lambda\Omega^r)\wedge\zeta$. Therefore, by equation \eqref{eq:model-operator}, $q_{2r+1}(E_J)(\Omega^r\wedge\zeta)=(q_{2r}(E_J)\Omega^r)\wedge\zeta=-4r(m-r)(\Omega^r\wedge\zeta)$.
\end{proof}

\subsection{Energy notation}
\label{sec:energy}

For an arbitrary $p$-form $\omega$, we use the homogeneous notation
\begin{equation}\label{eq:energy-notation}
    \begin{split}
        N&=|\omega|^2,\qquad S=\operatorname{Scal}(E),\\
        Q_p(E,\omega)&=\langle q_p(E)\omega,\omega\rangle,\qquad
        F_p(E,\omega)=\langle(\operatorname{Ric}(E)^\sharp)^{[p]}\omega,\omega\rangle.
    \end{split}
\end{equation}
Write $R=\operatorname{Ric}(E)^\sharp$. Then the contraction Gram matrix $T_{ij}=\langle\iota_i\omega,\iota_j\omega\rangle$ satisfies
\[
    F_p=\operatorname{tr}(RT), \qquad \operatorname{tr}T=pN, \qquad 0\le T\le N\operatorname{Id}.
\]
In particular $0\le F_p\le SN$ when $\sec(E)\ge0$. Indeed, since $\varepsilon_i^*=\iota_i$ and $R^{[p]}=\sum_{i,j}\operatorname{Ric}(E)_i^j\,\varepsilon_i\iota_j$, we have $F_p=\sum_{i,j}\operatorname{Ric}(E)_{ij}T_{ij}$ and
\[
    \operatorname{tr}T=\sum_i T_{ii}=\sum_i|\iota_i\omega|^2=\sum_i\langle\omega,\varepsilon_i\iota_i\omega\rangle=\langle\omega,p\omega\rangle=pN,
\]
where we used the identity $\sum_i\varepsilon_i\iota_i\omega=p\omega$. Recall that for each $i$, the elementary decomposition
\[
    \omega=\varepsilon_i\iota_i\omega+\omega_0
\]
is orthogonal and has $\iota_i\omega_0=0$ since $\langle \varepsilon_i\iota_i\omega, \omega_0\rangle = \langle \varepsilon_i\iota_i\omega, \omega-\varepsilon_i\iota_i\omega\rangle = |\iota_i\omega|^2 - |\varepsilon_i\iota_i\omega|^2 =0$ and $\iota_i\omega_0 = \iota_i\omega-\iota_i\varepsilon_i\iota_i\omega = \iota_i\omega-(\operatorname{Id}-\varepsilon_i\iota_i)\iota_i\omega =0$. The orthogonality of the decomposition gives $|\iota_i\omega|^2=|\varepsilon_i\iota_i\omega|^2\le|\omega|^2=N$. Thus, for any vector $v=\sum v_i e_i$, we have $\langle T v,v\rangle=|\iota_v\omega|^2\le |v|^2 N$, and therefore $0\le T\le N\operatorname{Id}$. That is, $T$ is a positive semidefinite matrix, and all its eigenvalues are at most $N$. Thus, when $\sec(E)\ge0$, one has $\operatorname{Ric}(E)\ge0$ and thus $0\le F_p\le SN$.

From Lemma \ref{lem:HodgeDual} and \eqref{eq:RicHodgeDual}, one obtains the algebraic identities
\begin{equation}\label{eq:hodge-energy}
    Q_{n-p}(E,*\omega)=Q_p(E,\omega), \qquad F_{n-p}(E,*\omega)=SN-F_p(E,\omega),
\end{equation}
since the Hodge star is an isometry: $\langle*\omega,*\eta\rangle=\langle\omega,\eta\rangle$.

\section{Four-plane estimates and induction in the degree}
\label{alg:section}

Throughout this section we use $Q_p,F_p,S,N$ as in \eqref{eq:energy-notation}. When the tensor and the degree are clear, these quantities are denoted by $Q$ and $F$.

\subsection{The four-vector inequality and its Ricci-weighted version}

We use the Ni--Wilking four-frame estimate \cite{NiWilking} in the form used by Ge \cite{Ge}. Interior multiplication reduces the mixed terms in each degree to this four-vector calculation.

\begin{lemma}[The NW inequality]\label{alg:NW}\label{alg:nw}
    Assume $\sec(E)\geq0$. Let $e_1,e_2,e_3,e_4$ be orthonormal. Set
    \[
        k=E_{1212},\quad l=E_{3434},\quad r=E_{1234}, \quad A=E_{1313}+E_{1414}+E_{2323}+E_{2424}.
    \]
    Then
    \begin{equation}\label{alg:NW-equation}
        6|r|\leq A+4\sqrt{kl}.
    \end{equation}
    Also, for arbitrary $u,v,w\in V$, the following Jacobi--Cauchy--Schwarz inequality holds:
    \begin{equation}\label{alg:Jacobi-CS}
        |E(u,v,u,w)|^2\leq E(u,v,u,v)E(u,w,u,w).
    \end{equation}
\end{lemma}

\begin{proof}
    The inequality \eqref{alg:NW-equation} is \cite[Lemma~3.1]{Ge}. For \eqref{alg:Jacobi-CS}, since
    \[
        E(u,v,u,v)+2tE(u,v,u,w)+t^2E(u,w,u,w)=E(u,v+tw,u,v+tw)\ge0,
    \]
    its discriminant is nonpositive; if the quadratic coefficient vanishes, nonnegativity forces the linear coefficient to vanish as well. These two cases give the stated estimate.
\end{proof}

\begin{proposition}[A sharp affine family in degree two]
\label{alg:two-affine}\label{alg:q2-family}
    Let $d\geq2$, $m=\lfloor d/2\rfloor$, and let $E$ have nonnegative sectional curvature on $\mathbb R^d$. For
    \[
        -\frac{2(m-1)}3\leq a\leq\frac23, \qquad c_m(a)=\frac{2(m-1)}{3m}+\frac{a}m,
    \]
    one has
    \begin{equation}\label{alg:two-affine-equation}
        q_2(E)-a\bigl(\operatorname{Ric}(E)^\sharp\bigr)^{[2]} \geq -c_m(a)\operatorname{Scal}(E)\operatorname{Id}_{\Lambda^2}.
    \end{equation}
    For every permitted $a$, the coefficient $c_m(a)$ is sharp.
\end{proposition}

\begin{proof}
    Let $\alpha\in\Lambda^2(\mathbb R^d)^*$ be a unit two-form. The orthogonal normal form of a real skew-adjoint endomorphism shows that $\alpha$ can be written as
    \[
        \alpha=\sum_{r=1}^m\lambda_r\theta_r, \qquad \theta_r=e^{2r-1}\wedge e^{2r}, \qquad \sum_{r=1}^m\lambda_r^2=1.
    \]
    Let $P_r=\operatorname{span}\{e_{2r-1},e_{2r}\}$. Write
    \[
        k_r=E_{2r-1,2r,2r-1,2r}, \qquad r_{rs}=E_{2r-1,2r,2s-1,2s}, \qquad A_{rs}=\sum_{i,j=0}^1E_{2r-1+i,\,2s-1+j,\,2r-1+i,\,2s-1+j}.
    \]
    If $d=2m+1$, write $A_{r0}=E_{2r-1,d,2r-1,d}+E_{2r,d,2r,d}$. In even dimension put $A_{r0}=0$. We claim
    \begin{align}
        F_2(E,\alpha) &= 2\sum_r k_r\lambda_r^2+\sum_{r<s}(\lambda_r^2+\lambda_s^2)A_{rs}+\sum_r\lambda_r^2 A_{r0},\label{alg:two-Ricci-expansion}\\
        \langle\mathcal R_E\alpha,\alpha\rangle &= \sum_r k_r\lambda_r^2+2\sum_{r<s}\lambda_r\lambda_s r_{rs}.\label{alg:two-curvature-expansion}
    \end{align}
    Indeed, $\langle\mathcal R_E\alpha,\alpha\rangle=\sum_{r,s}\lambda_r\lambda_sE_{2r-1,2r,2s-1,2s}$. Then splitting the summation $\sum_{r,s}$ into the three sums $\sum_{r=s}+\sum_{r<s}+\sum_{r>s}$ yields \eqref{alg:two-curvature-expansion}. Fix a plane $P_r = \operatorname{span}\{e_{2r-1}, e_{2r}\}$. We have $(R^{[2]}\alpha)_{2r-1,2r} = \operatorname{Ric}_{2r-1}^i \alpha_{i,2r} + \operatorname{Ric}_{2r}^i \alpha_{2r-1,i}$. Only the terms with $i\in\{2r-1,2r\}$ are nonzero since $\alpha=\sum_r\lambda_r\theta_r$. Thus
    \[
        (R^{[2]}\alpha)_{2r-1,2r} = (\operatorname{Ric}_{2r-1}^{2r-1} + \operatorname{Ric}_{2r}^{2r}) \lambda_r, \qquad F_2(E,\alpha) = \langle R^{[2]}\alpha,\alpha\rangle = \sum_r (\operatorname{Ric}_{2r-1}^{2r-1} + \operatorname{Ric}_{2r}^{2r}) \lambda_r^2,
    \]
    where
    \[
        \operatorname{Ric}_{2r-1}^{2r-1} + \operatorname{Ric}_{2r}^{2r} = \sum_{j=1}^d E_{2r-1,j,2r-1,j} + \sum_{j=1}^d E_{2r,j,2r,j} = 2k_r + \sum_{s \neq r} A_{rs} + A_{r0}.
    \]
    Using the symmetry of $A_{rs}$ proves \eqref{alg:two-Ricci-expansion}. Since $q_2(E)=R^{[2]}-2\mathcal R_E$, we have
    \begin{align*}
        Q_2-aF_2 &= (1-a)F_2-2\langle\mathcal R_E\alpha,\alpha\rangle\\
        &= -2a\sum_r\lambda_r^2 k_r+\sum_{r<s}(1-a)(\lambda_r^2+\lambda_s^2)A_{rs}-4\sum_{r<s}\lambda_r\lambda_s r_{rs}\\
        &\quad +(1-a)\sum_r\lambda_r^2A_{r0}.
    \end{align*}
    Applying Lemma \ref{alg:NW} to each pair of planes yields $6|r_{rs}|\le A_{rs}+4\sqrt{k_r k_s}$. Then
    \[
        -4\lambda_r\lambda_s r_{rs}\ge -\frac23|\lambda_r\lambda_s|A_{rs}-\frac83|\lambda_r\lambda_s|\sqrt{k_r k_s}.
    \]
    Substituting this and putting $x_r:=|\lambda_r|\sqrt{k_r}$ gives
    \begin{align}
        Q_2-aF_2 &\geq -2a\sum_r x_r^2-\frac83\sum_{r<s}x_rx_s \notag \\
        &\quad +\sum_{r<s}\left[(1-a)(\lambda_r^2+\lambda_s^2)-\frac23|\lambda_r\lambda_s|\right]A_{rs}+(1-a)\sum_r\lambda_r^2A_{r0}. \label{alg:estimation-of-q2}
    \end{align}
    The bracket is $\ge(2/3-a)(\lambda_r^2+\lambda_s^2)\ge0$. Also $(1-a)\sum\lambda_r^2A_{r0}\ge0$ since $1-a>0$. Hence
    \[
        Q_2-aF_2 \geq -2a\sum_r x_r^2-\frac83\sum_{r<s}x_rx_s.
    \]
    Using the identity $2\sum_{r<s}x_rx_s=(\sum_r x_r)^2-\sum_r x_r^2$ and Cauchy--Schwarz $m\sum_r x_r^2\ge (\sum_r x_r)^2$, we obtain $Q_2-aF_2\ge-2c_m(a)(\sum_r x_r)^2$. Here the direction of the inequality uses $a\leq2/3$. The second Cauchy--Schwarz inequality is
    \[
        \Bigl(\sum_r|\lambda_r|\sqrt{k_r}\Bigr)^2 \leq \Bigl(\sum_r\lambda_r^2\Bigr)\Bigl(\sum_r k_r\Bigr) = \sum_r k_r\leq\frac12\operatorname{Scal}(E).
    \]
    Substituting this into $Q_2-aF_2\ge-2c_m(a)(\sum_r x_r)^2$ proves \eqref{alg:two-affine-equation}.

    For sharpness take $E_J$ on a $2m$-plane $W$, extended by zero if $d$ is odd, and $\alpha=\Omega_J/\sqrt m$, where $\Omega_J=\sum_r\theta_r$. Then $|\alpha|^2=1$ and $R|_W=3\operatorname{Id}_W$. Thus $R^{[2]}|_{\Lambda^2W^*}=6\operatorname{Id}_{\Lambda^2W^*}$ and $F_2(E_J,\alpha)=\langle6\alpha,\alpha\rangle=6$. For $r\ne s$, we have $k_r=3$ and $r_{rs}=2$. Then by \eqref{alg:two-curvature-expansion}, $\langle\mathcal R_{E_J}\alpha,\alpha\rangle=2m+1$, so $Q_2(E_J,\alpha)=F_2(E_J,\alpha)-2\langle\mathcal R_{E_J}\alpha,\alpha\rangle=-4(m-1)$. Consequently, $Q_2-aF_2=-c_m(a)\operatorname{Scal}(E_J)$ since $\operatorname{Scal}(E_J)=6m$.
\end{proof}

\subsection{Exact contraction identities}

For $\omega\in\Lambda^pV^*$, let $\alpha_i=\iota_i\omega$ and $T_{ij}=\langle\alpha_i,\alpha_j\rangle$. Then $T$ is symmetric and positive semidefinite. As shown in Section \ref{sec:energy}, $\operatorname{tr}T=pN$ and $F=\operatorname{tr}(RT)$.

\begin{lemma}[Contraction and restriction]\label{alg:contractions}
    Suppose $p\geq2$. For an orthonormal basis, put
    \[
        H_i=e_i^\perp,\quad E_i=E|_{H_i}, \quad J_i(u,v)=E(e_i,u,e_i,v), \quad D=\sum_i\langle J_i^{[p-1]}\alpha_i,\alpha_i\rangle.
    \]
    We identify each symmetric form $J_i$ with its self-adjoint endomorphism $J_i^\sharp$. Then
    \begin{align}
        \sum_i Q_{p-1}(E,\alpha_i)&=(p-2)Q+F, \label{alg:full-contraction}\\
        \sum_i Q_{p-1}(E_i,\alpha_i)&=(p-2)Q+F-D, \label{alg:hyperplane-Q}\\
        \sum_i F_{p-1}(E_i,\alpha_i)&=(p-1)F-D. \label{alg:hyperplane-F}
    \end{align}
    In the first identity $q_{p-1}(E)$ acts on the full space $V$; in the other two it acts on $H_i$. These three identities do not require $E$ to have nonnegative sectional curvature. If in addition $E$ has nonnegative sectional curvature, then $D\ge0$.
\end{lemma}

\begin{proof}
    For $B=\varepsilon_{i_1}\cdots\varepsilon_{i_r}\iota_{j_1}\cdots\iota_{j_r}$, we have $\sum_k\varepsilon_k B\iota_k=(p-r)B$ on $\Lambda^p$. Indeed, commuting $\varepsilon_k$ past the exterior factors and $\iota_k$ past the interior factors gives two signs $(-1)^r$, leaving the degree operator $\sum_k\varepsilon_k\iota_k=(p-r)\operatorname{Id}|_{\Lambda^{p-r}}$. In particular $\sum_k\varepsilon_k(\varepsilon_i\iota_j)\iota_k = (p-1)\varepsilon_i\iota_j$ and $\sum_k\varepsilon_k(\varepsilon_i\varepsilon_j\iota_b\iota_a)\iota_k = (p-2)\varepsilon_i\varepsilon_j\iota_b\iota_a$ on $\Lambda^p$. These two identities yields
    \begin{equation}\label{alg:contraction-identities}
        \sum_i\langle R^{[p-1]}\alpha_i,\alpha_i\rangle=(p-1)F, \qquad \sum_i\langle\mathcal R_E^{(p-1)}\alpha_i,\alpha_i\rangle=(p-2)\langle\mathcal R_E^{(p)}\omega,\omega\rangle.
    \end{equation}
    This gives \eqref{alg:full-contraction} since $Q=Q_p(E,\omega)=F-2\langle\mathcal R_E^{(p)}\omega,\omega\rangle$.

    The form $\alpha_i$ has no $e^i$ factor and is therefore a $(p-1)$-form on $H_i$. For $u,v\in H_i$,
    \[
        \operatorname{Ric}(E_i)(u,v)=\operatorname{Ric}(E)(u,v)-J_i(u,v).
    \]
    This yields $F_{p-1}(E_i,\alpha_i) = \langle(\operatorname{Ric}(E)^\sharp)^{[p-1]}\alpha_i,\alpha_i\rangle-\langle J_i^{[p-1]}\alpha_i,\alpha_i\rangle$. Summing over $i$ and using the first identity in \eqref{alg:contraction-identities} proves \eqref{alg:hyperplane-F}.

    Since $\iota_i\alpha_i=\iota_i\iota_i\omega=0$, we have $\langle\varepsilon_i\beta,\alpha_i\rangle=0$ for any form $\beta$. This gives
    \[
        \langle\mathcal R_E^{(p-1)}\alpha_i,\alpha_i\rangle=\sum_{a<b,\ c<d}E_{abcd}\langle\varepsilon_a\varepsilon_b\iota_d\iota_c\alpha_i,\alpha_i\rangle = \sum_{\substack{a<b,\ c<d\\ a,b,c,d\neq i}}E_{abcd}\langle\varepsilon_a\varepsilon_b\iota_d\iota_c\alpha_i,\alpha_i\rangle = \langle\mathcal R_{E_i}^{(p-1)}\alpha_i,\alpha_i\rangle.
    \]
    Thus $Q_{p-1}(E_i,\alpha_i) = F_{p-1}(E_i,\alpha_i)-2\langle\mathcal R_E^{(p-1)}\alpha_i,\alpha_i\rangle$. By \eqref{alg:hyperplane-F} and the second identity in \eqref{alg:contraction-identities},
    \[
        \sum_i Q_{p-1}(E_i,\alpha_i) = \sum_i F_{p-1}(E_i,\alpha_i)-2\sum_i\langle\mathcal R_E^{(p-1)}\alpha_i,\alpha_i\rangle = (p-1)F-D-2(p-2)\langle\mathcal R_E^{(p)}\omega,\omega\rangle.
    \]
    This proves \eqref{alg:hyperplane-Q} since $\langle\mathcal R_E^{(p)}\omega,\omega\rangle=(F-Q)/2$.

    Finally, if $\sec(E)\ge0$, then $J_i\ge0$; hence $J_i^{[p-1]}\ge0$, and therefore $D\ge0$.
\end{proof}

For a fixed orthonormal basis, we have $\operatorname{Scal}(E_i)=S-2R_{ii}$, so
\begin{equation}\label{alg:basis-depend}
    \sum_i\operatorname{Scal}(E_i)|\alpha_i|^2 = \sum_i(S-2R_{ii})T_{ii} = pSN-2\sum_i R_{ii}T_{ii}.
\end{equation}
This expression depends on the chosen basis. Indeed, $F=\sum_i R_{ii}T_{ii}$ holds in an eigenbasis of $R$ or $T$ but may fail in a general frame.

\subsection{Spherical fourth moments and the averaged identities}

Let $d\mu=d\mathrm{vol}/\operatorname{Vol}(S^{n-1})$ be the normalized volume form on $S^{n-1}$. Let $dH$ be the Haar volume measure on $O(n)$, and set $dH_0=dH/\operatorname{Vol}(O(n))$. An overline denotes the average over ordered orthonormal frames:
\[
    \overline{\Phi}=\frac{1}{\operatorname{Vol}(O(n))}\int_{O(n)}\Phi(A)\,dH(A)=\int_{O(n)}\Phi(A)\,dH_0(A).
\]
Write $A=(e_1,\cdots,e_n)\in O(n)$. If $f$ is a function of a single vector variable, invariance gives
\[
    \overline{\sum_i f(e_i)}=\int_{O(n)}\sum_i f(e_i)\,dH_0=n\int_{S^{n-1}}f(v)\,d\mu(v).
\]
The second and fourth moments are
\begin{equation}\label{alg:sphere-moments}
    \int_{S^{n-1}} v_iv_j\,d\mu=\frac{\delta_{ij}}n, \qquad \int_{S^{n-1}} v_iv_jv_kv_l\,d\mu=\frac{\delta_{ij}\delta_{kl}+\delta_{ik}\delta_{jl}+\delta_{il}\delta_{jk}}{n(n+2)},
\end{equation}
where $v=(v_1,\cdots,v_n)\in S^{n-1}\subset\mathbb R^n$.

\begin{proposition}[Exact averaged formulas]\label{alg:averages}
    With the notation of Lemma \ref{alg:contractions}, let
    \[
        c=n-p,\qquad L=(p-2)n+2p-1,\qquad H=(p-1)n+p-4.
    \]
    Then
    \begin{align}
        \overline D&=\frac{(p+2)F-3Q}{n+2},\label{alg:average-D}\\
        \overline{\sum_i\operatorname{Scal}(E_i)|\alpha_i|^2}&=\frac{npSN-4F}{n+2},\label{alg:average-scalar}\\
        \overline{\sum_iQ_{p-1}(E_i,\alpha_i)}&=\frac{LQ+cF}{n+2},\label{alg:average-Q}\\
        \overline{\sum_iF_{p-1}(E_i,\alpha_i)}&=\frac{3Q+HF}{n+2}.\label{alg:average-F}
    \end{align}
    These are identities for every algebraic curvature tensor; sectional nonnegativity is needed only when the resulting quantities are used as nonnegative remainders.
\end{proposition}

\begin{proof}
    Put $U=\langle\mathcal R_E^{(p)}\omega,\omega\rangle$ and $\beta_{ia}=\iota_a\iota_i\omega$, so that $\beta_{ia}=-\beta_{ai}$. For $v=\sum_iv_ie_i$, the integrand defining $D$ is $\langle J_v^{[p-1]}\iota_v\omega,\iota_v\omega\rangle$. This equals
    \[
        \sum_{a,b}\langle(J_v)_{ab}\varepsilon_a\iota_b\iota_v\omega,\iota_v\omega\rangle = \sum_{a,b}E(v,e_a,v,e_b)\langle\iota_a\iota_v\omega,\iota_b\iota_v\omega\rangle = \sum_{i,j,k,l,a,b}E_{iajb}v_iv_jv_kv_l\langle\beta_{ka},\beta_{lb}\rangle.
    \]
    We evaluate each of the three Kronecker pairings in \eqref{alg:sphere-moments} separately. Pairing $i$ with $j$ and $k$ with $l$ gives $\sum_{k,a,b}R_{ab}\langle\iota_a\iota_k\omega,\iota_b\iota_k\omega\rangle$. By the first identity in \eqref{alg:contraction-identities},
    \[
        \sum_{k,a,b}R_{ab}\langle\iota_a\iota_k\omega,\iota_b\iota_k\omega\rangle = \sum_{k,a,b}R_{ab}\langle\varepsilon_b\iota_a\iota_k\omega,\iota_k\omega\rangle = \sum_{k}\langle R^{[p-1]}\iota_k\omega,\iota_k\omega\rangle = (p-1)F,
    \]
    Pairing $i$ with $k$ and $j$ with $l$ gives
    \[
        U_1 := \sum_{i,j,a,b}E_{iajb}\langle\beta_{ia},\beta_{jb}\rangle = 4\sum_{i<a,\ j<b}E_{iajb}\langle\beta_{ia},\beta_{jb}\rangle = 4\sum_{i<a,\ j<b}E_{iajb}\langle\varepsilon_j\varepsilon_b\iota_a\iota_i\omega,\omega\rangle = 4U.
    \]
    The remaining pairing gives $U_2 := \sum_{i,j,a,b}E_{iajb}\langle\beta_{ja},\beta_{ib}\rangle$. Interchanging $i$ and $j$ gives $U_2=\sum_{i,j,a,b}E_{jaib}\langle\beta_{ia},\beta_{jb}\rangle$. Alternatively, first interchange $i$ and $a$, use the skew symmetry of both $E$ and $\beta$, and then interchange $a$ and $j$. This gives $U_2=\sum_{i,j,a,b}E_{ijab}\langle\beta_{ia},\beta_{jb}\rangle$. Adding the last two expressions and using the first Bianchi identity $E_{jaib}+E_{ijab}=E_{iajb}$ shows that $2U_2=U_1$, so $U_2=2U$. The three contributions are thus $(p-1)F$, $4U$, and $2U$. Multiplying the spherical average by $n$,
    \[
        \overline D=\frac{(p-1)F+6U}{n+2}=\frac{(p+2)F-3Q}{n+2},
    \]
    since $2U=F-Q$. This proves \eqref{alg:average-D}.

    For the scalar term, $|\iota_v\omega|^2=T(v,v)$. Applying the same fourth moments to the product of the two symmetric quadratic forms $T(v,v)$ and $R(v,v)$ gives
    \[
        n\int|\iota_v\omega|^2R(v,v)\,d\mu(v) = \frac{\operatorname{tr}T\operatorname{tr}R+2\operatorname{tr}(TR)}{n+2} = \frac{pSN+2F}{n+2},
    \]
    where we used $\operatorname{tr}T=pN$, $\operatorname{tr}R=S$ and $\operatorname{tr}(TR)=F$. By \eqref{alg:basis-depend}, we have
    \[
        \overline{\sum_i\operatorname{Scal}(E_i)|\alpha_i|^2} = pSN-\frac{2pSN+4F}{n+2} = \frac{npSN-4F}{n+2}.
    \]
    Finally, by \eqref{alg:hyperplane-Q} and \eqref{alg:hyperplane-F}, we have
    \[
        \overline{\sum_i Q_{p-1}(E_i,\alpha_i)} = (p-2)Q+F-\overline D, \qquad \overline{\sum_i F_{p-1}(E_i,\alpha_i)} = (p-1)F-\overline D.
    \]
    Substituting \eqref{alg:average-D} gives the last two formulas.
\end{proof}

\subsection{The complete affine induction}

\begin{proposition}[Induction from degree $p-1$ to degree $p$]\label{alg:affine-recurrence}
    Assume that in dimension $n-1$, for every algebraic curvature tensor $E'$ with nonnegative sectional curvature and every $(p-1)$-form $\alpha$, it holds
    \[
        uQ_{p-1}(E',\alpha)+vF_{p-1}(E',\alpha)+w\operatorname{Scal}(E')|\alpha|^2\geq0.
    \]
    Then in dimension $n$, for every algebraic curvature tensor $E$ with nonnegative sectional curvature and every $p$-form $\omega$,
    \[
        u'Q_p+v'F_p+w'S|\omega|^2\geq0,
    \]
    where
    \begin{equation}\label{alg:affine-matrix}
        \begin{pmatrix}u'\\v'\\w'\end{pmatrix}
        =\begin{pmatrix}L&3&0\\c&H&-4\\0&0&np\end{pmatrix}
        \begin{pmatrix}u\\v\\w\end{pmatrix},
        \qquad c=n-p.
    \end{equation}
    Here $L=(p-2)n+2p-1$, $H=(p-1)n+p-4$.
\end{proposition}

\begin{proof}
    The restriction $E_i$ has nonnegative sectional curvature. Apply the assumed homogeneous inequality to each $(E_i,\alpha_i)$:
    \[
        uQ_{p-1}(E_i,\alpha_i)+vF_{p-1}(E_i,\alpha_i)+w\operatorname{Scal}(E_i)|\alpha_i|^2\ge0.
    \]
    Summing over all $i$ and taking the Haar average over all orthonormal frames, we obtain
    \[
        u\,\overline{\sum_iQ_{p-1}(E_i,\alpha_i)}+v\,\overline{\sum_iF_{p-1}(E_i,\alpha_i)}+w\,\overline{\sum_i\operatorname{Scal}(E_i)|\alpha_i|^2}\ge0.
    \]
    By Proposition \ref{alg:averages}, the result is $[u(LQ+cF)+v(3Q+HF)+w(npSN-4F)]/(n+2)\ge0$. Multiplying by $n+2$ proves the formula. Since no terms are discarded in this proof (in particular, the Jacobi term $D$ is computed exactly and retained), the derivation from the induction hypothesis to the conclusion is an exact chain of equalities, rather than a rough inequality estimate obtained by relaxations such as discarding nonnegative terms.
\end{proof}

If $w=0$, then the scalar-curvature term is absent, and the update formula for $w$ gives $w'=npw=0$. Hence the scalar-curvature term never reappears, and the $3\times 3$ matrix in \eqref{alg:affine-matrix} reduces to the $2\times 2$ matrix
\[
    M=
    \begin{pmatrix}
        L&3\\
        c&H
    \end{pmatrix},
\]
acting on the coefficient vector $\binom{u}{v}$. This matrix has two explicit eigenvectors:
\begin{equation}\label{alg:degree-eigenvectors}
    \binom{1}{c/3}, \qquad \binom{1}{-1},
\end{equation}
with eigenvalues
\[
    \lambda_+(n,p)=(p-1)(n+1), \qquad \lambda_-(n,p)=(p-2)(n+2).
\]

When $(n-1,p-1)$ changes to $(n,p)$, the codimension $(n-1)-(p-1)=n-p=c$ stays fixed. Since the two eigenvectors of $M$ in \eqref{alg:degree-eigenvectors} depend only on $c$, they are identical at every step of the induction. Hence, at every step $M$ can be diagonalized in exactly the same eigenbasis. Let the initial coefficient vector be
\[
    \binom{u_0}{v_0}=A\binom{1}{c/3}+B\binom{1}{-1}.
\]
At the $j$-th step (from dimension $n_0+j-1$ and degree $p_0+j-1$ to dimension $n_0+j$ and degree $p_0+j$), the eigenvalues of the matrix $M_j$ are
\[
    \lambda_+^{(j)}=(p_0+j-1)(n_0+j+1), \qquad \lambda_-^{(j)}=(p_0+j-2)(n_0+j+2).
\]
After $k$ steps, the coefficient vector becomes
\[
    \binom{u_k}{v_k}=A\left(\prod_{j=1}^{k}\lambda_+^{(j)}\right)\binom{1}{c/3}+B\left(\prod_{j=1}^{k}\lambda_-^{(j)}\right)\binom{1}{-1}.
\]
Therefore, the entire induction can be carried out explicitly by tracking the products of the coefficients of the two eigenmodes separately, without performing matrix multiplications step by step. These elementary eigenvectors make the entire induction explicit.

\begin{theorem}[Ricci-weighted bounds in every degree]\label{alg:Ricci-all}
    Let $E$ be an algebraic curvature tensor with nonnegative sectional curvature. Let $p\geq2$ and $c=n-p\geq2$. Define
    \begin{equation}\label{alg:kappa-raw}
        \kappa^0_{n,p}=
        \begin{cases}
            \dfrac c3,&c\text{ even},\\[6pt]
            \displaystyle\frac{c(c+2)-3\rho}{3(c+2+\rho)}, \quad \rho=\dfrac{n+2}{(p-1)(c+4)},&c\text{ odd}.
        \end{cases}
    \end{equation}
    Equivalently, the odd branch is
    \begin{equation}\label{alg:kappa-odd-simplified}
        \kappa^0_{n,p}=\frac c3-\frac{n+2}{3\bigl((p-1)(c+3)+1\bigr)}.
    \end{equation}
    Then
    \begin{equation}\label{alg:Ricci-all-equation}
        q_p(E)+\kappa^0_{n,p}R^{[p]}\geq0.
    \end{equation}
\end{theorem}

\begin{proof}
    The proof proceeds by induction on $p$ with $c=n-p$ fixed, starting from the base case $p=2$, in which case $n=c+2$. Let $m=\lfloor n/2\rfloor=\lfloor(c+2)/2\rfloor$. At the lower endpoint $a=-2(m-1)/3$ of Proposition \ref{alg:two-affine}, \eqref{alg:two-affine-equation} becomes $q_2(E)+2(m-1)R^{[2]}/3\ge0$ since $c_m(a)=0$. This is equivalent to that $Q_2(E,\omega)+\kappa^0_{c+2,2}F_2(E,\omega)\ge0$ for every $2$-form $\omega$, where $\kappa^0_{c+2,2}=2(m-1)/3=2\lfloor c/2\rfloor/3$. Hence the initial vector is exactly $(1,\kappa^0_{c+2,2})$.

    If $c$ is even, the initial vector is $(1,c/3)$. Equation \eqref{alg:degree-eigenvectors} shows that every inductive step multiplies this vector by the positive number $(p-1)(n+1)$. Dividing by that number preserves the coefficient $c/3$. This proves the even branch for every $p$.

    If $c$ is odd, the initial coefficient is $(c-1)/3$, and
    \[
        \binom{1}{(c-1)/3}=\frac{c+2}{c+3}\binom{1}{c/3}+\frac1{c+3}\binom{1}{-1}.
    \]
    At the step with degree $j$, the ambient dimension is $j+c$. After the steps $j=3,\cdots,p$, the ratio of the two accumulated eigenvalues is
    \[
        \prod_{j=3}^p\frac{(j-2)(j+c+2)}{(j-1)(j+c+1)}=\left(\prod_{j=3}^p\frac{j-2}{j-1}\right)\left(\prod_{j=3}^p\frac{j+c+2}{j+c+1}\right)=\frac1{p-1}\frac{p+c+2}{c+4}=\rho.
    \]
    For $p=2$, both products are empty and equal one; the displayed formula for $\rho$ also equals one. Dividing the resulting vector by its positive first component gives $\frac{(c+2)c/3-\rho}{c+2+\rho}$.
\end{proof}

\section{Three-forms and the sharp six-dimensional estimate}
\label{sec:three-forms}

Throughout this section $E$ is an algebraic curvature tensor with nonnegative sectional curvature. We continue to use $R=\operatorname{Ric}(E)^{\sharp}$ and $S=\operatorname{Scal}(E)$. For a unit three-form $\omega$, the two quadratic quantities that we shall keep separate are
\begin{equation}\label{six:QF}
    Q=\langle q_3(E)\omega,\omega\rangle, \qquad F=\langle R^{[3]}\omega,\omega\rangle.
\end{equation}
Keeping $F$ in the argument, rather than estimating it immediately in terms of $S$, will provide stronger input for the induction in dimension and degree.

\subsection{A six-dimensional normal form with its Hodge phase}

\begin{lemma}\label{six:normal-form}
    Let $V$ be an oriented Euclidean space of dimension six. Every unit real three-form $\omega$ can be written, in a suitable oriented orthonormal coframe, as
    \begin{equation}\label{six:normal-form-equation}
        \omega=(\cos\theta+\sin\theta\,*)\omega_0, \qquad \omega_0=a e^{123}+b e^{145}+c e^{246}+d e^{356}, \qquad a^2+b^2+c^2+d^2=1.
    \end{equation}
    The Hodge phase $\cos\theta+\sin\theta\,*$ is an orthogonal transformation of $\Lambda^3V^*$ and preserves the quadratic form of $q_3(E)$. It is not, in general, induced by an orthogonal transformation of $V$.
\end{lemma}

\begin{proof}
    We give the representation-theoretic reduction explicitly, since the extra phase is relevant to the Ricci term below. Now we divide the argument into steps.

    \medskip
    \noindent
    \textbf{1. The real spin module and its complex structure.}
    Let $\Sigma$ be the irreducible eight-dimensional real module for $\mathrm{Cl}_{0,6}$, and choose an invariant Euclidean inner product $g_\Sigma$ on it. If $e_1,\cdots,e_6$ is an oriented orthonormal basis of $V$, write $c_i=c(e_i)$. Thus
    \begin{equation}\label{eq:clifford-relations}
        c_i^*=-c_i, \qquad c_ic_j+c_jc_i=-2\delta_{ij}\operatorname{Id},
    \end{equation}
    where the adjoint is taken with respect to $g_\Sigma$. In particular, $c_i^2=-\operatorname{Id}$.

    Put $\mathsf J=c_1c_2c_3c_4c_5c_6$. By \eqref{eq:clifford-relations}, $\mathsf J^2=-\operatorname{Id}$ and $\mathsf J^*=-\mathsf J$. Thus $\mathsf J$ is an orthogonal complex structure and turns $\Sigma$ into a complex Hermitian space of dimension four. We regard $\Sigma$ as a complex vector space by declaring multiplication by $i$ to be $\mathsf J$. Using the convention that a Hermitian product is complex-linear in its first variable, set $h(x,y)=g_\Sigma(x,y)-i g_\Sigma(\mathsf Jx,y)$. This is a positive-definite Hermitian product whose real part is $g_\Sigma$.

    Each threefold product $c_ic_jc_k$ is symmetric and anticommutes with $\mathsf J$, since
    \[
        (c_ic_jc_k)^*=c_k^*c_j^*c_i^*=-c_kc_jc_i=c_ic_jc_k, \qquad c_ic_jc_k\mathsf J=-c_ic_j\mathsf Jc_k=c_i\mathsf Jc_jc_k=-\mathsf Jc_ic_jc_k.
    \]
    The second identity says precisely that $c_ic_jc_k$ is complex-antilinear. It is also a self-adjoint real map.

    For a real three-form $\omega=\sum_{i<j<k}\omega_{ijk}e^{ijk}$, define
    \begin{equation}\label{eq:clifford-map}
        \mathsf c(\omega)=\sum_{i<j<k}\omega_{ijk}c_ic_jc_k.
    \end{equation}
    From the preceding paragraph, $\mathsf c(\omega)$ is a self-adjoint complex-antilinear endomorphism of $\Sigma$.

    \medskip
    \noindent
    \textbf{2. The trace identity.}
    We next prove
    \begin{equation}\label{six:clifford-isometry}
        \operatorname{tr}_{\mathbb R}\!\bigl(\mathsf c(\eta)\mathsf c(\xi)\bigr)=8\langle\eta,\xi\rangle \qquad (\eta,\xi\in\Lambda^3V^*).
    \end{equation}
    For a three-element set $I=\{i<j<k\}$, abbreviate $c_I=c_ic_jc_k$. If $I\ne K$, using $c_i^2=-\operatorname{Id}$ to cancel the generators common to $I$ and $K$ reduces $c_Ic_K$, up to sign, to the Clifford monomial whose indices form $I\triangle K$. Since two distinct three-element subsets meet in zero, one, or two elements, this is a monomial of degree six, four, or two, respectively.

    Every nonconstant even monomial $M$ of degree two, four, or six has trace zero. Indeed, choose a generator $c_l$ occurring in $M$, then $c_l M=-Mc_l$. Since $c_l$ is invertible, this gives $c_l M c_l^{-1}=-M$. Similar real endomorphisms have the same trace, and therefore
    \[
        \operatorname{tr}_{\mathbb R}(M)=\operatorname{tr}_{\mathbb R}(-M)=-\operatorname{tr}_{\mathbb R}(M), \qquad \operatorname{tr}_{\mathbb R}(M)=0.
    \]
    Consequently $\operatorname{tr}_{\mathbb R}(c_Ic_K)=0$ whenever $I\ne K$.

    On the other hand, $c_I^2=\operatorname{Id}$ by \eqref{eq:clifford-relations}, so $\operatorname{tr}_{\mathbb R}(c_I^2)=\dim_{\mathbb R}\Sigma=8$. Expanding \eqref{six:clifford-isometry} now gives
    \[
        \operatorname{tr}_{\mathbb R}\!\bigl(\mathsf c(\eta)\mathsf c(\xi)\bigr) = \sum_{I,K}\eta_I\xi_K\operatorname{tr}_{\mathbb R}(c_Ic_K) = 8\sum_I\eta_I\xi_I = 8\langle\eta,\xi\rangle,
    \]
    as claimed. In particular, $\operatorname{tr}_{\mathbb R}\!\bigl(\mathsf c(\eta)^2\bigr)=8|\eta|^2$. Thus if $\mathsf c(\eta)=0$, then $\eta=0$. Hence $\mathsf c$ is injective.

    \medskip
    \noindent
    \textbf{3. Symmetric antilinear maps and Takagi factorization.}
    Choose a unitary basis of $(\Sigma,\mathsf J,h)$ and identify it with $\mathbb C^4$. Every complex-antilinear map $A$ has a unique expression
    \[
        A(z)=M\overline z, \qquad M\in M_4(\mathbb C).
    \]
    In these coordinates the real inner product is $g_\Sigma(z,w)=\operatorname{Re}(z^{\mathsf T}\overline w)$. A direct calculation gives
    \[
        g_\Sigma(Az,w)=\operatorname{Re}\!\left(z^{\mathsf T}\overline M^{\mathsf T}w\right), \qquad g_\Sigma(z,Aw)=\operatorname{Re}\!\left(z^{\mathsf T}\overline M w\right).
    \]
    Thus $A$ is self-adjoint as a real map if and only if $M^{\mathsf T}=M$. (To obtain the ``only if'' implication, vary $z,w$ and also replace either one by $iz$ or $iw$.)

    The space of complex symmetric $4\times4$ matrices has $10$ independent complex entries and hence real dimension $20$. Since $\dim_{\mathbb R}\Lambda^3V^*=\binom63=20$, the injective real-linear map \eqref{eq:clifford-map} is an isomorphism from $\Lambda^3V^*$ onto the space of self-adjoint complex-antilinear endomorphisms of $\Sigma$.

    Let $A=\mathsf c(\omega)$ and let $M$ be its symmetric matrix. Takagi factorization says that there are a unitary matrix $U$ and a real diagonal matrix
    \[
        D=\operatorname{diag}(s_1,s_2,s_3,s_4), \qquad s_\alpha\geq0,
    \]
    such that $M=UDU^{\mathsf T}$. Equivalently, if $D$ also denotes the antilinear map $z\mapsto D\overline z$, then
    \begin{equation}\label{eq:takagi-conjugation}
        U^{-1}AU=D.
    \end{equation}
    Indeed, $(U^{-1}AU)(z) = U^*M\overline U\,\overline z = U^*(UDU^{\mathsf T})\overline U\,\overline z=D\overline z$.

    \medskip
    \noindent
    \textbf{4. Why the special-unitary part is a change of coframe.}
    The even Clifford algebra commutes with $\mathsf J$, so the spin representation
    \[
        \rho:\operatorname{Spin}(6)\longrightarrow U(\Sigma,\mathsf J,h)\cong U(4)
    \]
    is complex-linear and unitary. Its infinitesimal image is spanned by $\frac12c_ic_j$, $1\leq i<j\leq6$. Each $c_ic_j$ is skew-adjoint and complex-linear. It also has complex trace zero. To see the latter directly, choose $k$ distinct from $i,j$. The complex-linear invertible operator $c_ic_k$ anticommutes with $c_ic_j$; hence $c_ic_j$ is similar, over $\mathbb C$, to its negative, and its complex trace vanishes.

    The fifteen bivectors $c_ic_j$ are linearly independent (the Clifford representation of $\mathrm{Cl}_{0,6}\cong M_8(\mathbb R)$ is faithful). They therefore form a fifteen-dimensional subalgebra of $\mathfrak u(4)$ contained in $\mathfrak{su}(4)$. Since $\dim_{\mathbb R}\mathfrak{su}(4)=15$, the infinitesimal image is all of $\mathfrak{su}(4)$. The group $\operatorname{Spin}(6)$ is connected, so the determinant of its image is identically one; hence the image lies in $SU(4)$. A connected Lie subgroup with Lie algebra $\mathfrak{su}(4)$ is open in $SU(4)$, and therefore equals $SU(4)$ because $SU(4)$ is connected. Thus
    \begin{equation}\label{eq:spin-su4}
        \rho(\operatorname{Spin}(6))=SU(4).
    \end{equation}

    Clifford multiplication is equivariant for this action. More precisely, if $s\in\operatorname{Spin}(6)$ projects to $g\in SO(V)$ and $S=\rho(s)$, then, after identifying $V$ and $V^*$ by the metric,
    \begin{equation}\label{eq:equivariance}
        S^{-1}\mathsf c(\alpha)S=\mathsf c(g^{-1}\!\cdot\alpha) \qquad (\alpha\in\Lambda^3V^*).
    \end{equation}
    Thus conjugation by an element of $SU(4)$ is exactly the spinorial form of an oriented orthonormal change of coframe.

    \medskip
    \noindent
    \textbf{5. A fixed real diagonal section.}
    We now identify, in exterior-form notation, the four-dimensional space of real diagonal antilinear maps used in the Takagi factorization. Set
    \[
        C_1=c_1c_2c_3,\quad C_2=c_1c_4c_5,\quad C_3=c_2c_4c_6,\quad C_4=c_3c_5c_6.
    \]
    They are complex-antilinear (i.e., $C_r\mathsf J=-\mathsf J C_r$), self-adjoint, pairwise commuting (i.e., $C_rC_s=C_sC_r$) involutions (i.e., $C_r^2=\operatorname{Id}$). Since $C_1$ is an antilinear isometry and $C_1^2=\operatorname{Id}$, it behaves like complex conjugation. Use $C_1$ as a conjugation on the complex vector space $\Sigma$. Its fixed space $W=\{x\in\Sigma \mid C_1x=x\}$ is a real form of $\Sigma$: every $x\in\Sigma$ has a unique decomposition
    \begin{equation}\label{eq:real-form-decomposition}
        x=u+\mathsf Jv, \qquad u,v\in W.
    \end{equation}
    Indeed, one may take $u=(x+C_1x)/2$ and $v=-\mathsf J(x-C_1x)/2$. The relations $C_1\mathsf J=-\mathsf JC_1$ show that $u,v\in W$. Consequently $\dim_{\mathbb R}W=4$.

    For $r=1,\cdots,4$, put $T_r=C_rC_1$. These maps are complex-linear (i.e., $T_r\mathsf J=\mathsf J T_r$), self-adjoint, pairwise commuting involutions. Each $T_r$ preserves $W$: if $x\in W$, then
    \[
        C_1(T_rx)=C_1C_rC_1x=C_rC_1^2x=T_rx.
    \]
    The restrictions $T_r|_W$ are therefore commuting self-adjoint operators on the four-dimensional real Euclidean space $W$. By the real spectral theorem they admit a common real orthonormal eigenbasis $u_1,\cdots,u_4$ of $W$, and every eigenvalue is $+1$ or $-1$ because $T_r^2=\operatorname{Id}$.

    By \eqref{eq:real-form-decomposition}, $u_1,\cdots,u_4$ is also a unitary complex basis of $\Sigma$. Indeed, $W$ and $\mathsf JW$ are respectively the $+1$ and $-1$ eigenspaces of the self-adjoint involution $C_1$, so they are orthogonal; hence $h(u_\alpha,u_\beta)=g_\Sigma(u_\alpha,u_\beta)=\delta_{\alpha\beta}$. In this basis $C_1$ is ordinary coordinatewise complex conjugation: any vector $z\in\Sigma$ has the form $z=z_1u_1+\cdots+z_4u_4$, $z_\alpha\in\mathbb C$; the antilinearity of $C_1$ implies $C_1(z)=C_1(\sum_\alpha z_\alpha u_\alpha)=\sum_\alpha\bar z_\alpha C_1u_\alpha=\sum_\alpha\bar z_\alpha u_\alpha$, i.e., in this basis, $C_1(z_1,z_2,z_3,z_4)=(\bar z_1,\bar z_2,\bar z_3,\bar z_4)$. Since $C_r=T_rC_1$, each $C_r$ has the form $C_r(z)=D_r\overline z$ since
    \[
        C_r(z)=T_rC_1z=T_r\bar z=D_r\overline z,
    \]
    where $D_r$ is a real diagonal matrix with diagonal entries in $\{+1,-1\}$. Thus all four $C_r$ belong to the four-dimensional real space of real diagonal antilinear maps.

    They are linearly independent since $\operatorname{tr}_{\mathbb R}(C_rC_s)=8\delta_{rs}$ by \eqref{six:clifford-isometry}. Consequently, they form a basis of the whole real diagonal space.  It follows that the Takagi diagonal map $D$ in \eqref{eq:takagi-conjugation} has a unique expression
    \begin{equation}\label{eq:diagonal-section}
        D=aC_1+bC_2+cC_3+dC_4=\mathsf c(\omega_0), \qquad \omega_0=a e^{123}+b e^{145}+c e^{246}+d e^{356}.
    \end{equation}

    \medskip
    \noindent
    \textbf{6. The one remaining phase is the Hodge phase.}
    We first record the exact relation between Clifford multiplication and the Hodge star. If $I=\{i<j<k\}$ and $I^c$ is its increasingly ordered complement, let $\varepsilon_I$ be determined by $*e^I=\varepsilon_I e^{I^c}$. Then $c_Ic_{I^c}=\varepsilon_I\mathsf J$. Since $c_I^2=\operatorname{Id}$, this gives
    \[
        c_I\mathsf J=c_I(\varepsilon_Ic_Ic_{I^c})=\varepsilon_Ic_{I^c}=\mathsf c(\varepsilon_I e^{I^c})=\mathsf c(*e^I).
    \]
    By linearity,
    \begin{equation}\label{eq:hodge-clifford}
        \mathsf c(*\alpha)=\mathsf c(\alpha)\mathsf J=-\mathsf J\mathsf c(\alpha) \qquad (\alpha\in\Lambda^3V^*).
    \end{equation}
    Here $\mathsf J$ is multiplication by $i$ on $\Sigma$. Therefore
    \begin{equation}\label{eq:hodge-phase-clifford}
        \mathsf c\bigl((\cos\theta+\sin\theta\,*)\alpha\bigr)=(\cos\theta-\sin\theta\,\mathsf J)\mathsf c(\alpha)=e^{-i\theta}\mathsf c(\alpha).
    \end{equation}

    Return to the unitary matrix $U$ in the Takagi factorization. Choose $\lambda\in U(1)$ with $\lambda^4=\det U$ and put $S=\lambda^{-1}U$. Then $S\in SU(4)$, so by \eqref{eq:spin-su4} it is induced by an oriented orthonormal change of coframe. Since $A$ is antilinear, scalar multiplication on its input is conjugated. Hence
    \begin{align*}
        D=U^{-1}AU&=(\lambda S)^{-1}A(\lambda S)\\
        &=\lambda^{-1}\overline\lambda\,S^{-1}AS=\lambda^{-2}S^{-1}AS.
    \end{align*}
    Thus $S^{-1}AS=\lambda^2D$. Choose $\theta$ so that $\lambda^2=e^{-i\theta}$. Combining \eqref{eq:diagonal-section} and \eqref{eq:hodge-phase-clifford}, we obtain
    \[
        S^{-1}\mathsf c(\omega)S=\mathsf c\bigl((\cos\theta+\sin\theta\,*)\omega_0\bigr).
    \]
    Injectivity of $\mathsf c$ and equivariance \eqref{eq:equivariance} show that, in the oriented orthonormal coframe corresponding to $S$, the form $\omega$ has exactly the expression \eqref{six:normal-form-equation}. The decomposition $U(4)=U(1)SU(4)$ also explains why only one continuous scalar phase remains after an oriented orthonormal change of coframe.

    It remains to check the norm statement and the final assertion. On $k$-forms in an oriented Euclidean $n$-space, $*^2=(-1)^{k(n-k)}\operatorname{Id}$. For $n=6$ and $k=3$, therefore, $*^2=-\operatorname{Id}$. The Hodge star is an isometry, so $*^{-1}=*^*$. Since $*^{-1}=-*$ on three-forms, one has $*^*=-*$. Consequently, if $U_\theta=\cos\theta+\sin\theta\,*$, then
    \[
        U_\theta^*U_\theta = (\cos\theta-\sin\theta\,*)(\cos\theta+\sin\theta\,*) = \cos^2\theta\,\operatorname{Id}-\sin^2\theta\,*^2 = \operatorname{Id}.
    \]
    Thus $U_\theta$ is orthogonal. Since $|\omega|=1$, the normal form gives $|\omega_0|=1$. The four displayed basis forms in $\omega_0$ are orthonormal, and hence $a^2+b^2+c^2+d^2=1$. Lemma \ref{lem:HodgeDual} gives $q_3(E)*=*q_3(E)$, so
    \[
        \langle q_3(E)U_\theta\omega_0,U_\theta\omega_0\rangle = \langle q_3(E)\omega_0,\omega_0\rangle.
    \]
    This is an invariance of the operator on forms, and does not assert that $U_\theta$ comes from a change of coframe on $V$.

    Finally, an orthogonal transformation of $V$ sends every decomposable three-form to a decomposable three-form. On the other hand,
    \[
        U_\theta(e^{123})=\cos\theta\,e^{123}+\sin\theta\,e^{456}.
    \]
    If both $\cos\theta$ and $\sin\theta$ are nonzero, this form is not decomposable. Indeed, the contraction map
    \[
        V\longrightarrow\Lambda^2V^*, \qquad v\longmapsto\iota_v\bigl(\cos\theta\,e^{123}+\sin\theta\,e^{456}\bigr)
    \]
    is injective: contractions by $e_1,e_2,e_3$ lie in $\Lambda^2\langle e^1,e^2,e^3\rangle$, whereas contractions by $e_4,e_5,e_6$ lie in $\Lambda^2\langle e^4,e^5,e^6\rangle$. More explicitly, their images are, up to the usual contraction signs,
    \[
        \cos\theta\,e^{23},\quad \cos\theta\,e^{13},\quad \cos\theta\,e^{12},\quad \sin\theta\,e^{56},\quad \sin\theta\,e^{46},\quad \sin\theta\,e^{45},
    \]
    which are linearly independent when both coefficients are nonzero. By contrast, the contraction map of a nonzero decomposable form $\alpha^1\wedge\alpha^2\wedge\alpha^3$ has a three-dimensional kernel, namely $\bigcap_{r=1}^3\ker\alpha^r$. Therefore $U_\theta$ cannot be the action on $\Lambda^3V^*$ of an orthogonal transformation of $V$ whenever $\sin\theta\cos\theta\ne0$. This proves the asserted ``not in general.''
\end{proof}

\subsection{Contractions and mixed terms}

Set $\alpha_s=\iota_s\omega$ and $T_{ab}=\langle\alpha_a,\alpha_b\rangle$. As shown in Section \ref{sec:energy}, $F = \sum_{a,b}R_{ab}T_{ab} = \sum_{a,b}R_{ab}\sum_{i<j}\omega_{aij}\omega_{bij}$. Specializing the second identity in \eqref{alg:contraction-identities} to $p=3$ gives
\begin{equation}\label{six:pair-contraction}
    \langle\mathcal R_E^{(3)}\omega,\omega\rangle=\sum_s\langle\mathcal R_E\alpha_s,\alpha_s\rangle.
\end{equation}
Since $Q=\langle q_3(E)\omega,\omega\rangle=\langle R^{[3]}\omega,\omega\rangle-2\langle\mathcal R_E^{(3)}\omega,\omega\rangle$, it follows that
\begin{equation}\label{eq:three-energy}
    Q=\sum_{a,b}R_{ab}\sum_{i<j}\omega_{aij}\omega_{bij}-2\sum_s\sum_{i<j,\ k<l}E_{ijkl}\omega_{sij}\omega_{skl}.
\end{equation}
Put $D_I:=\sum_{i\in I,\,j\notin I}K_{ij}$ with $|I|=3$.

\subsection{The full four-term expansion and the scalar estimate}
\label{six:scalar-estimate}

For the form $\omega_0$ in \eqref{six:normal-form-equation}, all six contractions are
\begin{equation}\label{six:contraction-table}
    \begin{array}{lll}
        \iota_1\omega_0=a e^{23}+b e^{45},&
        \iota_2\omega_0=-a e^{13}+c e^{46},&
        \iota_3\omega_0=a e^{12}+d e^{56},\\[2pt]
        \iota_4\omega_0=-b e^{15}-c e^{26},&
        \iota_5\omega_0=b e^{14}-d e^{36},&
        \iota_6\omega_0=c e^{24}+d e^{35}.
    \end{array}
\end{equation}
These contractions have mutually disjoint two-form supports, so their Gram matrix is diagonal, i.e., for any $a\neq b$, $\sum_{i<j}(\omega_0)_{aij}(\omega_0)_{bij}=0$. By Lemma \ref{six:normal-form}, $Q=\langle q_3(E)\omega,\omega\rangle=\langle q_3(E)\omega_0,\omega_0\rangle$. Then applying \eqref{eq:three-energy} and canceling the internal sectional curvatures gives
\begin{align}
    Q &= a^2D_{123}+b^2D_{145}+c^2D_{246}+d^2D_{356} \notag\\
    &\quad-4ab E_{2345}+4ac E_{1346}-4ad E_{1256} \notag\\
    &\quad-4bc E_{1526}+4bd E_{1436}-4cd E_{2435}. \label{six:four-term-expansion}
\end{align}

The following table records all the data needed for the six uses of the NW inequality. The column headed $A$ is the sum of the four cross-sectional curvatures in that application; it is not a coefficient of the form.
\begin{equation}\label{six:NW-table}
    \begin{array}{c|c|c|c}
        |\text{coefficient}|&r& (k,l)&A\\ \hline
        |ab|&E_{2345}&(K_{23},K_{45})&K_{24}+K_{25}+K_{34}+K_{35}\\
        |ac|&E_{1346}&(K_{13},K_{46})&K_{14}+K_{16}+K_{34}+K_{36}\\
        |ad|&E_{1256}&(K_{12},K_{56})&K_{15}+K_{16}+K_{25}+K_{26}\\
        |bc|&E_{1526}&(K_{15},K_{26})&K_{12}+K_{16}+K_{25}+K_{56}\\
        |bd|&E_{1436}&(K_{14},K_{36})&K_{13}+K_{16}+K_{34}+K_{46}\\
        |cd|&E_{2435}&(K_{24},K_{35})&K_{23}+K_{25}+K_{34}+K_{45}
    \end{array}
\end{equation}
In each row, Lemma \ref{alg:nw} gives $4|\text{coefficient}|\,|r| \le \frac23|\text{coefficient}|A + \frac83|\text{coefficient}|\sqrt{kl}$. Thus every curvature variable and every mixed term in \eqref{six:four-term-expansion} is accounted for by this table.

To see the scalar estimate directly, collect the opposite rows $ab$ and $cd$. The contribution of the diagonal term $a^2D_{123}+b^2D_{145}+c^2D_{246}+d^2D_{356}$ on $K_{23},K_{45},K_{24},K_{35}$ is
\[
    (c^2+d^2)(K_{23}+K_{45})+(a^2+b^2)(K_{24}+K_{35}).
\]
The contribution of the mixed term on $K_{23},K_{45},K_{24},K_{35}$ is bounded below by
\[
    -\frac23|cd|(K_{23}+K_{45})-\frac23|ab|(K_{24}+K_{35})-\frac83|ab|\sqrt{K_{23}K_{45}}-\frac83|cd|\sqrt{K_{24}K_{35}}.
\]
Hence their contribution on $K_{23},K_{45},K_{24},K_{35}$ is bounded below by
\begin{align}
    &(c^2+d^2-\frac23|cd|)(K_{23}+K_{45})+(a^2+b^2-\frac23|ab|)(K_{24}+K_{35})\notag\\
    &\hspace{8mm} -\frac83|ab|\sqrt{K_{23}K_{45}}-\frac83|cd|\sqrt{K_{24}K_{35}}. \label{six:opposite-pair}
\end{align}
Using $2\sqrt{xy}\le x+y$, the coefficient of either $K_{23}$ or $K_{45}$ is at least
\begin{equation}\label{six:opposite-pair-final}
    c^2+d^2-\frac23|cd|-\frac43|ab| \ge \frac23(c^2+d^2)-\frac23(a^2+b^2) = \frac23-\frac43(a^2+b^2) \ge -\frac23.
\end{equation}
The coefficient for $K_{24}$ and $K_{35}$ is obtained by interchanging $(a,b)$ and $(c,d)$, and satisfies the same bound. The opposite row pairs $\{ac,bd\}$ and $\{ad,bc\}$ account, respectively, for
\[
    \{K_{13},K_{46},K_{14},K_{36}\}, \qquad \{K_{12},K_{56},K_{15},K_{26}\},
\]
with exactly the same estimate. The coefficients of the three sectional curvatures not yet listed arise from the diagonal term $a^2D_{123}+b^2D_{145}+c^2D_{246}+d^2D_{356}$ and the term $A$ in the NW inequality. Specifically,
\begin{equation}\label{six:remaining-coefficients}
    \begin{array}{c|l}
        K_{16}&1-\frac23(|a|+|b|)(|c|+|d|)\\[0.8ex]
        K_{25}&1-\frac23(|a|+|c|)(|b|+|d|)\\[0.8ex]
        K_{34}&1-\frac23(|a|+|d|)(|b|+|c|).
    \end{array}
\end{equation}
For example, $(|a|+|b|)(|c|+|d|) \le 2\sqrt{(a^2+b^2)(c^2+d^2)} \le 1$. Their coefficients are consequently at least $1/3$. Since $S=2\sum_{i<j}K_{ij}$, we have proved
\begin{equation}\label{six:scalar-bound}
    Q\ge-\frac13S, \qquad {\rm i.e.,} ~ q_3(E)\ge-\frac13\operatorname{Scal}(E)\operatorname{Id}.
\end{equation}

\subsection{Two affine estimates retaining the Ricci energy}

\begin{theorem}\label{thm:six-affine}
    Let $E$ be an algebraic curvature tensor on a six-dimensional Euclidean space with nonnegative sectional curvature. For every unit three-form, the quantities in \eqref{six:QF} satisfy
    \begin{equation}\label{six:affine-pair}
        Q\ge-\frac23F, \qquad Q\ge\frac23(F-S).
    \end{equation}
    Equivalently,
    \[
        q_3(E)\ge-\frac23(\operatorname{Ric}(E)^\sharp)^{[3]}, \qquad q_3(E)\ge\frac23\bigl((\operatorname{Ric}(E)^\sharp)^{[3]}-\operatorname{Scal}(E)\operatorname{Id}\bigr).
    \]
    In particular $q_3(E)\ge-\frac13\operatorname{Scal}(E)\operatorname{Id}$, and this scalar coefficient is sharp.
\end{theorem}

\begin{proof}
    Write $\omega$ as in \eqref{six:normal-form-equation}, and put
    \[
        A=|a|,\quad B=|b|,\quad C=|c|,\quad D=|d|, \qquad c_\theta=\cos\theta,\quad s_\theta=\sin\theta.
    \]
    For the Gram matrix $T_{ij}=\langle\iota_i\omega,\iota_j\omega\rangle$, we have
    \begin{equation}\label{six:phase-Gram}
        T=c_\theta^2 T^0+s_\theta^2(\operatorname{Id}-T^0)+\sin(2\theta)U,
    \end{equation}
    where
    \[
        T^0=\operatorname{diag}(a^2+b^2,a^2+c^2,a^2+d^2,b^2+c^2,b^2+d^2,c^2+d^2)
    \]
    and $U$ is symmetric, with its only nonzero entries
    \begin{equation}\label{six:phase-U}
        U_{16}=U_{61}=ab+cd,\qquad U_{25}=U_{52}=-(ac+bd),\qquad U_{34}=U_{43}=ad+bc.
    \end{equation}
    Indeed, the Gram matrix of $\omega_0$ is $T^0$ by \eqref{six:contraction-table}. With the chosen orientation,
    \[
        *\omega_0=a e^{456}+b e^{236}+c e^{135}+d e^{124},
    \]
    whose contractions are
    \[
        \begin{array}{lll}
            \iota_1(*\omega_0)=ce^{35}+de^{24}, & \iota_2(*\omega_0)=be^{36}-de^{14}, & \iota_3(*\omega_0)=-be^{26}-ce^{15},\\
            \iota_4(*\omega_0)=ae^{56}+de^{12}, & \iota_5(*\omega_0)=-ae^{46}+ce^{13}, & \iota_6(*\omega_0)=ae^{45}+be^{23}.
        \end{array}
    \]
    Thus the Gram matrix of $*\omega_0$ is $\operatorname{Id}-T^0$. Taking inner products of these contractions with \eqref{six:contraction-table} gives zero except at the three unordered index pairs $16,25,34$. For instance,
    \[
        T_{16}=\langle\iota_1(c_\theta\omega_0),\iota_6(s_\theta*\omega_0)\rangle+\langle\iota_1(s_\theta*\omega_0),\iota_6(c_\theta\omega_0)\rangle=c_\theta s_\theta(2ab+2cd)=\sin(2\theta)U_{16}.
    \]

    For any three-form, $F=\sum_{i,j}R_{ij}T_{ij}$. Write $F=F_{\mathrm{diag}}+F_{\mathrm{cross}}$, where $F_{\mathrm{diag}}=\sum_{i}R_{ii}T_{ii}$ and $F_{\mathrm{cross}}=2\sum_{i<j}R_{ij}T_{ij}$. Then
    \begin{equation}\label{six:phase-Ricci-cross}
        F_{\mathrm{cross}} = 2\sin(2\theta)\bigl((ab+cd)R_{16}-(ac+bd)R_{25}+(ad+bc)R_{34}\bigr).
    \end{equation}

    We first estimate $Q+\frac23F_{\mathrm{diag}}$. Put $x_{ij}=\sqrt{K_{ij}}$. The six pairs of sectional curvatures in \eqref{six:NW-table} will be used with the following labels:
    \begin{equation}\label{six:remainder-table}
        \begin{array}{c|c|c}
            \text{sectional pair}&\text{coefficient labels}&h\\ \hline
            (23,45)&(a,b)&c^2+d^2\\
            (13,46)&(a,c)&b^2+d^2\\
            (12,56)&(a,d)&b^2+c^2\\
            (15,26)&(b,c)&a^2+d^2\\
            (14,36)&(b,d)&a^2+c^2\\
            (24,35)&(c,d)&a^2+b^2
        \end{array}
    \end{equation}
    Consider the first row. After the NW estimates, the coefficient on $K_{23}$ in $Q+\frac23F_{\mathrm{diag}}$ is
    \[
        \frac43(c_\theta^2a^2+s_\theta^2b^2)+\frac53(c^2+d^2)-\frac23CD.
    \]
    Indeed, since $F_{\mathrm{diag}}=\sum_{i}R_{ii}T_{ii}$, the coefficient of $K_{23}$ in $F_{\mathrm{diag}}$ is $T_{22}+T_{33}$. By \eqref{six:phase-Gram},
    \[
        T_{22}+T_{33} = 2(c_\theta^2a^2+s_\theta^2b^2)+c^2+d^2.
    \]
    Adding $\frac23(T_{22}+T_{33})$ to the coefficient $c^2+d^2-\frac23CD$ from \eqref{six:opposite-pair} gives the displayed coefficient on $K_{23}$ in $Q+\frac23F_{\mathrm{diag}}$. The coefficient on $K_{45}$ is the same with $a,b$ interchanged. There is a negative square-root term $-\frac83ABx_{23}x_{45}$ in \eqref{six:opposite-pair}. Adding this to the displayed coefficient and completing squares gives
    \[
        \frac43\bigl(c_\theta^2(Ax_{23}-Bx_{45})^2+s_\theta^2(Bx_{23}-Ax_{45})^2\bigr)+\left(\frac53(c^2+d^2)-\frac23CD\right)(K_{23}+K_{45}).
    \]
    Since $2CD\le c^2+d^2$, the last coefficient is at least $\frac43(c^2+d^2)=\frac43h$. Substitution of each of the remaining five rows of \eqref{six:remainder-table} gives the corresponding completed-square formula: the two coefficient labels supply the square terms, and the complementary square sum supplies the remainder $\frac43h$.

    Let $\mathcal P$ denote the twelve sectional variables in \eqref{six:remainder-table}, and give each of the two entries in a row the listed weight $h_e$. For the remaining three variables $K_{16},K_{25},K_{34}$, their coefficient in $F_{\mathrm{diag}}$ is exactly one: for example $T_{11}+T_{66}=1$. By \eqref{six:remaining-coefficients}, their coefficient in $Q+\frac23F_{\mathrm{diag}}$ is at least one: for example
    \[
        1-\frac23(A+B)(C+D)+\frac23 \ge 1-\frac23+\frac23 = 1.
    \]
    Dropping only the nonnegative completed squares and excess diagonal terms, we obtain
    \begin{equation}\label{six:retained-remainder}
        Q+\frac23F_{\mathrm{diag}} \ge \frac43\sum_{e\in\mathcal P}h_eK_e+K_{16}+K_{25}+K_{34}.
    \end{equation}

    It remains to estimate \eqref{six:phase-Ricci-cross}. The Jacobi form $E(e_k,\cdot,e_k,\cdot)$ is positive semidefinite. By \eqref{alg:Jacobi-CS},
    \[
        |R_{ij}| \le \sum_{k\ne i,j}|E_{kikj}| \le \sum_{k\ne i,j}\sqrt{K_{ik}K_{jk}}.
    \]
    Consequently
    \begin{equation}\label{six:cross-graph-bound}
        \frac23|F_{\mathrm{cross}}| \le \frac43|\sin(2\theta)|\sum_{\{e,f\}\in\mathcal G}u_{ef}x_ex_f,
    \end{equation}
    where the complete list of edges and their nonnegative weights is
    \begin{equation}\label{six:cross-graph-table}
        \begin{array}{c|c}
            u_{ef}&\{e,f\}\\ \hline
            AB+CD&(12,26),\ (13,36),\ (14,46),\ (15,56)\\
            AC+BD&(12,15),\ (23,35),\ (24,45),\ (26,56)\\
            AD+BC&(13,14),\ (23,24),\ (35,45),\ (36,46)\rlap{.}
        \end{array}
    \end{equation}
    For every edge in this table,
    \begin{equation}\label{six:edge-CS}
        u_{ef}\le\sqrt{h_eh_f}, \qquad 2u_{ef}x_ex_f\le h_ex_e^2+h_fx_f^2.
    \end{equation}
    For example, on $(12,26)$ this is
    \[
        AB+CD\le\sqrt{(B^2+C^2)(A^2+D^2)}.
    \]
    On $(13,36)$ it is $AB+CD\le\sqrt{(B^2+D^2)(A^2+C^2)}$. The other two edges with this weight interchange the factors in these two inequalities. The next two rows are the same two Cauchy--Schwarz inequalities after permuting the four letters. Thus \eqref{six:edge-CS} checks every edge and does not require any sign condition on $a,b,c,d$.

    Every vertex of $\mathcal G$ occurs exactly twice in \eqref{six:cross-graph-table}. Summing \eqref{six:edge-CS} over the twelve edges therefore gives
    \[
        \sum_{\{e,f\}\in\mathcal G}u_{ef}x_ex_f \le \sum_{e\in\mathcal P}h_eK_e.
    \]
    This together with \eqref{six:cross-graph-bound} implies
    \[
        \frac23|F_{\mathrm{cross}}| \le \frac43|\sin(2\theta)|\sum_{e\in\mathcal P}h_eK_e \le \frac43\sum_{e\in\mathcal P}h_eK_e.
    \]
    Hence
    \[
        Q+\frac23F \ge Q+\frac23F_{\mathrm{diag}}-\frac23|F_{\mathrm{cross}}| \ge K_{16}+K_{25}+K_{34} \ge 0.
    \]
    This proves the first inequality of \eqref{six:affine-pair}.

    For the second one, applying the first inequality to $*\omega$ gives
    \[
        Q_3(E,*\omega)\ge-\frac23F_3(E,*\omega).
    \]
    From Lemma \ref{lem:HodgeDual} and \eqref{eq:RicHodgeDual}, $Q_3(E,*\omega)=Q_3(E,\omega)=Q$, and $*R^{[3]}*^{-1}=S\operatorname{Id}-R^{[3]}$ on three-forms in dimension six; hence
    \[
        F_3(E,*\omega) = \langle R^{[3]}(*\omega),*\omega\rangle = \langle *^{-1}R^{[3]}*\omega,\omega\rangle = \langle(S\operatorname{Id}-R^{[3]})\omega,\omega\rangle = S-F.
    \]
    We obtain $Q\ge-\frac23(S-F)$ as asserted. Taking the arithmetic mean of the two inequalities gives $Q\ge-S/3$.

    For sharpness, decompose $V=W^4\oplus\mathbb R^2$ and choose an orthogonal complex structure $J$ on $W$. Let $t>0$ and let $E=tE_J$ on $W$. Extend $E$ by zero. Let $\Omega_J$ be its K\"ahler form and let $\zeta$ be a unit covector in the flat summand. By Lemma \ref{lem:model-operator}, $\operatorname{Ric}(E)^\sharp|_W=3t\operatorname{Id}_W$, $S=12t$ and $|\Omega_J|^2=2$. Lemma \ref{lem:model-powers} gives $q_3(E)(\zeta\wedge\Omega_J)=-4t\zeta\wedge\Omega_J$. By the Leibniz rule,
    \[
        R^{[3]}(\zeta\wedge\Omega_J)=\zeta\wedge R^{[2]}\Omega_J=6t\zeta\wedge\Omega_J.
    \]
    Thus the unit form $\omega=\zeta\wedge\Omega_J/\sqrt2$ has $(S,F,Q)=(12t,6t,-4t)$. It realizes equality in both affine inequalities and in the scalar estimate, proving sharpness.
\end{proof}

\begin{corollary}\label{six:low-dimensional-sharp}\label{six:sharp}
    Let $E$ be an algebraic curvature tensor with nonnegative sectional curvature on an $n$-dimensional Euclidean space $V$. Then
    \[
        q_3(E)\ge
        \begin{cases}
            0, & n=4,\\[2pt]
            -\dfrac13\operatorname{Scal}(E)\mathrm{Id}_{\Lambda^3 V^*}, & n=5,6.
        \end{cases}
    \]
    Moreover, these coefficients are sharp.
\end{corollary}

\begin{proof}
    In dimension four, Hodge duality identifies $q_3(E)$ with $q_1(E)=\operatorname{Ric}(E)^\sharp\ge0$. A nonzero constant-curvature tensor on a two-plane, extended by zero to its orthogonal complement, has positive scalar curvature and a Ricci-null covector. Its Hodge dual is a nonzero three-form of zero energy, so the coefficient zero is sharp. In dimension five it identifies $q_3(E)$ with $q_2(E)$, for which the two-form estimate \eqref{intro:ge} gives the sharp coefficient $1/3$. Alternatively, the preceding model with a single flat direction already realizes $Q/S=-1/3$ in dimension five. Dimension six is Theorem \ref{thm:six-affine}.
\end{proof}

\section{Dimension lifting and scalar estimates}\label{sec:scalar-bounds}

\subsection{A fixed-degree lift from dimension six}

The six-dimensional affine estimate improves one branch of the preceding family. We state explicitly the dependence on that estimate, which is proved separately in Theorem \ref{thm:six-affine}.

\begin{theorem}[The dimension lift for three-forms]\label{alg:three-lift}\label{thm:three-lift}
    For $n\geq6$, define
    \begin{equation}\label{alg:tn}
        t_n=\frac{19n-82}{21n-78}.
    \end{equation}
    Every algebraic curvature tensor $E$ with nonnegative sectional curvature satisfies
    \begin{equation}\label{alg:three-lift-equation}
        Q_3(E,\omega)\geq t_n\bigl(F_3(E,\omega)-S|\omega|^2\bigr).
    \end{equation}
\end{theorem}

\begin{proof}
    We first calculate the dimension-lifting rule for an arbitrary affine inequality $uQ_3+vF_3+wSN\ge0$ known in dimension $n-1$. Let $\eta=*\omega$, of degree $p=n-3$. The form $\alpha_i=\iota_i\eta$ belongs to $\Lambda^{n-4}H_i^*$. Its Hodge dual in $H_i$ is a three-form $\beta_i$. By \eqref{eq:hodge-energy} on $H_i$,
    \[
        Q_{n-4}(E_i,\alpha_i)=Q_3(E_i,\beta_i), \qquad F_3(E_i,\beta_i)=\operatorname{Scal}(E_i)|\alpha_i|^2-F_{n-4}(E_i,\alpha_i).
    \]
    Likewise, on the full space, $Q_{n-3}(E,\eta)=Q_3(E,\omega)$ and $F_{n-3}(E,\eta)=SN-F_3(E,\omega)$.

    Put $A_n=n^2-3n-7$. Substitution of $p=n-3$ into Proposition \ref{alg:averages} gives
    \begin{align*}
        \overline{\sum_iQ_3(E_i,\beta_i)}&=\frac{A_nQ+3(SN-F)}{n+2},\\
        \overline{\sum_iF_3(E_i,\beta_i)}&=\frac{3SN+(A_n+4)F-3Q}{n+2},\\
        \overline{\sum_i\operatorname{Scal}(E_i)|\beta_i|^2}&=\frac{(A_n+3)SN+4F}{n+2}.
    \end{align*}
    Since we have assumed $uQ_3+vF_3+wSN\ge0$,
    \[
        uQ_3(E_i,\beta_i)+vF_3(E_i,\beta_i)+w\operatorname{Scal}(E_i)|\beta_i|^2\ge0.
    \]
    Summing over all $i$ and taking the Haar average as in the proof of Proposition \ref{alg:affine-recurrence}, the coefficient vector evolves by
    \begin{equation}
        \begin{pmatrix}u'\\v'\\w'\end{pmatrix}=(A_n\operatorname{Id}+M)\begin{pmatrix}u\\v\\w\end{pmatrix},
        \qquad
        M=\begin{pmatrix}0&-3&0\\-3&4&4\\3&3&3\end{pmatrix}.
    \end{equation}

    Theorem \ref{thm:six-affine} gives the initial vector $(1,-2/3,2/3)$ in dimension six. Decomposition gives
    \[
        (1,-2/3,2/3)=\frac56(1,-1,1)+\frac16(1,1,-1)+0\cdot(3,-7,-3).
    \]
    The three displayed vectors are eigenvectors of $M$, with eigenvalues $3$, $-3$ and $7$. Our induction is on $n$, and $p=3$ is fixed. But $M$ is independent of $n$, and $A_n\operatorname{Id}$ is a scalar matrix, it follows that in the induction on $n$, the directions of the three eigenvectors of $A_n\operatorname{Id}+M$ remain invariant, while only their lengths change. Moreover, since the initial vector $(1,-2/3,2/3)$ has no component along the eigenvector $(3,-7,-3)$, it suffices to consider the remaining two eigenvectors and their eigenvalues. At dimension $j$, the two related eigenvalues of $A_j\operatorname{Id}+M$ are
    \[
        A_j+3=(j-4)(j+1),\qquad A_j-3=(j-5)(j+2).
    \]
    All are positive for $j\geq7$. Their accumulated ratio is
    \[
        r_n=\prod_{j=7}^n\frac{(j-5)(j+2)}{(j-4)(j+1)}=\frac{n+2}{4(n-4)}.
    \]
    The empty product at $n=6$ is one, in agreement with this expression. After dividing by the positive first coordinate, the resulting inequality is
    \[
        Q-\frac{5-r_n}{5+r_n}F+\frac{5-r_n}{5+r_n}SN\geq0.
    \]
    Since $(5-r_n)/(5+r_n)=(19n-82)/(21n-78)$, this proves the theorem.
\end{proof}

Applying \eqref{alg:three-lift-equation} to $*\omega$ gives $Q_{n-3}(E,\omega)\geq-t_nF_{n-3}(E,\omega)$. In operator form, this is $q_{n-3}(E)+t_nR^{[n-3]}\ge0$. This proves the improved Ricci-weighted inequality in codimension three. For subsequent statements, we can therefore define the improved coefficient
\begin{equation}\label{alg:kappa-improved}
    \kappa_{n,p}=
    \begin{cases}
        t_n,&n-p=3\text{ and }n\ge6,\\
        \kappa^0_{n,p},&\text{otherwise}.
    \end{cases}
\end{equation}
It is stronger than the raw odd-codimension branch.

\subsection{Scalar estimates and the sharp even branch}

\begin{theorem}[Explicit estimates in all degrees]\label{alg:scalar-all}\label{alg:even-sharp}
    For $2\leq p\leq n-2$, put
    \begin{equation}\label{alg:Cnp}
        C_{n,p}=\frac{\kappa_{n,p}\kappa_{n,n-p}}{\kappa_{n,p}+\kappa_{n,n-p}}.
    \end{equation}
    Then
    \begin{equation}\label{alg:scalar-all-equation}
        q_p(E)\geq-C_{n,p}S\operatorname{Id}_{\Lambda^p}.
    \end{equation}
    For $p=0,1,n-1,n$, one can take $C_{n,p}=0$. If $n$ and $p$ are both even, the coefficient in \eqref{alg:Cnp} equals $p(n-p)/(3n)$, and is sharp on the entire cone of tensors with nonnegative sectional curvature.
\end{theorem}

\begin{proof}
    Write $a=\kappa_{n,p}$ and $b=\kappa_{n,n-p}$; both are positive in the indicated interior range. The Ricci-weighted inequalities for $\omega$ and $*\omega$ are $Q\geq-aF$ and $Q\geq-b(SN-F)$, respectively. Multiply the first by $b/(a+b)$ and the second by $a/(a+b)$. Adding cancels $F$ and yields $Q\geq-abSN/(a+b)$. For degrees zero and $n$, the curvature term vanishes. In degree one it is $R\geq0$, and degree $n-1$ follows by Hodge duality.

    For the sharpness when $n$ and $p$ are both even, let $n=2m$, $p=2r$, and take $E=E_J$. Put
    \[
        \omega=\frac{\Omega_J^r}{r!\sqrt{\binom mr}}, \qquad \Omega_J=\sum_{i=1}^m\theta_i, \qquad \theta_i=e^{2i-1}\wedge e^{2i}.
    \]
    By Lemma \ref{lem:model-powers}, $|\omega|=1$ and $q_{2r}(E_J)\omega=-4r(m-r)\omega$. Thus $Q=-4r(m-r)$. Lemma \ref{lem:model-operator} gives $S(E_J)=6m$, so the ratio is $4r(m-r)/(6m)=p(n-p)/(3n)$. Hence the coefficient cannot be decreased. The cases $p=0,n$ have the trivial sharp coefficient zero.
\end{proof}

In degree two, \eqref{alg:Cnp} recovers \cite{Ge}'s coefficient. In degree three, for $n\geq7$, it becomes
\begin{equation}\label{alg:Gamma3}
    \Gamma_n=\frac{t_n\kappa_n}{t_n+\kappa_n}, \qquad \kappa_n=
    \begin{cases}
        (n-3)/3,&n\text{ odd},\\[3pt]
        (2n^2-6n-5)/(3(2n+1)),&n\text{ even}.
    \end{cases}
\end{equation}

\subsection{An exact Jacobi remainder}

The preceding induction follows the four-plane route. There is also a useful independent identity that verifies the even-degree constant and describes exactly what would be needed for a stronger odd-degree estimate.

\begin{proposition}[Jacobi defect identity]\label{alg:Jacobi-defect}
    For $1\leq p\leq n-1$, define
    \begin{equation}\label{alg:defect-definition}
        \mathscr A_p(E,\omega)=n\int_{S^{n-1}}\left[|\iota_v\omega|^2\operatorname{tr}J_v-\langle J_v^{[p-1]}\iota_v\omega,\iota_v\omega\rangle\right]d\mu(v),
    \end{equation}
    where $J_v=E(v,\cdot,v,\cdot)|_{v^\perp}$. Then
    \begin{equation}\label{alg:defect-identity}
        (n+2)\mathscr A_p=3Q+p(SN-F)
    \end{equation}
    and
    \begin{equation}\label{alg:scalar-defect-identity}
        Q+\frac{p(n-p)}{3n}SN=\frac{n+2}{3n}\left[(n-p)\mathscr A_p(E,\omega)+p\mathscr A_{n-p}(E,*\omega)\right].
    \end{equation}
    In particular, under sectional nonnegativity, $\mathscr A_p\geq0$ and
    \begin{equation}\label{alg:all-cases-estimate}
        Q\ge-\frac{p(n-p)}{3n}SN.
    \end{equation}
\end{proposition}

\begin{proof}
    Fix $v$ and choose an eigenbasis of $J_v$ on $v^\perp$, with eigenvalues $\lambda_j$. For $\alpha=\iota_v\omega$,
    \[
        |\alpha|^2\operatorname{tr}J_v-\langle J_v^{[p-1]}\alpha,\alpha\rangle = \sum_{j=1}^{n-1}\lambda_j\left(|\alpha|^2-|\iota_j\alpha|^2\right) = \sum_{j=1}^{n-1}\lambda_j|e^j\wedge\alpha|^2.
    \]
    The second equality follows from $\varepsilon_j\iota_i+\iota_i\varepsilon_j=\delta_{ij}\operatorname{Id}$, $\iota_i^*=\varepsilon_i$ and $\varepsilon_i^*=\iota_i$.

    The two integrals in \eqref{alg:defect-definition} were computed in the proof of Proposition \ref{alg:averages}. They are respectively $(pSN+2F)/(n+2)$ and $((p+2)F-3Q)/(n+2)$. Their difference gives \eqref{alg:defect-identity}. For $p=1$ the second integral is zero and $Q=F$, giving the same identity directly. Applying \eqref{alg:defect-identity} to $*\omega$, we have
    \[
        (n+2)\mathscr A_{n-p}(E,*\omega)=3Q+(n-p)\bigl(SN-(SN-F)\bigr)=3Q+(n-p)F.
    \]
    This together with \eqref{alg:defect-identity} gives \eqref{alg:scalar-defect-identity}. Under sectional nonnegativity, each $\lambda_j\ge0$; therefore, the integrand in the definition of $\mathscr A_p$ is nonnegative, which together with \eqref{alg:scalar-defect-identity} gives \eqref{alg:all-cases-estimate}.
\end{proof}

Thus $p(n-p)/(3n)$ is a valid scalar coefficient for all parities; the earlier coefficients improve it in the odd branches. Equality in this rough bound has an exact intrinsic necessary and sufficient condition: equality in \eqref{alg:all-cases-estimate} holds if and only if
\begin{align}
    \operatorname{im}J_v&\subset\{y\in v^\perp:y^\flat\wedge\iota_v\omega=0\},\label{alg:defect-equality-1}\\
    \operatorname{im}J_v&\subset\{y\in v^\perp:y^\flat\wedge\iota_v(*\omega)=0\}.\label{alg:defect-equality-2}
\end{align}
Indeed, if the equality holds, both nonnegative defects must vanish: $\mathscr A_p(E,\omega)=\mathscr A_{n-p}(E,*\omega)=0$. Since the integrand is continuous and nonnegative, it is identically zero, so each term $\lambda_j|e^j\wedge\alpha|^2=0$. In particular, if $\lambda_j>0$, then $e^j\wedge\alpha=0$. In the eigenbasis, $J_v(y)=\sum_{j=1}^{n-1}\lambda_j\langle y,e_j\rangle e_j$. Therefore, $\operatorname{im}J_v$ is spanned precisely by those eigenvectors $e_j$ corresponding to $\lambda_j>0$. That is, any $y\in\operatorname{im}J_v$ can be written as $y=\sum_{j\in\{\lambda_j>0\}}c_je_j$. Thus $y^\flat\wedge\alpha=\sum_{j\in\{\lambda_j>0\}}c_j(e^j\wedge\alpha)=0$, which gives \eqref{alg:defect-equality-1}. \eqref{alg:defect-equality-2} is similar. Conversely, if these two inclusions hold, then for all $\lambda_j>0$, the corresponding $e^j$ lies in $\operatorname{im}J_v$. Hence
\[
    e^j\wedge\iota_v\omega=0 \quad \text{and} \quad e^j\wedge\iota_v(*\omega)=0.
\]
Consequently, each term $\lambda_j|e^j\wedge\iota_v\omega|^2$ vanishes, the integrand vanishes, the integral vanishes, and equality holds. Turning this criterion into a classification requires additional algebra; the identity by itself does not supply such a classification.

For three-forms, \eqref{alg:defect-identity} gives $(n+2)\mathscr A_3=3Q+3(SN-F)$. Hence the proposed stronger support inequality $Q\geq-\frac23(SN-F)$ is equivalent to the precise additional requirement
\[
    (n+2)\mathscr A_3(E,\omega)\geq SN-F.
\]
Nonnegativity of $\mathscr A_3$ alone does not prove this requirement. This distinction isolates an actual missing estimate in the unresolved higher-dimensional sharp problem.

\section{Seven-dimensional spinor energies and an explicit improvement}
\label{sec:seven-dimensions}

In dimension seven, retaining the full spinor data gives an improvement over the estimate obtained by averaging the six-dimensional theorem. The resulting scalar coefficient is explicit, but is not asserted to be sharp. We continue to use
\[
    S=\operatorname{Scal}(E), \qquad F=\langle(\operatorname{Ric}(E)^{\sharp})^{[3]}\omega,\omega\rangle, \qquad Q=\langle q_3(E)\omega,\omega\rangle, \qquad |\omega|=1.
\]

\subsection{The Clifford model and its normalization}

Let $\Sigma$ be an eight-dimensional real irreducible Clifford module for $\mathrm{Cl}_{0,7}$. Choose an invariant Euclidean inner product on $\Sigma$, and write $c_i=c(e_i)$ for Clifford multiplication. Our conventions are
\begin{equation}\label{spin:Clifford-relations}
    c_i^*=-c_i, \qquad c_ic_j+c_jc_i=-2\delta_{ij}\operatorname{Id}.
\end{equation}
In particular, $c_i^2=-\operatorname{Id}$, and every $c_i$ is invertible with $c_i^{-1}=-c_i$.

As in the proof of Lemma \ref{six:normal-form}, the map
\begin{equation}\label{spin:Clifford-map}
    \mathsf c:\Lambda^3V^*\longrightarrow\operatorname{Sym}^2_0(\Sigma), \qquad \mathsf c(\omega)=\sum_{i<j<k}\omega_{ijk}c_ic_jc_k
\end{equation}
is an isomorphism and satisfies
\begin{equation}\label{spin:trace-pairing}
    \operatorname{tr}\bigl(\mathsf c(\eta)\mathsf c(\xi)\bigr)=8\langle\eta,\xi\rangle,
\end{equation}
where $\operatorname{Sym}^2_0(\Sigma)$ denotes the real vector space of symmetric trace-free endomorphisms of $\Sigma$ and has real dimension $\frac{8\cdot(8+1)}2-1=35$. Also, $\dim\Lambda^3V^*=\binom73=35$. Every threefold monomial is symmetric, squares to the identity, and has trace zero: taking $l\notin I=\{i<j<k\}$ yields $c_lc_Ic_l^{-1}=-c_I$, so $\operatorname{tr}(c_I)=0$. Hence $\mathsf c(\omega)$ is indeed in $\operatorname{Sym}^2_0(\Sigma)$.

For a unit form $\omega$, put $A=\mathsf c(\omega)$, and choose any orthonormal eigenbasis $f_1,\cdots,f_8$ of this real symmetric matrix:
\begin{equation}\label{spin:eigenvalues}
    Af_a=s_af_a,\qquad \sum_a s_a=\operatorname{tr}(A)=0,\qquad \sum_a s_a^2=\operatorname{tr}(A^2)=8|\omega|^2=8.
\end{equation}
This diagonalization is a change of basis of the spinor module. It is not generally induced by a change of orthonormal frame on $V$; the Clifford matrices must therefore be retained in this same spinor basis.

For $i<j$, set $\Gamma_{ij}=c_ic_j$ and $B_{ij}=\Gamma_{ij}/2$. Direct use of \eqref{spin:Clifford-relations} gives
\[
    [B_{ij},c_i]=c_j, \qquad [B_{ij},c_j]=-c_i, \qquad [B_{ij},c_k]=0\quad(k\ne i,j).
\]
On forms, define $\rho_{ij}=\varepsilon_j\iota_i-\varepsilon_i\iota_j$. For example, $\rho_{ij}e^i=e^j$, $\rho_{ij}e^j=-e^i$, and $\rho_{ij}e^k=0$ for $k\ne i,j$. The commutator is a derivation: $[B,XYZ]=[B,X]YZ+X[B,Y]Z+XY[B,Z]$. Consequently, $[B_{ij},c_ac_bc_c]=\mathsf c(\rho_{ij}e^{abc})$, and by linearity $\mathsf c(\rho_{ij}\omega)=[B_{ij},A]$. Applying this twice and using \eqref{eq:Casimir-formula},
\begin{equation}\label{spin:double-commutator}
    \mathsf c(q_3(E)\omega)=-\sum_{i<j,\,k<l}E_{ijkl}[B_{ij},[B_{kl},A]].
\end{equation}

We will also use the scalar Clifford contraction
\begin{equation}\label{spin:scalar-Clifford}
    \sum_{i<j,\,k<l}E_{ijkl}\Gamma_{ij}\Gamma_{kl}=-\frac S2\operatorname{Id}.
\end{equation}
To prove this, group the terms according to the intersection of the two unordered pairs $\{i,j\}$ and $\{k,l\}$. If the pairs coincide, then $\Gamma_{ij}^2=c_ic_jc_ic_j=-\operatorname{Id}$, so these terms contribute $-\sum_{i<j}E_{ijij}\operatorname{Id}=-S\operatorname{Id}/2$. If the pairs have exactly one common index, the two possible orders cancel. For example, $\Gamma_{ij}\Gamma_{ik}=c_jc_k$ and $\Gamma_{ik}\Gamma_{ij}=-c_jc_k$, while $E_{ikij}=E_{ijik}$; thus their sum is zero. It remains to consider four distinct indices. Fix $i<j<k<l$ and let $C=c_ic_jc_kc_l$. The three pairings of these indices give
\[
    \Gamma_{ij}\Gamma_{kl}=C, \qquad \Gamma_{ik}\Gamma_{jl}=-C, \qquad \Gamma_{il}\Gamma_{jk}=C.
\]
The reversed pair orders give the same products because two bivectors on disjoint index sets commute. Hence the total coefficient of $C$ is a multiple of $E_{ijkl}-E_{ikjl}+E_{iljk}$, which vanishes by Bianchi. Thus all four-distinct-index terms vanish, proving \eqref{spin:scalar-Clifford}.

\subsection{The two energies in a common spinor eigenbasis}

For $a\ne b$, define a vector and a two-form by
\begin{equation}\label{spin:u-Theta}
    u_{ab}=\sum_i\langle c_i f_b,f_a\rangle e_i, \qquad \Theta_{ab}=\sum_{i<j}\langle c_ic_j f_b,f_a\rangle e^{ij}.
\end{equation}
By \eqref{spin:Clifford-relations}, $u_{ba}=-u_{ab}$ and $\Theta_{ba}=-\Theta_{ab}$. Therefore the quantities
\begin{equation}\label{spin:edge-weights}
    k_{ab}=\langle\mathcal R_E\Theta_{ab},\Theta_{ab}\rangle, \qquad r_{ab}=\operatorname{Ric}(E)(u_{ab},u_{ab})
\end{equation}
are symmetric in $a,b$.

\begin{lemma}\label{spin:Theta-geometry}
    For every pair $a\ne b$, $u_{ab}$ is a unit vector, $\iota_{u_{ab}}\Theta_{ab}=0$, and $\Theta_{ab}$ is, up to its overall sign, the K\"ahler form of an orthogonal complex structure on $u_{ab}^\perp$. In particular $|\Theta_{ab}|^2=3$.
\end{lemma}

\begin{proof}
    We first establish the standard cross-product identities associated with a unit spinor. Let $f\in\Sigma$ be a unit vector. For $i,j\in\{1,\cdots,7\}$ with $i\neq j$, $c_ic_j$ is skew-symmetric, so $\langle c_if,c_jf\rangle=-\langle f,c_ic_jf\rangle=0$. When $i=j$, we have $\langle c_if,c_if\rangle=-\langle f,c_i^2f\rangle=\langle f,f\rangle=1$. Thus $c_1f,\cdots,c_7f$ are orthonormal. Since $\langle c_if,f\rangle=-\langle f,c_if\rangle=0$ and $\dim\Sigma=8$, the ordered list $\{f,c_1f,\cdots,c_7f\}$ is an orthonormal basis of $\Sigma$.

    Define
    \[
        \varphi_f(x,y,z)=\langle c(x)c(y)c(z)f,f\rangle, \qquad x,y,z\in V,
    \]
    where $c(x)=\sum_i x_i c_i$. We check that $\varphi_f$ is alternating. Using $c(y)c(x)=-c(x)c(y)-2\langle x,y\rangle\operatorname{Id}$,
    \begin{align*}
        \varphi_f(y,x,z)&=\langle c(y)c(x)c(z)f,f\rangle\\
        &=-\langle c(x)c(y)c(z)f,f\rangle-2\langle x,y\rangle\langle c(z)f,f\rangle\\
        &=-\varphi_f(x,y,z),
    \end{align*}
    where $\langle c(z)f,f\rangle=0$ by the skew-symmetry of $c(z)$. The same argument for the other adjacent transposition proves full alternation.

    Define $\times_f$ by $\langle x\mathbin{\times_f}y,z\rangle=\varphi_f(x,y,z)$. We now expand $c(x)c(y)f$ in the orthonormal basis $\{f,c_1f,\cdots,c_7f\}$. Its coefficient along $f$ is
    \[
        \langle c(x)c(y)f,f\rangle=\frac12\Bigl(\langle c(x)c(y)f,f\rangle+\langle f,c(y)c(x)f\rangle\Bigr)=\frac12\langle-2\langle x,y\rangle\operatorname{Id}(f),f\rangle=-\langle x,y\rangle,
    \]
    For a unit vector $z$, its coefficient along $c(z)f$ is
    \[
        \langle c(x)c(y)f,c(z)f\rangle=-\langle c(z)c(x)c(y)f,f\rangle=-\varphi_f(z,x,y)=-\varphi_f(x,y,z)=-\langle x\mathbin{\times_f}y,z\rangle.
    \]
    Hence
    \begin{equation}\label{spin:cross-product}
        c(x)c(y)f=-\langle x,y\rangle f-c(x\mathbin{\times_f}y)f.
    \end{equation}

    Taking norms in \eqref{spin:cross-product} and using the orthogonality of $f$ and $c(V)f$ gives
    \begin{equation}\label{spin:cross-norm}
        |x\mathbin{\times_f}y|^2=|x|^2|y|^2-\langle x,y\rangle^2.
    \end{equation}
    Here $|c(x)c(y)f|=|x||y|$ since $c(x)^*c(x)=|x|^2\operatorname{Id}$ and similarly for $y$. Alternation also shows that the cross product is orthogonal to both inputs; for example, $\langle x\mathbin{\times_f}y,x\rangle=\varphi_f(x,y,x)=0$ and likewise $\langle x\mathbin{\times_f}y,y\rangle=0$.

    Let $v$ be a unit vector. On $v^\perp$, define $J_v(y)=v\mathbin{\times_f}y$. The alternating property shows that $J_v$ is skew-adjoint:
    \[
        \langle J_v y,z\rangle=\varphi_f(v,y,z)=-\varphi_f(v,z,y)=-\langle y,J_vz\rangle.
    \]
    Equation \eqref{spin:cross-norm} gives $|J_vy|=|y|$ for $y\perp v$. Thus $J_v$ is an isometry of $v^\perp$. Since it is skew-adjoint, $J_v^*=-J_v$, while isometry gives $J_v^*J_v=\operatorname{Id}$; hence $J_v^2=-\operatorname{Id}_{v^\perp}$. Therefore $J_v$ is an orthogonal complex structure on $v^\perp$, and its associated K\"ahler form is
    \begin{equation}\label{spin:Kahler-form}
        \langle J_vy,z\rangle=\varphi_f(v,y,z)=\iota_v\varphi_f(y,z).
    \end{equation}

    Now take $f=f_a$. Since $f_b\perp f_a$, and the map $V\to f_a^\perp$, $v\mapsto c(v)f_a$, is an isometry, there is a unique unit vector $v$ such that $f_b=c(v)f_a$. Substituting this into \eqref{spin:u-Theta}, and using \eqref{spin:cross-product}, gives
    \[
        \langle u_{ab},e_i\rangle=\langle c_i f_b,f_a\rangle=\langle c(e_i)c(v)f_a,f_a\rangle=-\langle e_i,v\rangle-\langle c(e_i\mathbin{\times_{f_a}}v)f_a,f_a\rangle=-\langle e_i,v\rangle,
    \]
    and
    \[
        \langle\Theta_{ab},e^{ij}\rangle=\langle c_ic_jf_b,f_a\rangle=\langle c(e_i)c(e_j)c(v)f_a,f_a\rangle=\varphi_{f_a}(e_i,e_j,v)=\varphi_{f_a}(v,e_i,e_j)=\langle\iota_v\varphi_{f_a},e^{ij}\rangle.
    \]
    Hence $u_{ab}=-v$ and $\Theta_{ab}=\iota_v\varphi_{f_a}$. In particular, $\iota_{u_{ab}}\Theta_{ab}=-\iota_v\iota_v\varphi_{f_a}=0$. By \eqref{spin:Kahler-form}, $\Theta_{ab}$ is the K\"ahler form of the orthogonal complex structure $J_v$ on $u_{ab}^\perp$. Replacing $v$ by $-v$ changes the form by an overall sign, which explains the stated sign ambiguity.

    Finally, if $e_1',\cdots,e_6'$ is an orthonormal basis of $u_{ab}^\perp$ adapted to $J_v$, then $\Theta_{ab}=\pm(e^{1'2'}+e^{3'4'}+e^{5'6'})$, so $|\Theta_{ab}|^2=1+1+1=3$.
\end{proof}

\begin{lemma}\label{spin:energy-identities}
    With the above notation,
    \begin{align}
        Q&=\frac1{16}\sum_{a<b}(s_a-s_b)^2 k_{ab},\label{spin:Q-energy}\\
        S-F&=\frac1{16}\sum_{a<b}(s_a-s_b)^2 r_{ab}.\label{spin:Ricci-energy}
    \end{align}
    The coefficients $k_{ab}$ and $r_{ab}$ satisfy the row-sum identities
    \begin{equation}\label{spin:row-sums}
        \sum_{b\ne a}k_{ab}=\frac S2, \qquad \sum_{b\ne a}r_{ab}=S.
    \end{equation}
\end{lemma}

\begin{proof}
    Using \eqref{spin:trace-pairing} and \eqref{spin:double-commutator}, we have
    \[
        Q = \langle q_3(E)\omega,\omega\rangle = \frac18\operatorname{tr}\bigl(\mathsf c(q_3(E)\omega)\mathsf c(\omega)\bigr) = -\frac18\sum_{i<j,\,k<l}E_{ijkl}\operatorname{tr}\bigl([B_{ij},[B_{kl},A]]A\bigr).
    \]
    For arbitrary endomorphisms $X,Y,A$, cyclicity of trace gives
    \[
        \operatorname{tr}([X,[Y,A]]A)=\operatorname{tr}\bigl((X[Y,A]-[Y,A]X)A\bigr)=\operatorname{tr}\bigl([Y,A](AX-XA)\bigr)=-\operatorname{tr}\bigl([Y,A][X,A]\bigr).
    \]
    Applying this with $X=B_{ij}$ and $Y=B_{kl}$ yields
    \begin{equation}\label{spin:Q-trace}
        Q=\frac18\sum_{i<j,\,k<l}E_{ijkl}\operatorname{tr}\bigl([B_{ij},A][B_{kl},A]\bigr).
    \end{equation}

    In the eigenbasis $Af_a=s_af_a$,
    \begin{align*}
        [B_{ij},A]_{ab}&=\langle f_a,[B_{ij},A]f_b\rangle=\langle f_a,B_{ij}Af_b\rangle-\langle Af_a,B_{ij}f_b\rangle\\
        &=(s_b-s_a)\langle f_a,B_{ij}f_b\rangle=(s_b-s_a)(B_{ij})_{ab}.
    \end{align*}
    Because $B_{ij}$ is skew-symmetric and $A$ is symmetric, $[B_{ij},A]$ is symmetric. Expanding the trace in \eqref{spin:Q-trace} gives
    \[
        Q=\frac18\sum_{a,b}(s_a-s_b)^2\sum_{i<j,\,k<l}E_{ijkl}(B_{ij})_{ab}(B_{kl})_{ab}=\frac14\sum_{a<b}(s_a-s_b)^2\sum_{i<j,\,k<l}E_{ijkl}(B_{ij})_{ab}(B_{kl})_{ab}.
    \]
    Since $B_{ij}=\Gamma_{ij}/2$,
    \begin{align*}
        \sum_{i<j,\,k<l}E_{ijkl}(B_{ij})_{ab}(B_{kl})_{ab}&=\frac14\sum_{i<j,\,k<l}E_{ijkl}\langle\Gamma_{ij}f_b,f_a\rangle\langle\Gamma_{kl}f_b,f_a\rangle\\
        &=\frac14\sum_{i<j,\,k<l}E_{ijkl}(\Theta_{ab})_{ij}(\Theta_{ab})_{kl}\\
        &=\frac14\langle\mathcal R_E\Theta_{ab},\Theta_{ab}\rangle=\frac14k_{ab}.
    \end{align*}
    This proves \eqref{spin:Q-energy}.

    For any symmetric endomorphism $B$ of $V$, there is an identity
    \begin{equation}\label{spin:Ricci-Clifford}
        \mathsf c(B^{[3]}\omega)=\frac12\Bigl((\operatorname{tr}B)A-\sum_{i,j}B_{ij}c_iAc_j\Bigr).
    \end{equation}
    Both sides are equivariant under orthogonal changes of basis, so it is enough to verify the formula when $B$ is diagonal, say $Be_r=\lambda_re_r$. For a threefold monomial $C=c_ic_jc_k$, conjugation by $c_r$ gives
    \[
        c_rCc_r=
        \begin{cases}
            -C,&r\in\{i,j,k\},\\
            +C,&r\notin\{i,j,k\}.
        \end{cases}
    \]
    Consequently,
    \[
        \frac12\Bigl((\operatorname{tr}B)C-\sum_r\lambda_rc_rCc_r\Bigr) = \frac12\Bigl(\sum_r\lambda_r-\sum_{r\notin\{i,j,k\}}\lambda_r+\sum_{r\in\{i,j,k\}}\lambda_r\Bigr)C = (\lambda_i+\lambda_j+\lambda_k)C.
    \]
    This is exactly $\mathsf c(B^{[3]}e^{ijk})$, proving \eqref{spin:Ricci-Clifford}.

    Now put $B=\operatorname{Ric}(E)^\sharp$. Pair \eqref{spin:Ricci-Clifford} with $A$ and use the trace pairing:
    \begin{align*}
        F&=\frac18\operatorname{tr}(\mathsf c(B^{[3]}\omega)A)=\frac1{16}\Bigl((\operatorname{tr}B)\operatorname{tr}(A^2)-\sum_{i,j}B_{ij}\operatorname{tr}(c_iAc_jA)\Bigr)\\
        &=\frac{S}2-\frac1{16}\sum_{i,j}B_{ij}\operatorname{tr}(c_iAc_jA),
    \end{align*}
    where the last equality uses $\operatorname{tr}B=S$ and $\operatorname{tr}(A^2)=8$. We compute the remaining trace in the eigenbasis. Since $A f_a=s_af_a$,
    \[
        \operatorname{tr}(c_iAc_jA)=\sum_{a,b}s_as_b\langle f_a,c_if_b\rangle\langle f_b,c_jf_a\rangle=-\sum_{a,b}s_as_b\langle c_if_b,f_a\rangle\langle c_jf_b,f_a\rangle.
    \]
    The terms with $a=b$ vanish because $\langle c_if_a,f_a\rangle=0$. Hence
    \begin{align*}
        \sum_{i,j}B_{ij}\operatorname{tr}(c_iAc_jA)&=-2\sum_{a<b}s_as_b\sum_{i,j}B_{ij}\langle c_if_b,f_a\rangle\langle c_jf_b,f_a\rangle\\
        &=-2\sum_{a<b}s_as_b\operatorname{Ric}(E)(u_{ab},u_{ab})=-2\sum_{a<b}s_as_b r_{ab}.
    \end{align*}
    Thus
    \begin{equation}\label{spin:Ricci-bilinear}
        F=\frac S2+\frac18\sum_{a<b}s_as_b r_{ab}.
    \end{equation}

    We now prove the row sum for $r_{ab}$. Fix $a$. Since
    \[
        \langle u_{ab},e_i\rangle=\langle c_if_b,f_a\rangle=-\langle f_b,c_if_a\rangle,
    \]
    the coordinate matrix of the vectors $u_{ab}$, $b\ne a$, is the transpose, up to an overall sign, of the coordinate matrix of the orthonormal basis $c_1f_a,\cdots,c_7f_a$ of $f_a^\perp$ with respect to the orthonormal basis $f_b$, $b\ne a$. Therefore $\{u_{ab} \mid b\ne a\}$ is an orthonormal basis of $V$. Summing the Ricci quadratic form over this basis gives
    \[
        \sum_{b\ne a}r_{ab}=\sum_{b\ne a}\operatorname{Ric}(E)(u_{ab},u_{ab})=\operatorname{tr}(\operatorname{Ric}(E))=S.
    \]
    Using this and \eqref{spin:Ricci-bilinear},
    \begin{align*}
        \sum_{a<b}(s_a-s_b)^2r_{ab}&=\sum_a s_a^2\sum_{b\ne a}r_{ab}-2\sum_{a<b}s_as_b r_{ab}\\
        &=S\sum_a s_a^2-2\sum_{a<b}s_as_b r_{ab}\\
        &=8S-2\sum_{a<b}s_as_b r_{ab}=16(S-F).
    \end{align*}
    This proves \eqref{spin:Ricci-energy}.

    It remains to prove the row sum for $k_{ab}$. For fixed $a$,
    \[
        \sum_{b\ne a}k_{ab}=\sum_{i<j,\,k<l}E_{ijkl}\sum_{b\ne a}\langle\Gamma_{ij}f_b,f_a\rangle\langle\Gamma_{kl}f_b,f_a\rangle=\sum_{i<j,\,k<l}E_{ijkl}\sum_b\langle f_b,\Gamma_{ij}^*f_a\rangle\langle f_b,\Gamma_{kl}^*f_a\rangle.
    \]
    The omitted $b=a$ term is zero because every $\Gamma_{ij}$ is skew-symmetric. Completeness of the orthonormal basis $f_1,\cdots,f_8$ gives
    \[
        \sum_{b\ne a}k_{ab}=\sum_{i<j,\,k<l}E_{ijkl}\langle\Gamma_{ij}^*f_a,\Gamma_{kl}^*f_a\rangle=-\sum_{i<j,\,k<l}E_{ijkl}\langle f_a,\Gamma_{ij}\Gamma_{kl}f_a\rangle=\frac S2,
    \]
    where the last equality uses the scalar Clifford contraction \eqref{spin:scalar-Clifford}.
\end{proof}

\subsection{A new affine inequality from the original NW lemma}

\begin{lemma}\label{spin:edge-bounds}
    If $E$ has nonnegative sectional curvature, then for all $a\ne b$
    \begin{equation}\label{spin:edge-interval}
        0\le r_{ab}\le\frac S2, \qquad -\frac S6+\frac13r_{ab}\le k_{ab}\le\frac{7S}{6}-\frac73r_{ab}.
    \end{equation}
\end{lemma}

\begin{proof}
    Fix $a\ne b$, and abbreviate $u=u_{ab}$, $\Theta=\Theta_{ab}$, $r=r_{ab}$, $k=k_{ab}$. Let $H=u^\perp$. By Lemma \ref{spin:Theta-geometry}, $\Theta$ is the K\"ahler form of an orthogonal complex structure on $H$, up to an overall sign. Therefore we may choose an orthonormal coframe of $H$ such that $\Theta=e^{12}+e^{34}+e^{56}$. This choice of frame is made separately for the fixed pair $\{a,b\}$; no compatibility between the frames chosen for different pairs is required. Let
    \[
        \kappa_i=K_{2i-1,2i}, \qquad A_{ij}=\sum_{a\in\{2i-1,2i\},\,b\in\{2j-1,2j\}}K_{ab}, \qquad t_{ij}=E_{2i-1,2i,2j-1,2j}.
    \]
    Then
    \[
        k=\sum_{i=1}^3\kappa_i+2\sum_{i<j}t_{ij}, \qquad S_H:=\operatorname{Scal}(E|_H)=2\Bigl(\sum_i\kappa_i+\sum_{i<j}A_{ij}\Bigr).
    \]
    The NW inequality gives $6|t_{ij}|\le A_{ij}+4\sqrt{\kappa_i\kappa_j}$, hence
    \[
        k \ge \sum_i\kappa_i-\frac13\sum_{i<j}A_{ij}-\frac43\sum_{i<j}\sqrt{\kappa_i\kappa_j} \ge -\frac13\Bigl(\sum_i\kappa_i+\sum_{i<j}A_{ij}\Bigr) = -\frac{S_H}{6}.
    \]
    The second step uses $\sum_{i<j}\sqrt{\kappa_i\kappa_j}\le\sum_i\kappa_i$. Using instead the upper bound for each $t_{ij}$ gives
    \[
        k \le \sum_i\kappa_i+\frac13\sum_{i<j}A_{ij}+\frac43\sum_{i<j}\sqrt{\kappa_i\kappa_j} \le \frac73\sum_i\kappa_i+\frac13\sum_{i<j}A_{ij} \le \frac73\Bigl(\sum_i\kappa_i+\sum_{i<j}A_{ij}\Bigr) = \frac76S_H.
    \]
    Since $S_H=S-2\operatorname{Ric}(E)(u,u)=S-2r$, substituting this into the two estimates for $k$ yields
    \[
        k\ge-\frac{S-2r}{6}=-\frac S6+\frac r3, \qquad k\le\frac76(S-2r)=\frac{7S}{6}-\frac73r.
    \]
    Finally, by sectional nonnegativity, $r\ge0$; by $S-2r=S_H\ge0$, we have $r\le S/2$.
\end{proof}

\begin{theorem}\label{spin:affine-theorem}
    For every algebraic curvature tensor with nonnegative sectional curvature in dimension seven and every unit three-form,
    \begin{equation}\label{spin:new-affine}
        Q\ge-\frac{S+F}{3}.
    \end{equation}
    Consequently
    \begin{equation}\label{spin:seventeen-thirtyseven}
        q_3(E)\ge-\frac{17}{37}S\operatorname{Id}.
    \end{equation}
\end{theorem}

\begin{proof}
    By \eqref{spin:Q-energy}, \eqref{spin:edge-interval} and \eqref{spin:Ricci-energy},
    \[
        Q=\frac1{16}\sum_{a<b}(s_a-s_b)^2 k_{ab} \ge \frac1{16}\sum_{a<b}(s_a-s_b)^2\left(-\frac S6+\frac13r_{ab}\right) = -\frac{S}{96}\sum_{a<b}(s_a-s_b)^2+\frac13(S-F).
    \]
    By \eqref{spin:eigenvalues},
    \[
        \sum_{a<b}(s_a-s_b)^2=8\sum_a s_a^2-\Bigl(\sum_a s_a\Bigr)^2=8\cdot8-0=64.
    \]
    Hence $Q\ge-(64/96)S+(S-F)/3=-(S+F)/3$. The dimension-lifting estimate in Theorem \ref{alg:three-lift}, specialized to $n=7$, is
    \begin{equation}\label{spin:old-affine}
        Q\ge\frac{17}{23}(F-S).
    \end{equation}
    Multiplying \eqref{spin:new-affine} by $51/74$ and \eqref{spin:old-affine} by $23/74$, and adding, gives \eqref{spin:seventeen-thirtyseven}.
\end{proof}

\subsection{Four--four spinor splittings and sums of curvature coefficients}

To improve the scalar bound, we need an estimate relating the coefficients $k_{ab}$ for different pairs of indices. We obtain it by considering Clifford matrices with eigenvalues $1$ and $-1$, each of multiplicity four.

\begin{lemma}\label{spin:four-four}
    If a unit three-form satisfies $\mathsf c(\omega)^2=\operatorname{Id}$, then, in some orthonormal coframe of $V$, it has the form
    \begin{equation}\label{spin:Lie-normal-form}
        \omega=\lambda e^{123}+\mu e^{456}, \qquad \lambda^2+\mu^2=1.
    \end{equation}
    In particular $\langle q_3(E)\omega,\omega\rangle\ge0$ whenever $E$ has nonnegative sectional curvature.
\end{lemma}

\begin{proof}
    We first prove the Clifford-square identity
    \begin{equation}\label{spin:Clifford-square}
        \mathsf c(\omega)^2=|\omega|^2\operatorname{Id}-\mathsf c_4\Bigl(\sum_i(\iota_i\omega)\wedge(\iota_i\omega)\Bigr),
    \end{equation}
    where $\mathsf c_4$ denotes degree-four Clifford multiplication, i.e., for a four-form $\beta=\sum_{i<j<k<l}\beta_{ijkl}e^{ijkl}$, write its degree-four Clifford multiplication as
    \[
        \mathsf c_4(\beta):=\sum_{i<j<k<l}\beta_{ijkl}c_ic_jc_kc_l.
    \]
    Write $C_I=c_ic_jc_k$ when $I=\{i<j<k\}$, then
    \[
        \mathsf c(\omega)^2=\Bigl(\sum_I\omega_IC_I\Bigr)\Bigl(\sum_J\omega_JC_J\Bigr)=\sum_{I,J}\omega_I\omega_JC_IC_J.
    \]
    The diagonal terms in $\mathsf c(\omega)^2$ give
    \[
        \sum_I\omega_I^2 C_I^2=\sum_I\omega_I^2\operatorname{Id}=|\omega|^2\operatorname{Id},
    \]
    since each $C_I^2=\operatorname{Id}$. Consider two distinct triples $I$ and $J$. In $\mathsf c(\omega)^2$, their contribution is the sum of the two ordered products
    \[
        \omega_I\omega_J(C_IC_J+C_JC_I).
    \]
    If $|I\cap J|=0,2$, then $C_IC_J+C_JC_I=0$. If $|I\cap J|=1$, then $C_IC_J+C_JC_I$ is a degree-four Clifford monomial. For example, the pair $e^{123}$ and $e^{145}$ gives
    \[
        (c_1c_2c_3)(c_1c_4c_5)+(c_1c_4c_5)(c_1c_2c_3)=-2c_2c_3c_4c_5.
    \]
    On the exterior-form side,
    \[
        \bigl(\iota_1(e^{123}+e^{145})\bigr)\wedge\bigl(\iota_1(e^{123}+e^{145})\bigr)=2e^{2345}.
    \]
    This confirms the minus sign in \eqref{spin:Clifford-square}. The same calculation applies to every pair of triples with exactly one common index, proving the identity.

    The map $\mathsf c_4$ is injective. Indeed, every degree-four Clifford monomial is symmetric and squares to the identity. The product of two distinct degree-four monomials is, up to sign, a nonconstant even Clifford monomial, and has trace zero: conjugating it by a generator appearing in the monomial changes its sign. Consequently, $\operatorname{tr}(C_IC_J)=8\delta_{I,J}$, where $I$ and $J$ are increasing quadruples of indices, and $\operatorname{tr}\bigl(\mathsf c_4(\beta)\mathsf c_4(\gamma)\bigr)=8\langle\beta,\gamma\rangle$. Hence $\mathsf c_4$ is injective. For $|\omega|=1$, identity \eqref{spin:Clifford-square} therefore shows that $\mathsf c(\omega)^2=\operatorname{Id}$ if and only if
    \begin{equation}\label{spin:Jacobi-four-form}
        \sum_i(\iota_i\omega)\wedge(\iota_i\omega)=0.
    \end{equation}

    Define a bracket on $V$ by $\langle[x,y],z\rangle=\omega(x,y,z)$. Clearly the bracket is skew-symmetric, and its inner product is invariant: $\langle[x,y],z\rangle=-\langle y,[x,z]\rangle$. In particular, $\operatorname{ad}_x:y\mapsto[x,y]$ is skew-adjoint for every $x$.

    We now evaluate the four-form in \eqref{spin:Jacobi-four-form}. For any two-form $\alpha$,
    \[
        (\alpha\wedge\alpha)(x,y,z,w)=2\Bigl(\alpha(x,y)\alpha(z,w)-\alpha(x,z)\alpha(y,w)+\alpha(x,w)\alpha(y,z)\Bigr).
    \]
    Apply this to $\alpha=\iota_{e_i}\omega$ and sum over $i$. Since $(\iota_{e_i}\omega)(x,y)=\omega(e_i,x,y)=\langle[x,y],e_i\rangle$, we obtain
    \begin{align*}
        \Bigl(\sum_i(\iota_i\omega)\wedge(\iota_i\omega)\Bigr)(x,y,z,w)&=2\Bigl(\langle[x,y],[z,w]\rangle-\langle[x,z],[y,w]\rangle+\langle[x,w],[y,z]\rangle\Bigr)\\
        &=2\left\langle[[x,y],z]+[[y,z],x]+[[z,x],y],w\right\rangle,
    \end{align*}
    where the last equality uses the invariance of the inner product: for example, $\langle[x,y],[z,w]\rangle=\langle[[x,y],z],w\rangle$. Since this holds for all $w$ and $\mathsf c(\omega)^2=\operatorname{Id}$, by \eqref{spin:Jacobi-four-form}, we have the Jacobi identity $[[x,y],z]+[[y,z],x]+[[z,x],y]=0$. Therefore $(V,[\cdot,\cdot])$ is a Lie algebra with a positive-definite invariant inner product.

    We next use the dimension of this Lie algebra. Recall that the Killing form is defined as the bilinear form $\operatorname{Kill}(x,y) = \operatorname{tr}(\operatorname{ad}_x \circ \operatorname{ad}_y)$. For $B\in\operatorname{End}(V)$, its Hilbert–Schmidt norm is defined by $|B|_{\rm HS}^2=\operatorname{tr}(B^*B)$. Since every $\operatorname{ad}_x$ is skew-adjoint, $\operatorname{Kill}(x,x)=\operatorname{tr}(\operatorname{ad}_x^2)=-|\operatorname{ad}_x|_{\rm HS}^2$. Hence the kernel of the Killing form is precisely the center $\mathfrak z$. Invariance also gives $\mathfrak z^\perp=[V,V]$. Indeed, $x\perp[V,V]$ means $\langle x,[y,z]\rangle=0$ for all $y,z$. By invariance this is equivalent to $\langle[x,y],z\rangle=0$ for all $y,z$, hence to $[x,y]=0$ for every $y$, which says $x\in\mathfrak z$. Thus $[V,V]$ is compact semisimple and is an orthogonal sum of compact simple ideals.

    A compact simple Lie algebra of rank one has two roots, so its dimension is $1+2=3$ and it is isomorphic to $\mathfrak{su}(2)$. A compact simple Lie algebra of rank two has at least six roots and therefore dimension at least $2+6=8$. If its rank is at least three, it has at least twice its rank in roots, so its dimension is at least $3r\ge9$. Since $\dim V=7$, every simple factor of $[V,V]$ must therefore have dimension three, and there can be at most two such factors.

    On a three-dimensional simple factor, every alternating three-form is a scalar multiple of its volume form. Distinct simple factors are orthogonal ideals, and the bracket vanishes between different factors. The center contributes nothing to $\omega$, because $\omega(x,y,z)=\langle[x,y],z\rangle$ and central vectors have zero bracket. Choosing an adapted orthonormal coframe therefore gives
    \[
        \omega=\lambda e^{123}+\mu e^{456},
    \]
    where one coefficient may be zero if there is only one simple factor. The displayed basis forms are orthonormal, so $|\omega|=1$ implies $\lambda^2+\mu^2=1$. This proves \eqref{spin:Lie-normal-form}.

    It remains to prove the curvature assertion. The Ricci term $(\operatorname{Ric}(E)^\sharp)^{[3]}$ changes at most one index of a basis three-form. By the definition of $\mathcal R_E^{(3)}$, its terms apply two interior multiplications and then two wedge multiplications, so they change at most two indices. Since the index sets $\{1,2,3\}$ and $\{4,5,6\}$ are disjoint, neither term can transform $e^{123}$ into $e^{456}$. Thus
    \[
        \langle q_3(E)e^{123},e^{456}\rangle=0.
    \]
    For any distinct indices $i,j,k$, we have
    \[
        \langle q_3(E)e^{ijk},e^{ijk}\rangle=\sum_{p\notin\{i,j,k\}}(K_{ip}+K_{jp}+K_{kp}).
    \]
    In dimension seven, this gives
    \[
        D_{123}:=\langle q_3(E)e^{123},e^{123}\rangle=\sum_{p=4}^7(K_{1p}+K_{2p}+K_{3p})\ge0
    \]
    and
    \[
        D_{456}:=\langle q_3(E)e^{456},e^{456}\rangle=\sum_{p\in\{1,2,3,7\}}(K_{4p}+K_{5p}+K_{6p})\ge0.
    \]
    Therefore $\langle q_3(E)\omega,\omega\rangle=\lambda^2D_{123}+\mu^2D_{456}\ge0$.
\end{proof}

For $A\subset\{1,\cdots,8\}$, denote the sums over $A\times A^c$ by
\[
    k(A)=\sum_{a\in A,\,b\notin A}k_{ab}, \qquad r(A)=\sum_{a\in A,\,b\notin A}r_{ab}.
\]

\begin{corollary}\label{spin:four-cut}
    In every orthonormal spinor basis,
    \begin{equation}\label{spin:cut-bounds}
        |A|=4\ \Longrightarrow\ k(A)\ge0, \qquad |A|=3\ \Longrightarrow\ k(A)\ge-\frac56S.
    \end{equation}
\end{corollary}

\begin{proof}
    First suppose $|A|=4$. Define a symmetric endomorphism $T$ of the spinor module by
    \[
        Tf_a=
        \begin{cases}
            f_a,&a\in A,\\
            -f_a,&a\notin A.
        \end{cases}
    \]
    It is trace-free, and $\operatorname{tr}(T^2)=\operatorname{tr}(\operatorname{Id})=8$. By the isomorphism \eqref{spin:Clifford-map}, there is a three-form $\omega_A$ such that $T=\mathsf c(\omega_A)$. The trace-pairing identity implies $8|\omega_A|^2=\operatorname{tr}(T^2)=8$, so $|\omega_A|=1$. Also $\mathsf c(\omega_A)^2=T^2=\operatorname{Id}$, so Lemma \ref{spin:four-four} gives $\langle q_3(E)\omega_A,\omega_A\rangle\ge0$.

    Across the partition $A\sqcup A^c$, the eigenvalue difference is $s_a-s_b=\pm2$, while it is zero for pairs lying on the same side. The energy identity \eqref{spin:Q-energy} therefore gives
    \[
        \langle q_3(E)\omega_A,\omega_A\rangle=\frac1{16}\sum_{a<b}(s_a-s_b)^2k_{ab}=\frac14k(A).
    \]
    Thus $k(A)\ge0$ whenever $|A|=4$.

    Now suppose $|A|=3$ and let $B=A^c$, so $|B|=5$. For every $j\in B$, the set $A\cup\{j\}$ has four elements, hence $k(A\cup\{j\})\ge0$. Since $(A\cup\{j\})^c=B\setminus\{j\}$,
    \[
        k(A\cup\{j\})=\sum_{a\in A}\sum_{b\in B\setminus\{j\}}k_{ab}+\sum_{b\in B\setminus\{j\}}k_{jb}=k(A)-\sum_{a\in A}k_{aj}+\sum_{b\in B\setminus\{j\}}k_{jb}.
    \]
    Summing over $j\in B$,
    \begin{align*}
        0\le\sum_{j\in B}k(A\cup\{j\})&=\sum_{j\in B}k(A)-\sum_{j\in B}\sum_{a\in A}k_{aj}+\sum_{j\in B}\sum_{b\in B\setminus\{j\}}k_{jb}\\
        &=5k(A)-k(A)+2\sum_{\substack{i<j\\i,j\in B}}k_{ij}=4k(A)+2\sum_{\substack{i<j\\i,j\in B}}k_{ij}.
    \end{align*}

    For each $i\in B$, the row-sum identity \eqref{spin:row-sums} gives $\sum_{b\ne i}k_{ib}=S/2$. Summing over $i\in B$,
    \[
        k(A)+2\sum_{\substack{i<j\\i,j\in B}}k_{ij}=\frac{5S}{2}.
    \]
    Hence
    \[
        0\le4k(A)+2\sum_{\substack{i<j\\i,j\in B}}k_{ij}=4k(A)+\left(\frac{5S}{2}-k(A)\right)=3k(A)+\frac{5S}{2}.
    \]
    Therefore $k(A)\ge-5S/6$ whenever $|A|=3$.
\end{proof}

\subsection{A scalar bound from consecutive eigenvalue differences}

\begin{theorem}\label{spin:rational-improvement}\label{spin:rational}
    In dimension seven every algebraic curvature tensor with nonnegative sectional curvature satisfies
    \begin{equation}\label{spin:rational-bound}
        q_3(E)\ge-\widehat\gamma_7 S\operatorname{Id}, \qquad \widehat\gamma_7=\frac{5168}{11253}=\frac{17}{37}\left(1-\frac{5}{11253}\right).
    \end{equation}
    If $\gamma_7^*$ denotes the optimal scalar coefficient for three-forms on the cone of tensors with nonnegative sectional curvature, then
    \begin{equation}\label{spin:optimal-interval}
        \frac49\le\gamma_7^*\le\frac{5168}{11253}.
    \end{equation}
\end{theorem}

\begin{proof}
    If $S=0$, sectional nonnegativity forces all sectional curvatures to vanish, and polarization gives $E=0$. We may therefore rescale and assume $S=1$. Fix an orthonormal spinor basis $f_1,\cdots,f_8$. Let $K=(k_{ab})$ and $R_\Sigma=(r_{ab})$ be the symmetric $8\times8$ matrices whose off-diagonal entries are given by \eqref{spin:edge-weights} and whose diagonal entries are zero. Their row sums are $1/2$ and $1$, respectively, by \eqref{spin:row-sums}. Thus $L_K=\frac12\operatorname{Id}-K$ and $L_R=\operatorname{Id}-R_\Sigma$ both annihilate $\mathbf1=(1,\cdots,1)^{\mathsf T}$. We will use the identity
    \begin{equation}\label{spin:row-form}
        z^{\mathsf T}(d\operatorname{Id}-B)z=\sum_{a<b}b_{ab}(z_a-z_b)^2,
    \end{equation}
    obtained by expanding the squares whenever $B=(b_{ab})$ is symmetric with zero diagonal and constant row sum $d$.

    For a unit vector $x\in\mathbf1^\perp$, put $s_a=\sqrt8\,x_a$. Since $\sum_a s_a=0$ and $\sum_a s_a^2=8$, \eqref{spin:Clifford-map} and \eqref{spin:trace-pairing} identify $\operatorname{diag}(s_1,\cdots,s_8)$ with a unit three-form. Conversely, every unit three-form occurs in this way in an eigenbasis of its Clifford matrix. Applying \eqref{spin:row-form} to \eqref{spin:Q-energy} and \eqref{spin:Ricci-energy}, we have
    \begin{equation}\label{spin:normalized-energies}
        Q=\frac12 x^{\mathsf T}L_Kx,\qquad 1-F=\frac12 x^{\mathsf T}L_Rx.
    \end{equation}
    The affine estimate \eqref{spin:old-affine} now gives $x^{\mathsf T}Yx\ge0$ for $Y=L_K+\frac{17}{23}L_R$. Since this holds for every unit $x\in\mathbf1^\perp$ and $Y\mathbf1=0$, the symmetric matrix $Y$ is positive semidefinite.

    Define $G=(g_{ab})$ by
    \begin{equation}\label{spin:g-definition}
        g_{ab}=k_{ab}+\frac16-\frac13 r_{ab}\quad(a\ne b),\qquad g_{aa}=0.
    \end{equation}
    Lemma \ref{spin:edge-bounds} says $g_{ab}\ge0$, and \eqref{spin:row-sums} gives $\sum_{b\ne a}g_{ab}=1/2+7/6-1/3=4/3$. Write $L_G=\frac43\operatorname{Id}-G$ and set
    \[
        D=x^{\mathsf T}L_Gx=\sum_{a<b}g_{ab}(x_a-x_b)^2,\qquad H=x^{\mathsf T}Yx.
    \]
    Since $G=K+(\mathbf1\mathbf1^{\mathsf T}-\operatorname{Id})/6-R_\Sigma/3$, we have
    \begin{equation}\label{spin:LG-identity}
        L_G=L_K+\frac43\operatorname{Id}-\frac13 L_R \quad {\rm on} ~ \mathbf1^\perp.
    \end{equation}
    By \eqref{spin:normalized-energies}, $D/2=Q+(1+F)/3$ and $H/2=Q+\frac{17}{23}(1-F)$. Multiplying the first identity by $51/74$, the second by $23/74$, and adding, gives
    \begin{equation}\label{spin:combined-remainder-explicit}
        T:=Q+\frac{17}{37}=\alpha D+\beta H, \qquad \alpha=\frac{51}{148}, \quad \beta=\frac{23}{148}.
    \end{equation}

    We use $D$ to estimate the variation of the coordinates within two groups and $H$ to estimate the difference of their means. Since $x=\sum_{a=1}^8x_af_a$, we can relabel the basis $f_1,\cdots,f_8$ so that $x_1\le\cdots\le x_8$. Put $d_j=x_{j+1}-x_j\ge0$ and $A_j=\{1,\cdots,j\}$. Then for $a<b$, $x_b-x_a=\sum_{j=a}^{b-1}d_j$ and $(x_b-x_a)^2\ge\sum_{j=a}^{b-1}d_j^2$. For $A\subset\{1,\cdots,8\}$, write $g(A)=\sum_{a\in A,\,b\notin A}g_{ab}$. Then
    \begin{equation}\label{spin:ordered-gaps}
        D=\sum_{a<b}g_{ab}(x_a-x_b)^2\ge\sum_{a<b}g_{ab}\sum_{j=a}^{b-1}d_j^2=\sum_{j=1}^7 g(A_j)d_j^2.
    \end{equation}
    To estimate these coefficients, let $m=|A|$. By \eqref{spin:g-definition}, $g(A)=k(A)+m(8-m)/6-r(A)/3$. Since $r_{ab}\ge0$ and $\sum_{b\ne a}r_{ab}=1$, we have
    \[
        r(A) = \sum_{a\in A}\sum_{b\notin A}r_{ab} \le \sum_{a\in A}\sum_{b\ne a}r_{ab} = \sum_{a\in A}1 = m.
    \]
    For $m=1$ the row sum of $G$ gives $g(A)=4/3$. For $m=3,4$, Corollary \ref{spin:four-cut} gives, respectively,
    \[
        g(A)\ge-\frac56+\frac{15}6-1=\frac23, \qquad g(A)\ge\frac{16}6-\frac43=\frac43.
    \]
    Since $g(A)=g(A^c)$, this proves $g(A_j)\ge2/3$ for $j=1,3,4,5,7$.

    For the two remaining positions, consider the disjoint pairs $P_1=\{1,2\}$ and $P_2=\{7,8\}$. Their union has four elements, so
    \[
        g(P_1)+g(P_2)=g(P_1\cup P_2)+2\sum_{a\in P_1,\,b\in P_2}g_{ab}\ge\frac43.
    \]
    Thus at least one of $g(A_2)=g(P_1)$ and $g(A_6)=g(P_2)$ is at least $2/3$. Let $j_0\in\{2,6\}$ correspond to the smaller value, choosing either in case of equality. This splits the ordered indices into two groups, $A_{j_0}$ and $A_{j_0}^c$. Also $\{|A_{j_0}|,|A_{j_0}^c|\}=\{2,6\}$. Let $P$ denote the smaller group, so $|P|=2$. By \eqref{spin:ordered-gaps}, $D\ge\frac23\sum_{j\ne j_0}d_j^2=\frac23\sum_{j\ne j_0}(x_{j+1}-x_j)^2$. The omitted term $d_{j_0}=x_{j_0+1}-x_{j_0}$ is the unique consecutive difference that crosses the partition $P\sqcup P^c$. Thus the sum $\sum_{j\ne j_0}(x_{j+1}-x_j)^2$ contains precisely the consecutive differences within $P$ and within $P^c$.

    Write $\mathbf1_P$ for the indicator vector of $P$, and decompose $x$ orthogonally as
    \[
        x=\rho v+w, \qquad v=\frac{4\mathbf1_P-\mathbf1}{\sqrt{24}}, \qquad \rho=\langle x,v\rangle, \qquad \rho^2+|w|^2=1.
    \]
    The vectors $\mathbf1$ and $v$ span the subspace of vectors constant on each of $P$ and $P^c$. If $C=P$ or $P^c$, then writing $\bar x_C=\sum_{b\in C}x_b/|C|$, we have $w_a=x_a-\bar x_C$ for any $a\in C$, since $x\perp\mathbf1$. For $t_1\le\cdots\le t_m$, let $\bar t$ be their mean, then $|t_a-\bar t|\le t_m-t_1$ for any $a$. By Cauchy--Schwarz, $(t_a-\bar t)^2\le(t_m-t_1)^2\le(m-1)\sum_{i=1}^{m-1}(t_{i+1}-t_i)^2$. Applying this and using $2(2-1),6(6-1)\le30$,
    \begin{equation}\label{spin:within-subsets}
        |w|^2=\sum_{a\in P}(x_a-\bar x_P)^2+\sum_{a\in P^c}(x_a-\bar x_{P^c})^2\le30\sum_{j\ne j_0}d_j^2\le45D.
    \end{equation}

    It remains to control the direction $v$. Write $P=\{p,q\}$. By \eqref{spin:row-sums}, $k(P)=1-2k_{pq}$ and $r(P)=2-2r_{pq}$. Across $P$ and $P^c$ the squared coordinate difference of $v$ is $2/3$, and within either set it is zero. Thus \eqref{spin:row-form} and the upper bound $k_{pq}\le7/6-7r_{pq}/3$ from Lemma \ref{spin:edge-bounds} give
    \[
        v^{\mathsf T}Yv = \sum_{a<b}\left(k_{ab}+\frac{17}{23}r_{ab}\right)(v_a-v_b)^2 = \frac23\left(k(P)+\frac{17}{23}r(P)\right) \ge \frac{20+440r_{pq}}{207}\ge\frac{20}{207} =: y.
    \]
    For an upper bound, if $B$ in \eqref{spin:row-form} has nonnegative entries, then by $(z_a-z_b)^2\le2(z_a^2+z_b^2)$,
    \[
        z^{\mathsf T}(d\operatorname{Id}-B)z \le 2\sum_{a<b}b_{ab}(z_a^2+z_b^2) = 2\sum_a z_a^2\sum_{b\ne a}b_{ab} = 2d|z|^2.
    \]
    Hence $L_G\le\frac83\operatorname{Id}$ and $L_R\le2\operatorname{Id}$. On $\mathbf1^\perp$, the definition of $Y$ and \eqref{spin:LG-identity} give
    \[
        Y=L_G+\frac{74}{69}L_R-\frac43\operatorname{Id}\le M\operatorname{Id}, \qquad M=\frac83+\frac{148}{69}-\frac43=\frac{80}{23}.
    \]

    Write $|z|_Y=(z^{\mathsf T}Yz)^{1/2}$. Since $v^{\mathsf T}Yv\ge y$ and $0\le Y\le M\operatorname{Id}$ on $\mathbf1^\perp$, the triangle inequality for this seminorm and \eqref{spin:within-subsets} give
    \begin{align*}
        \sqrt{y}\,|\rho| &\le |\rho|(v^{\mathsf T}Yv)^{1/2} = |\rho v|_Y = |x-w|_Y \le |x|_Y+|w|_Y \\
        &= (x^{\mathsf T}Yx)^{1/2}+(w^{\mathsf T}Yw)^{1/2} \le \sqrt{H}+\sqrt{M}|w| \le \sqrt{H}+\sqrt{45MD}.
    \end{align*}
    Then using Cauchy--Schwarz and \eqref{spin:combined-remainder-explicit}, we obtain
    \[
        \rho^2 \le \frac1y(\beta H+\alpha D)\left(\frac1\beta+\frac{45M}{\alpha}\right) = \frac{T}{y}\left(\frac1\beta+\frac{45M}{\alpha}\right), \qquad |w|^2 \le \frac{45}{\alpha}T.
    \]
    Adding the two bounds and using $\rho^2+|w|^2=1$ yields
    \begin{equation}\label{spin:explicit-gap}
        T \ge \left[\frac1{\beta y}+\frac{45}{\alpha}\left(1+\frac My\right)\right]^{-1} = \frac{85}{416361}.
    \end{equation}
    Since $T=Q+17/37$, this gives $Q\ge 85/416361-17/37=-5168/11253$ for every unit three-form when $S=1$. By homogeneity in $E$ and in the form, \eqref{spin:rational-bound} follows, and hence $\gamma_7^*\le\widehat\gamma_7$.

    For the lower bound in \eqref{spin:optimal-interval}, decompose $V=W^6\oplus\mathbb R$, choose an orthogonal complex structure $J$ on $W$, and extend $E=tE_J$ by zero for $t>0$. Let $\zeta$ be a unit covector on the flat summand. By Lemmas \ref{lem:model-operator} and \ref{lem:model-powers}, $S=18t$, $|\Omega_J|^2=3$, and $q_3(E)(\zeta\wedge\Omega_J)=-8t\,\zeta\wedge\Omega_J$. The unit form $\zeta\wedge\Omega_J/\sqrt3$ thus has $-Q/S=8t/(18t)=4/9$, proving $\gamma_7^*\ge4/9$.
\end{proof}

\section{Sharp estimates under strong nonnegative curvature}
\label{strong:section}

We now impose an additional algebraic hypothesis. It permits an exact answer in all degrees, and supplies comparison models for the full sectionally nonnegative problem.

\subsection{The Bianchi projection and rank-one generators}

A symmetric operator on $\Lambda^2V^*$ can be identified with a four-tensor having the skew symmetries in each pair and symmetry under interchanging the pairs, but not necessarily satisfying the first Bianchi identity. Since a four-form $\eta\in\Lambda^4V^*$ has these symmetry properties, it determines a symmetric operator $T_\eta$ on $\Lambda^2V^*$ by
\[
    \langle T_\eta(X^\flat\wedge Y^\flat),Z^\flat\wedge W^\flat\rangle=\eta(X,Y,Z,W).
\]
In particular, $\langle T_\eta\xi,\xi\rangle=0$ whenever $\xi$ is decomposable. By an abuse of notation, we write $T_\eta=\eta$.

\begin{definition}\label{strong:definition}
    An algebraic curvature tensor $E$ has \emph{strongly nonnegative curvature} if there exists $\eta\in\Lambda^4V^*$ such that $\mathcal R_E+\eta$ is positive semidefinite on $\Lambda^2V^*$.
\end{definition}

This condition implies $K_E\geq0$, since
\[
    E(X,Y,X,Y)=\langle\mathcal R_E(X^\flat\wedge Y^\flat),X^\flat\wedge Y^\flat\rangle=\langle(\mathcal R_E+\eta)(X^\flat\wedge Y^\flat),X^\flat\wedge Y^\flat\rangle\ge0.
\]
In contrast, nonnegative sectional curvature only guarantees nonnegativity on decomposable bivectors and does not guarantee positive semidefiniteness on all of $\Lambda^2V^*$; hence strongly nonnegative curvature is a stronger condition. The two cones coincide in dimension four and differ in every dimension $n\ge5$.

\begin{lemma}[Rank-one Bianchi decomposition]\label{strong:decomposition}
    An algebraic curvature tensor $E$ is strongly nonnegative if and only if it can be written as a finite sum
    \begin{equation}\label{strong:rank-one-sum}
        E=\frac13\sum_{\nu=1}^N E_{A_\nu},
    \end{equation}
    where each $A_\nu$ is skew-adjoint and
    \begin{align}
        E_A(X,Y,Z,W)&=\Omega_A(X,Z)\Omega_A(Y,W)-\Omega_A(X,W)\Omega_A(Y,Z) \notag\\
        &\quad+2\Omega_A(X,Y)\Omega_A(Z,W), \qquad \Omega_A(X,Y)=\langle AX,Y\rangle. \label{strong:EA}
    \end{align}
\end{lemma}

\begin{proof}
    For a symmetric operator $B$ on $\Lambda^2V^*$, write
    \[
        b(B)_{ijkl}=\frac13\bigl(B_{ijkl}+B_{jkil}+B_{kijl}\bigr).
    \]
    Direct verification shows that $b(B)$ is a four-form. Conversely, if $B$ itself comes from a four-form, then $B_{ijkl}=B_{jkil}=B_{kijl}$, so $b(B)=B$. Therefore $\pi(B):=B-b(B)$ satisfies the first Bianchi identity and is the projection onto algebraic curvature tensors. In particular, $\pi(\eta)=0$ for every four-form $\eta$ and $\pi(\mathcal R_E)=\mathcal R_E$.

    For a two-form $\Omega=\Omega_A$, the rank-one operator $B=\Omega\otimes\Omega$ has entries $B_{ijkl}=\Omega_{ij}\Omega_{kl}$. The projection is
    \begin{equation}\label{strong:projection}
        \pi(B)_{ijkl} = \frac13\bigl(2\Omega_{ij}\Omega_{kl}-\Omega_{jk}\Omega_{il}-\Omega_{ki}\Omega_{jl}\bigr) = \frac13\bigl(\Omega_{ik}\Omega_{jl}-\Omega_{il}\Omega_{jk}+2\Omega_{ij}\Omega_{kl}\bigr) = \frac13(E_A)_{ijkl}.
    \end{equation}
    This also verifies the Bianchi identity for \eqref{strong:EA}.

    If $\mathcal R_E+\eta\geq0$, then by the orthogonal diagonalization of real symmetric matrices, $\mathcal R_E+\eta=\sum_\nu\Omega_\nu\otimes\Omega_\nu^\flat$, i.e., $E+\eta=\sum_\nu\Omega_\nu\otimes\Omega_\nu$, where the square roots of the nonnegative eigenvalues have been absorbed into the eigenforms, and $\Omega_\nu^\flat\in(\Lambda^2V^*)^*=\Lambda^2V$. For each $\Omega_\nu$, there exists a unique skew-adjoint operator $A_\nu$ such that $\Omega_{A_\nu}=\Omega_\nu$. Applying $\pi$ to the decomposition and using \eqref{strong:projection} proves \eqref{strong:rank-one-sum}.

    Conversely, since $(\Omega_A\wedge\Omega_A)_{ijkl}=2\bigl(\Omega_{ij}\Omega_{kl}-\Omega_{ik}\Omega_{jl}+\Omega_{il}\Omega_{jk}\bigr)$, we have $(E_A)_{ijkl}+(\Omega_A\wedge\Omega_A)_{ijkl}/2=3\Omega_{ij}\Omega_{kl}$ by \eqref{strong:EA}. Hence $E_A+\frac12\Omega_A\wedge\Omega_A=3\Omega_A\otimes\Omega_A$, or equivalently,
    \begin{equation}\label{strong:EA-positive-modification}
        \mathcal R_{E_A}+\frac12T_{\Omega_A\wedge\Omega_A}=3\Omega_A\otimes\Omega_A^\flat\geq0.
    \end{equation}
    Using this and \eqref{strong:rank-one-sum}, $\mathcal R_E+\frac16\sum_{\nu=1}^NT_{\Omega_{A_\nu}\wedge\Omega_{A_\nu}}=\sum_{\nu=1}^N\Omega_{A_\nu}\otimes\Omega_{A_\nu}^\flat \ge 0$.
\end{proof}

\subsection{Exterior-power blocks and two Cauchy--Schwarz inequalities}

\begin{theorem}[Sharp bound on the strongly nonnegative cone]\label{strong:sharp}
    Let $n\geq4$, $0\leq p\leq n$, and put
    \[
        a=\left\lfloor\frac p2\right\rfloor, \qquad b=\left\lfloor\frac{n-p}{2}\right\rfloor, \qquad C^{\mathrm{str}}_{n,p}=\frac{2ab}{3(a+b)}.
    \]
    For every strongly nonnegative algebraic curvature tensor $E$,
    \begin{equation}\label{strong:sharp-equation}
        q_p(E)\geq-C^{\mathrm{str}}_{n,p}\operatorname{Scal}(E)\operatorname{Id}_{\Lambda^p}.
    \end{equation}
    The coefficient is sharp for every $n$ and $p$ in this range.
\end{theorem}

\begin{proof}
    By Lemma \ref{strong:decomposition} and linearity it suffices to consider $E_A$. Let $m=\lfloor n/2\rfloor$ and $\nu=n-2m\in\{0,1\}$. Choose orthonormal two-planes $P_i$ and orthonormal bases $\{e_{2i-1},e_{2i}\}$ of $P_i$ such that
    \[
        A|_{P_i}=t_iJ_i, \qquad t_i\ge0, \qquad J_i e_{2i-1}=e_{2i}, \quad J_i e_{2i}=-e_{2i-1}.
    \]
    If $\nu=1$, denote the remaining unit covector by $z$. Set $\theta_i=e^{2i-1,2i}$, so $\Omega_A=\sum_i t_i\theta_i$. On each $P_i$, $A^2=-t_i^2\operatorname{Id}$; hence $\operatorname{tr} A^2=-2\sum_i t_i^2$ and $\operatorname{Scal}(E_A)=-3\operatorname{tr} A^2=6\sum_i t_i^2$.

    Recall that if $T$ is real and self-adjoint, then its complexification is Hermitian self-adjoint. Moreover, if $u=x+iy$ with $x,y$ real, then $\langle T_{\mathbb C}u,u\rangle_{\mathbb C}=\langle Tx,x\rangle+\langle Ty,y\rangle$. Thus a real quadratic-form inequality for $T$ is equivalent to the corresponding Hermitian inequality for its complexification. Therefore, we may complexify the exterior algebra: $(\Lambda^pV^*)_{\mathbb C}=(\Lambda^pV^*)\otimes_{\mathbb R}\mathbb C$. On each $P_i$, consider the following orthonormal basis of $(\Lambda^*P_i^*)_{\mathbb C}$:
    \[
        1, \qquad u_i^+=\frac{e^{2i-1}+\sqrt{-1}e^{2i}}{\sqrt2}, \qquad u_i^-=\frac{e^{2i-1}-\sqrt{-1}e^{2i}}{\sqrt2}, \qquad \theta_i.
    \]
    They satisfy $\deg(1)=0$, $\deg(u_i^\pm)=1$, $\deg(\theta_i)=2$, and
    \[
        D_A(1)=0, \qquad D_Au_i^+=\sqrt{-1}t_i u_i^+, \qquad D_Au_i^-=-\sqrt{-1}t_i u_i^-, \qquad D_A(\theta_i)=0.
    \]
    If $\nu=1$, then $D_Az=0$.

    Fix a set $J\subset\{1,\cdots,m\}$ of $s$ two-planes from which we choose one-form basis vectors. For every $j\in J$, choose a sign $\sigma_j\in\{+1,-1\}$, and use the one-form $u_j^{\sigma_j}$. If $\nu=1$, choose whether $z$ occurs. Denote this choice by $\epsilon\in\{0,1\}$, so $\epsilon\le\nu$. If $\nu=0$, then $\epsilon=0$. Among the remaining $d=m-s$ two-planes, choose exactly $r$ planes to contribute their full two-form $\theta_i$, where $r$ is determined by the degree condition
    \begin{equation}\label{strong:degree-count}
        p=2r+s+\epsilon.
    \end{equation}
    Let $I\subset\{1,\cdots,m\}\setminus J$ with $|I|=r$. Define
    \[
        f_I=\Bigl(\bigwedge_{j\in J}u_j^{\sigma_j}\Bigr)\wedge\Bigl(\bigwedge_{i\in I}\theta_i\Bigr)\wedge z^\epsilon.
    \]
    We order the factors in a fixed way. Since every $\theta_i$ has even degree, permuting the $\theta_i$ does not introduce signs. The vectors $f_I$ form an orthonormal basis of the subspace
    \[
        \mathcal H_{J,\sigma,\epsilon,r}:=\operatorname{span}_{\mathbb C}\left\{f_I \mid I\subset\{1,\cdots,m\}\setminus J,\ |I|=r\right\} \subset (\Lambda^pV^*)_{\mathbb C}.
    \]
    Moreover, $(\Lambda^pV^*)_{\mathbb C}$ can be decomposed as the orthogonal direct sum of all such subspaces:
    \begin{equation}\label{strong:space-decomposition}
        (\Lambda^pV^*)_{\mathbb C}=\bigoplus_{J,\sigma,\epsilon,r}\mathcal H_{J,\sigma,\epsilon,r}.
    \end{equation}
    We claim that the three operators in \eqref{eq:model-operator} preserve every $\mathcal H_{J,\sigma,\epsilon,r}$.

    We compute those terms separately. First, the Leibniz rule gives
    \[
        -2(A^2)^{[p]}f_I=\Bigl(4\sum_{i\in I}t_i^2+2\sum_{j\in J}t_j^2\Bigr)f_I.
    \]
    Second, from $D_A\theta_i=0$ and $D_Au_j^{\sigma_j}=\sqrt{-1}\sigma_jt_ju_j^{\sigma_j}$, we have $D_Af_I=\sqrt{-1}(\sum_{j\in J}\sigma_jt_j)f_I$ and
    \[
        -D_A^2f_I=\Bigl(\sum_{j\in J}\sigma_jt_j\Bigr)^2 f_I.
    \]
    Third, since $(\theta_i\wedge)^*=(\varepsilon_{2i-1}\varepsilon_{2i})^*=\iota_{2i}\iota_{2i-1}$,
    \[
        \Lambda_A=\Bigl(\sum_{i=1}^m t_i\theta_i\wedge\Bigr)^*=\sum_{i=1}^m t_i\iota_{2i}\iota_{2i-1}.
    \]
    If $i\notin I$, then $f_I$ has no $\theta_i$ factor, so $\iota_{2i}\iota_{2i-1}f_I=0$. If $i\in I$, then $f_I$ has a $\theta_i$ factor, so $\iota_{2i}\iota_{2i-1}$ exactly removes this factor, yielding $f_{I\setminus\{i\}}$. Hence $\Lambda_A f_I=\sum_{i\in I}t_i f_{I\setminus\{i\}}$. Then
    \[
        L_A\Lambda_Af_I = \sum_{i\in I}t_i\sum_{k=1}^m t_k\theta_k\wedge f_{I\setminus\{i\}} = \sum_{i\in I}t_i^2 f_I+\sum_{i\in I}\sum_{k\notin I\cup J}t_it_kf_{(I\setminus\{i\})\cup\{k\}}.
    \]
    If $k\notin I\cup J$, then $(I\setminus\{i\})\cup\{k\}\subset\{1,\cdots,m\}\setminus J$ and has exactly $r$ elements. Hence $L_A\Lambda_Af_I\in\mathcal H_{J,\sigma,\epsilon,r}$. This proves the claim.

    From the claim and the decomposition \eqref{strong:space-decomposition}, if $q_p(E_A)\big|_{\mathcal H}\ge-C\operatorname{Scal}(E_A)\operatorname{Id}_{\mathcal H}$ for all $\mathcal H$, then for any $\omega=\sum\omega_{\mathcal H}$,
    \[
        \langle q_p(E_A)\omega,\omega\rangle_{\mathbb C} = \sum_{\mathcal H}\langle q_p(E_A)\omega_{\mathcal H},\omega_{\mathcal H}\rangle_{\mathbb C} \ge -C\operatorname{Scal}(E_A)\sum_{\mathcal H}|\omega_{\mathcal H}|^2 = -C\operatorname{Scal}(E_A)|\omega|^2.
    \]
    Therefore, it suffices to prove \eqref{strong:sharp-equation} on a fixed $\mathcal H_{J,\sigma,\epsilon,r}$.

    Define $A_r$ on $\mathcal H_{J,\sigma,\epsilon,r}$ by
    \[
        A_rf_I:=\sum_{i\in I}\sum_{k\notin I\cup J}t_it_kf_{(I\setminus\{i\})\cup\{k\}}.
    \]
    Then $L_A\Lambda_Af_I=\sum_{i\in I}t_i^2f_I+A_rf_I$. By \eqref{eq:model-operator},
    \begin{equation}\label{strong:occupancy-operator}
        q_p(E_A)\big|_{\mathcal H_{J,\sigma,\epsilon,r}}=\Bigl[2\sum_{j\in J}t_j^2+\Bigl(\sum_{j\in J}\sigma_jt_j\Bigr)^2\Bigr]\operatorname{Id}-4A_r.
    \end{equation}
    If $r=0$, then $I=\varnothing$. If $r=d$, then $I=\{1,\cdots,m\}\setminus J$ and $I\cup J=\{1,\cdots,m\}$. Both cases imply $A_r=0$, and \eqref{strong:occupancy-operator} is nonnegative. It remains to consider $0<r<d$.

    Note that $A_r$ is self-adjoint. Indeed, let $I'\subset\{1,\cdots,m\}\setminus J$ with $|I'|=r$, then
    \[
        \langle A_r f_I,f_{I'}\rangle_{\mathbb C}=
        \begin{cases}
            t_it_k, & I'=(I\setminus\{i\})\cup\{k\},\quad i\in I,\ k\notin I\cup J,\\
            0, & {\rm otherwise}.
        \end{cases}
    \]
    Since $I'=(I\setminus\{i\})\cup\{k\}$ is equivalent to $I=(I'\setminus\{k\})\cup\{i\}$, $A_r$ is self-adjoint. Let $f\in\mathcal H_{J,\sigma,\epsilon,r}$ be a unit vector and write
    \[
        f=\sum_{\substack{I\subset\{1,\cdots,m\}\setminus J\\ |I|=r}}a_If_I,
    \]
    where $a_I\in\mathbb C$. Define
    \[
        T=\sum_{i\in\{1,\cdots,m\}\setminus J}t_i^2, \qquad X=\sum_{\substack{I\subset\{1,\cdots,m\}\setminus J\\ |I|=r}}|a_I|^2\sum_{i\in I}t_i^2.
    \]
    We now compute the quadratic form of $A_r$.

    Group the terms according to the common $(r-1)$-element subset $K$. For $K\subset\{1,\cdots,m\}\setminus J$ with $|K|=r-1$, the possible endpoints are $K\cup\{i\}$, where $i\notin J\cup K$. Expanding the square gives
    \begin{equation}\label{strong:subset-energy}
        \langle A_rf,f\rangle_{\mathbb C} = \sum_{\substack{K\subset\{1,\cdots,m\}\setminus J\\ |K|=r-1}}\left(\left|\sum_{i\notin J\cup K}t_ia_{K\cup\{i\}}\right|^2-\sum_{i\notin J\cup K}t_i^2|a_{K\cup\{i\}}|^2\right).
    \end{equation}
    Since $|\{1,\cdots,m\}\setminus(J\cup K)|=d-r+1$, the Cauchy--Schwarz inequality gives
    \[
        \left|\sum_{i\notin J\cup K}t_ia_{K\cup\{i\}}\right|^2 \leq (d-r+1)\sum_{i\notin J\cup K}t_i^2|a_{K\cup\{i\}}|^2.
    \]
    Substituting this into \eqref{strong:subset-energy} yields
    \begin{equation}\label{strong:first-CS}
        \langle A_rf,f\rangle_{\mathbb C} \leq (d-r)\sum_{\substack{K\subset\{1,\cdots,m\}\setminus J\\ |K|=r-1}}\sum_{i\notin J\cup K}t_i^2|a_{K\cup\{i\}}|^2 = (d-r)X .
    \end{equation}

    We now apply the same estimate to complementary subsets. For $I\subset\{1,\cdots,m\}\setminus J$ with $|I|=r$, let
    \[
        I^c:=\{1,\cdots,m\}\setminus(I\cup J), \qquad b_{I^c}:=a_I.
    \]
    The map $I\mapsto I^c$ is a bijection between $r$-element and $(d-r)$-element subsets. Under this relabeling, the weighted replacement operator $A_r$ becomes the corresponding weighted replacement operator on $(d-r)$-element subsets. Applying \eqref{strong:first-CS} with $r$ replaced by $d-r$ gives
    \[
        \langle A_rf,f\rangle_{\mathbb C} \leq (d-(d-r))\sum_{\substack{I^c\subset\{1,\cdots,m\}\setminus J\\ |I^c|=d-r}}|b_{I^c}|^2\sum_{i\in I^c}t_i^2 = r\sum_{\substack{I\subset\{1,\cdots,m\}\setminus J\\ |I|=r}}|a_I|^2\sum_{i\notin I\cup J}t_i^2.
    \]
    Since $\sum_{i\notin I\cup J}t_i^2=T-\sum_{i\in I}t_i^2$, we obtain
    \begin{equation}\label{strong:second-CS}
        \langle A_rf,f\rangle_{\mathbb C}\leq r(T-X).
    \end{equation}
    Multiplying \eqref{strong:first-CS} by $r/d$ and \eqref{strong:second-CS} by $(d-r)/d$, we obtain $\langle A_rf,f\rangle_{\mathbb C}\leq\frac{r(d-r)}dT$.

    By \eqref{strong:degree-count}, $r\leq\lfloor p/2\rfloor=a$. On the other hand,
    \[
        n-p=2m+\nu-(2r+s+\epsilon)=2(d-r)+s+(\nu-\epsilon),
    \]
    so $d-r\leq\lfloor(n-p)/2\rfloor=b$. Since $\frac{xy}{x+y}$ is nondecreasing in each variable whenever $x+y>0$,
    \[
        \frac{r(d-r)}d=\frac{r(d-r)}{r+(d-r)}\leq\frac{ab}{a+b}, \qquad T\leq\sum_{i=1}^m t_i^2.
    \]
    Discarding the nonnegative scalar term in \eqref{strong:occupancy-operator}, we conclude that
    \[
        q_p(E_A)\big|_{\mathcal H_{J,\sigma,\epsilon,r}} \geq -\frac{4ab}{a+b}\left(\sum_{i=1}^m t_i^2\right)\operatorname{Id} = -\frac{2ab}{3(a+b)}\operatorname{Scal}(E_A)\operatorname{Id}.
    \]

    To prove sharpness, assume first that $a,b\geq1$ and put $m_0=a+b$. Let $W\subset V$ be a $2m_0$-dimensional subspace carrying an orthogonal complex structure $J$. Extend $E_J$ by zero to $V$, then $E_J$ is indeed strongly nonnegative by \eqref{strong:EA-positive-modification}. The remaining dimension is
    \[
        n-2m_0=
        \begin{cases}
            0,&n,p\text{ both even},\\
            1,&n\text{ odd},\\
            2,&n\text{ even and }p\text{ odd}.
        \end{cases}
    \]
    Use
    \[
        \omega=\frac{\Omega_J^a}{a!\sqrt{\binom{m_0}{a}}} ~ \text{if $p$ is even}; \qquad \omega=\frac{\Omega_J^a\wedge\zeta}{a!\sqrt{\binom{m_0}{a}}} ~ \text{if $p$ is odd},
    \]
    where $\zeta$ is a unit covector in the remaining space. In either case, $\omega$ is a $p$-form. By Lemmas \ref{lem:model-operator} and \ref{lem:model-powers}, $|\omega|=1$, $\operatorname{Scal}(E_J)=6m_0$ and $q_p(E_J)\omega=-4a(m_0-a)\omega=-4ab\omega$. This proves sharpness for the interior cases.

    If $a=0$ or $b=0$, then $C^{\mathrm{str}}_{n,p}=0$ and $p\in\{0,1,n-1,n\}$. For $p=0$ or $p=n$, $q_p\equiv0$. For $p=1$, $q_1=(\mathrm{Ric}^\sharp)^{[1]}$. In this case, take a two-dimensional plane $P\subset V$ and define
    \[
        E(X,Y,Z,W)=
        \begin{cases}
            \langle X,Z\rangle\langle Y,W\rangle-\langle X,W\rangle\langle Y,Z\rangle,&X,Y,Z,W\in P,\\
            0,&\text{if any one of }X,Y,Z,W\text{ lies in }P^\perp.
        \end{cases}
    \]
    This $E$ is strongly nonnegative since $\mathcal R_E\geq0$. Because $n\ge4$, we have $\dim P^\perp=n-2\ge2$, so there is a nonzero unit vector $u\in P^\perp$. Taking $\omega=u^\flat$, we obtain $\langle q_1(E)\omega,\omega\rangle=\mathrm{Ric}(E)(u,u)=0$. The case of $p=n-1$ follows by Lemma \ref{lem:HodgeDual}. Thus zero is also sharp in these cases.
\end{proof}

For geometric applications based on $E=\operatorname{Rm}-K_{\min}I$, the extra assumption in Theorem \ref{strong:sharp} must be imposed on this tensor $E$. Strong nonnegativity of $\operatorname{Rm}$ alone is not sufficient to imply strong nonnegativity of $E=\operatorname{Rm}-K_{\min}I$.

\section{Algebraic equality cases}\label{eqcase:section}

In the previous theorems, it was often proved that a certain estimate attains equality by constructing a model tensor $E_J$ and a form $\omega$, thereby showing that the coefficient is sharp. However, the existence of one example for which equality holds does not mean that all examples attaining equality have been found. There may be other tensors and forms that also attain equality. We write $S=\operatorname{Scal}(E)$, $Q_p=\langle q_p(E)\omega,\omega\rangle$ and $F_p=\langle(\operatorname{Ric}(E)^\sharp)^{[p]}\omega,\omega\rangle$. Unless otherwise stated, $|\omega|=1$.

\subsection{Null directions and a four-dimensional determination}

\begin{lemma}\label{eqcase:nullity}\label{eqcase:ricci-kernel}
    Suppose that $E$ has nonnegative sectional curvature. Then
    \begin{equation}\label{eqcase:ricci-null}
        \operatorname{Ric}(E)(v,v)=0 \quad \Longrightarrow \quad E(v,x,y,z)=0 \quad \text{for all }x,y,z\in V.
    \end{equation}
    In particular, $S=0$ implies $E=0$.
\end{lemma}

\begin{proof}
    If $\operatorname{Ric}(E)(v,v)=0$, then $E(v,x,v,x)=0$ for every $x$. For arbitrary $u,x$, consider
    \[
        E(v+tu,x,v+tu,x)=2tE(v,x,u,x)+t^2E(u,x,u,x).
    \]
    Nonnegativity for both signs of sufficiently small $t$ gives $E(v,x,u,x)=0$. Define $T(x,u,y)=E(v,x,u,y)$. Polarizing the identity $E(v,x,u,x)=0$ in $x$ gives $T(x,u,y)=-T(y,u,x)$. Also, $T(x,u,y)=-T(x,y,u)$. These two transpositions generate all permutations, so $T$ is a three-form. Apply the first Bianchi identity and use pair symmetry:
    \[
        0=E(v,x,u,y)+E(x,u,v,y)+E(u,v,x,y)=T(x,u,y)+T(y,x,u)-T(u,x,y)=3T(x,u,y).
    \]
    
    Finally, if $S=0$, then $\operatorname{Ric}(E)(v,v)=0$ for all $v$. Hence $E(v,x,y,z)=0$ for all $v,x,y,z$.
\end{proof}

\begin{lemma}\label{eqcase:four-dimensional}
    Let $E$ be a sectional-nonnegative algebraic curvature tensor on an orthonormal four-frame. Suppose that
    \[
        K_{12}=k>0, \quad K_{34}=l>0, \quad K_{13}=K_{14}=K_{23}=K_{24}=0, \qquad E_{1234}=\frac23\epsilon\sqrt{kl},
    \]
    where $\epsilon\in\{1,-1\}$. Then $E=E_A$, where $A$ is determined by
    \[
        \Omega_A=\sqrt{k/3}\,e^{12}+\epsilon\sqrt{l/3}\,e^{34}.
    \]
    If $k=l$, then $E=(k/3)E_J$, with $\Omega_J=e^{12}+\epsilon e^{34}$.
\end{lemma}

\begin{proof}
    Reverse $e_4$ if necessary to assume $\epsilon=1$ without loss of generality. Define
    \[
        T:(e_1,e_2,e_3,e_4)\mapsto(se_1,se_2,s^{-1}e_3,s^{-1}e_4), \quad (T^*E)(X,Y,Z,W)=E(TX,TY,TZ,TW)
    \]
    with $s^8=l/k$. Then \cite[Lemma 4.2(2)]{Ge} gives $T^*E=(\sqrt{kl}/3)E_J$, where $J$ is determined by $\Omega_J=e^{12}+e^{34}$. Hence $E=(\sqrt{kl}/3)(T^{-1})^*E_J=E_A$ with $\Omega_A=\sqrt{k/3}\,e^{12}+\sqrt{l/3}\,e^{34}$.
\end{proof}

\subsection{Complete equality classification in dimensions five and six}

\begin{theorem}\label{eqcase:five-six}
    Let $n=5$ or $6$, let $E$ have nonnegative sectional curvature, and let $\omega$ be a unit three-form. Then
    \[
        Q_3=-\frac13S
    \]
    holds precisely in the following cases: either $E=0$, or there are an orthogonal splitting $V=W\oplus F$, with $\dim W=4$, an orthogonal complex structure $J$ on $W$, a number $t>0$, and a unit $\zeta\in F^*$ such that
    \begin{equation}\label{eqcase:five-six-model}
        E=tE_J,\qquad \omega=\frac1{\sqrt2}\Omega_J\wedge\zeta.
    \end{equation}
    Here $J$ is extended by zero on $F$. In the nonzero case,
    \[
        S=12t, \qquad q_3(E)\omega=-4t\omega, \qquad \ker\bigl(q_3(E)+S/3\bigr)=\Omega_J\wedge F^*.
    \]
    Thus the last eigenspace has dimension one in dimension five and two in dimension six.
\end{theorem}

\begin{proof}
    We first prove necessity in dimension six. The normal form used in Section \ref{six:scalar-estimate} is
    \[
        \omega=(\cos\theta+\sin\theta*)\omega_0, \qquad \omega_0=a e^{123}+b e^{145}+c e^{246}+d e^{356}, \qquad a^2+b^2+c^2+d^2=1.
    \]
    By Lemma \ref{six:normal-form}, the Hodge phase preserves the $q_3$ quadratic form. Thus it is enough to classify equality for $\omega_0$ and then restore the phase.

    The six NW applications in Section \ref{six:scalar-estimate}, followed by $2\sqrt{xy}\le x+y$, give
    \begin{equation}\label{eqcase:six-budget}
        Q_3\ge\sum_{i<j}C_{ij}K_{ij}\ge-\frac23\sum_{i<j}K_{ij}.
    \end{equation}
    For clarity, the equality information in this coefficient estimate is recorded explicitly. By \eqref{six:opposite-pair-final},
    \[
        C_{23}=C_{45}=c^2+d^2-\frac23|cd|-\frac43|ab|\ge\frac23(c^2+d^2)-\frac23(a^2+b^2)=-\frac23+\frac43(c^2+d^2)\ge-\frac23.
    \]
    Hence
    \begin{equation}\label{eqcase:six-coefficient-equality}
        C_{23}=-\frac23 \quad \Longleftrightarrow \quad c=d=0, \quad |a|=|b|=1/\sqrt2.
    \end{equation}
    The coefficients of the other ten nonexceptional sectional curvatures have precisely this form after permuting the four coefficient labels. The coefficients of the three exceptional sectional curvatures $K_{16},K_{25},K_{34}$ are at least $1/3$; see Section \ref{six:scalar-estimate}.

    Assume now $E\ne0$. Lemma \ref{eqcase:nullity} gives $S>0$, so at least one coordinate sectional curvature is positive. Equality in \eqref{eqcase:six-budget} forces its coefficient to be $-2/3$; it cannot be exceptional. Permuting the four monomials and changing frame signs, we may take this sectional curvature to be $K_{23}$. Such permutations are frame changes: the six coordinates correspond to the six pairwise intersections of the four triples $123,145,246,356$, and every permutation of those four triples permutes these intersections. Hence taking that sectional curvature to be $K_{23}$ is indeed without loss of generality. Now by \eqref{eqcase:six-coefficient-equality}, $c=d=0$ and $|a|=|b|=1/\sqrt2$. After substitution, the coefficients of the other coordinate sectional curvatures are strictly greater than $-2/3$. Therefore all other coordinate sectional curvatures must vanish and the only sectional coefficients equal to $-2/3$ are $C_{23},C_{45}$.

    In Section \ref{six:scalar-estimate}, to obtain the left-hand side of \eqref{six:opposite-pair-final}, we applied $2\sqrt{xy}\le x+y$ to \eqref{six:opposite-pair}. Equality in this step gives $K_{23}=K_{45}=k>0$. We also made the following estimate:
    \[
        -4abE_{2345} \ge -4|ab||E_{2345}| \ge -\frac{2}{3}|ab|(K_{24}+K_{25}+K_{34}+K_{35})-\frac{8}{3}|ab|\sqrt{K_{23}K_{45}}.
    \]
    The first equality gives $\operatorname{sgn}(ab)=\operatorname{sgn}(E_{2345})$. The second equality gives $|E_{2345}|=\frac23\sqrt{K_{23}K_{45}}=\frac{2k}{3}$ since $K_{24}+K_{25}+K_{34}+K_{35}=0$. Thus $E_{2345}=(2k/3)\operatorname{sgn}(ab)$. Since the only nonvanishing coordinate sectional curvatures are $K_{23}$ and $K_{45}$, directions $e_1$ and $e_6$ have zero Ricci curvature and hence are full curvature-null directions by Lemma \ref{eqcase:nullity}. Thus $E$ is supported on $W=\operatorname{span}\{e_2,e_3,e_4,e_5\}$. Lemma \ref{eqcase:four-dimensional} gives
    \[
        E=(k/3)E_J, \qquad \Omega_J=e^{23}+\operatorname{sgn}(ab)e^{45}.
    \]
    After absorbing an overall sign in the flat covector, $\omega_0=\Omega_J\wedge e^1/\sqrt2$. Then
    \begin{align*}
        (\cos\theta+\sin\theta\,*)\omega_0 &= \frac1{\sqrt2}\left(\cos\theta\,\Omega_J\wedge e^1+\sin\theta\,*(\Omega_J\wedge e^1)\right)\\
        &= \frac1{\sqrt2}\left(\cos\theta\,\Omega_J\wedge e^1\pm\sin\theta\,\Omega_J\wedge e^6\right)\\
        &= \frac1{\sqrt2}\Omega_J\wedge\left(\cos\theta\,e^1\pm\sin\theta\,e^6\right) = \frac1{\sqrt2}\Omega_J\wedge\zeta,
    \end{align*}
    where $\zeta=\cos\theta\,e^1\pm\sin\theta\,e^6\in\operatorname{span}\{e^1,e^6\}=F^*$.

    For $n=5$, consider the six-dimensional space $V\oplus\mathbb Rz$, where $z\perp V$ and $|z|=1$. Extend $E$ and $\omega$ by zero to $V\oplus\mathbb Rz$, so $z$ is curvature-null and $\iota_z\omega=0$. The extension preserves $S$ and $Q_3$ and preserves sectional nonnegativity by orthogonal projection of both arguments of the sectional numerator. The six-dimensional result gives $W,F,J,t,\zeta$ for $V\oplus\mathbb Rz$. It remains to show that $W\subset V$ and $\zeta\in(F\cap V)^*$. Since $z$ is curvature-null, $\mathbb Rz\subset\ker R$, i.e., $(\ker R)^\perp\subset V$, where $R=\operatorname{Ric}(E)^\sharp$. But $R=3tP_W$, so $W=(\ker R)^\perp\subset V$. Since $z\perp W$, also $\iota_z\omega=\zeta(z)\Omega_J/\sqrt2$. Thus $\iota_z\omega=0$ gives $\zeta(z)=0$, so $\zeta\in V^*$. This together with $\zeta\in F^*$ implies $\zeta\in(F\cap V)^*$.

    For sufficiency and the entire eigenspace, on $W^4$ the curvature operator is
    \[
        \mathcal R_{E_J}(X^\flat\wedge Y^\flat)=(JX)^\flat\wedge(JY)^\flat+2\Omega_J(X,Y)\Omega_J,
    \]
    where $X,Y\in W$. It has eigenvalues $5,1,-1$; the eigenspace for $5$ is $\mathbb R\Omega_J$, that for $1$ is $\operatorname{span}\{e^{12}-e^{34},e^{13}+e^{24},e^{14}-e^{23}\}$, and that for $-1$ is $\operatorname{span}\{e^{13}-e^{24},e^{14}+e^{23}\}$. Since $\operatorname{Ric}(E_J)=3g_W$, $q_2(E_J)=6\operatorname{Id}-2\mathcal R_{E_J}$ has eigenvalues $-4,4,8$ on those spaces. The splitting $V=W\oplus F$ gives
    \[
        \Lambda^3V^*=\bigoplus_{j=0}^3\operatorname{span}\{\alpha\wedge\beta \mid \alpha\in\Lambda^jW^*,\ \beta\in\Lambda^{3-j}F^*\}, \qquad q_3(tE_J)(\alpha\wedge\beta)=tq_j(E_J|_W)\alpha\wedge\beta.
    \]
    Since $q_0(E_J|_W)\equiv0$, all its eigenvalues are $0$. For $j=1$ and $j=3$, $q_j(E_J|_W)=3\operatorname{Id}$, so all its eigenvalues are $3$. Hence the eigenvalue $-4$ occurs only when $j=2$. The $-4t$ eigenspace is $\operatorname{span}\{\alpha\wedge\beta \mid \alpha\in\mathbb R\Omega_J,\ \beta\in F^*\}$. Finally $S=12t$, as asserted.
\end{proof}

\subsection{The two-form Ricci endpoint}

The next lemma is stronger than the scalar equality statement in degree two: it allows unequal singular values of the model endomorphism. This distinction is needed when equality is descended through hyperplane contractions.

\begin{lemma}\label{eqcase:two-form-endpoint}
    Let $\dim H=2m$ or $2m+1$, with $m\ge2$, and let $E$ have nonnegative sectional curvature on $H$. Put $c=2(m-1)/3$. Suppose a nonzero $\alpha\in\Lambda^2H^*$ satisfies
    \[
        \langle(q_2(E)+c(\operatorname{Ric}(E)^\sharp)^{[2]})\alpha,\alpha\rangle=0, \qquad F_2(E,\alpha)>0.
    \]
    Then $E=E_A$ for a skew-adjoint $A$ of rank $2m$. In an orthonormal frame $\{e_i\}$ in which $\alpha$ is in orthogonal normal form $\alpha=\sum_{a=1}^m\lambda_a e^{p_aq_a}$, every $\lambda_a$ is nonzero and
    \[
        \Omega_A=\sum_{a=1}^m t_a e^{p_aq_a}, \qquad t_a=\operatorname{sgn}(\lambda_a)\sqrt{K_{p_aq_a}/3}, \qquad \lambda_at_a\text{ is independent of }a.
    \]
    In odd dimension the unused direction is in the full curvature nullity. In particular, every Jacobi operator $J_x^\sharp$ of $E$ has rank at most one, where $J_x(y,z):=E(x,y,x,z)$.
\end{lemma}

\begin{proof}
    Use the notation $k_a=K_{p_aq_a}$, $A_{ab}=K_{p_ap_b}+K_{p_aq_b}+K_{q_ap_b}+K_{q_aq_b}$, and $r_{ab}=E_{p_aq_ap_bq_b}$. If there is an unused unit direction $z$, put $Z_a=K_{p_az}+K_{q_az}$; otherwise put $Z_a=0$. Write $x_a=|\lambda_a|\sqrt{k_a}$. At the lower endpoint $a=-2(m-1)/3$ of Proposition \ref{alg:two-affine}, the scalar coefficient $c_m(a)=0$. Then \eqref{alg:two-affine-equation} becomes $q_2(E)+c(\operatorname{Ric}(E)^\sharp)^{[2]}\ge0$. From \eqref{alg:estimation-of-q2},
    \begin{align}
        Q_2+cF_2 &\ge 2c\sum_ax_a^2 - \frac83\sum_{a<b}x_ax_b + \sum_{a<b}\left[(1+c)(\lambda_a^2+\lambda_b^2)-\frac23|\lambda_a\lambda_b|\right]A_{ab} + (1+c)\sum_a\lambda_a^2Z_a \notag\\
        &= \frac43\sum_{a<b}(x_a-x_b)^2 + \sum_{a<b}\left[(1+c)(\lambda_a^2+\lambda_b^2)-\frac23|\lambda_a\lambda_b|\right]A_{ab} + (1+c)\sum_a\lambda_a^2Z_a, \label{eqcase:ricci-remainder}
    \end{align}
    where we used $(m-1)\sum_a x_a^2-2\sum_{a<b}x_ax_b=\sum_{a<b}(x_a-x_b)^2$. Furthermore, the bracket is $\ge(c+2/3)(\lambda_a^2+\lambda_b^2)=(2m/3)(\lambda_a^2+\lambda_b^2)$. Thus every term on the second line of \eqref{eqcase:ricci-remainder} is nonnegative, so all vanish. In particular, the $x_a$ have a common value $x_0$. By \eqref{alg:two-Ricci-expansion}, $0<F_2=2mx_0^2$, so $x_a>0$ for every $a$. Hence all $\lambda_a,k_a$ are nonzero, and $A_{ab}=Z_a=0$, so all coordinate sectional curvatures are zero except $k_a$. Jacobi--Cauchy--Schwarz inequality \eqref{alg:Jacobi-CS} then eliminates all components with exactly three different indices. For a component $E_{ijkl}$ with four distinct indices, if its four indices are not the union of two complete blocks, then there exists at least one index $i_0\in\{i,j,k,l\}$ such that the index paired with $i_0$ is not among these four indices; hence $\operatorname{Ric}(E|_U)(e_{i_0},e_{i_0})=0$, where $U:=\operatorname{span}\{e_i,e_j,e_k,e_l\}$; Lemma \ref{eqcase:nullity} therefore gives $E_{ijkl}=0$. In odd dimension, $z$ is Ricci-null on all of $H$, so every component involving $z$ vanishes by Lemma \ref{eqcase:nullity}.

    In the proof of Proposition \ref{alg:two-affine}, we used NW inequality to obtain
    \[
        -4\lambda_a\lambda_b r_{ab}\ge -4|\lambda_a\lambda_b| |r_{ab}|\ge -\frac23|\lambda_a\lambda_b|A_{ab}-\frac83|\lambda_a\lambda_b|\sqrt{k_a k_b}.
    \]
    Since $A_{ab}=0$, equality in both inequalities yields $r_{ab}=(2/3)\operatorname{sgn}(\lambda_a\lambda_b)\sqrt{k_ak_b}$. Hence on $U_{ab}:=\operatorname{span}\{e_{p_a},e_{q_a},e_{p_b},e_{q_b}\}$, Lemma \ref{eqcase:four-dimensional} gives $E|_{U_{ab}}=E_{B_{ab}}|_{U_{ab}}$, where $B_{ab}$ is determined by
    \[
        \Omega_{B_{ab}}=\sqrt{k_a/3}\,e^{p_aq_a}+\operatorname{sgn}(\lambda_a\lambda_b)\sqrt{k_b/3}\,e^{p_bq_b}.
    \]
    Thus $\Omega_A|_{U_{ab}}=\operatorname{sgn}(\lambda_a)\Omega_{B_{ab}}$, which gives $E_A|_{U_{ab}}=E_{B_{ab}}|_{U_{ab}}$. Therefore $E|_{U_{ab}}=E_A|_{U_{ab}}$. Since all other components of $E$ and $E_A$ are zero, $E=E_A$. Moreover $\lambda_at_a=x_a/\sqrt3=x_0/\sqrt3$ is independent of $a$. Finally, the Jacobi operator $J_x^\sharp(y)=E_A(x,y,x,\cdot)^\sharp=3\langle Ax,y\rangle Ax$ has rank one (when $Ax\ne0$) or is the zero operator (when $Ax=0$).
\end{proof}

\subsection{Rank-one Jacobi operators determine a single model tensor}

\begin{lemma}\label{eqcase:rank-one-jacobi}
    Suppose $E$ has nonnegative sectional curvature and every Jacobi operator has rank at most one. Then $E=E_A$ for some skew-adjoint endomorphism $A$.
\end{lemma}

\begin{proof}
    In a fixed basis let $J(x)$ be the matrix of $J_x$, i.e., $J_{ij}(x)=E(x,e_i,x,e_j)$. Its entries are homogeneous quadratic polynomials in $x=(x_1,\cdots,x_n)$, and it is positive semidefinite for every $x$. If $J(x)=0$ for all $x$, then ${\rm sec}(E)=0$ and hence $E=0=E_0$. Otherwise some diagonal entry is a nonzero nonnegative quadratic polynomial in $x$; without loss of generality, assume it is $J_{11}=:f$. Put $g_i:=J_{i1}$. Since $J(x)$ has rank at most one for every $x$, all $2\times2$ minors vanish. In particular,
    \begin{equation}\label{eqcase:polynomial-minors}
        fJ_{ij}=g_ig_j.
    \end{equation}
    Work in the unique factorization domain $\mathbb R[x_1,\cdots,x_n]$. If $f$ is irreducible, it is prime. From $f\mid g_i^2$ we get $f\mid g_i$. Since both $f$ and $g_i$ are homogeneous of degree two, $g_i=a_if$ for constants $a_i$. Here $a_1=1$. Cancel $f$ in \eqref{eqcase:polynomial-minors} to obtain $J=faa^T$, with the nonzero constant vector $a=(a_i)$. But $J(x)x=0$ would imply $f(x)a\langle a,x\rangle=0$, i.e., for every $i$, $f(x)a_i\langle a,x\rangle=0$. Take $i=1$ to obtain $f(x)\langle a,x\rangle=0$, impossible in an integral domain. Hence $f$ is reducible over $\mathbb R$.

    A nonzero homogeneous nonnegative quadratic that is reducible must have the form $f=b\ell^2$, with $b>0$ and $\ell$ a real linear form: if $f=\ell_1\ell_2$, where $\ell_1$ and $\ell_2$ are two nonproportional linear factors, then $\ell_1$ and $\ell_2$ take opposite signs on some open set, contradicting $f\ge0$. Since $\ell^2\mid g_i^2$, we have $g_i=\ell h_i$ with $h_i$ linear. Equation \eqref{eqcase:polynomial-minors} now gives $J_{ij}=h_ih_j/b$. Thus
    \[
        J(x)=(Bx)(Bx)^T,
    \]
    where $B$ is a linear map with $(Bx)_i=h_i(x)/\sqrt b$. The identity $J(x)x=0$ says $(Bx)\langle Bx,x\rangle=0$ as a polynomial vector identity. Since $J\not\equiv0$, $B\neq0$, so at least one coordinate of $Bx$ is a nonzero polynomial. Thus $\langle Bx,x\rangle\equiv0$. Polarization gives $\langle Bx,y\rangle+\langle By,x\rangle=0$, i.e., $B+B^T=0$.

    It follows that $E(x,y,x,y)=\langle J(x)y,y\rangle=\langle Bx,y\rangle^2$ for all $x,y$. By the definition of $E_B$ in \eqref{eq:EJ}, $E_B(x,y,x,y)=3\langle Bx,y\rangle^2$. Since the sectional curvature determines the algebraic curvature tensor, $E=E_B/3=E_{B/\sqrt3}$, as required.
\end{proof}

\subsection{Complete equality in even dimension and even degree}

\begin{theorem}\label{eqcase:even}
    Let $n=2m$, $p=2a$, and $1\le a\le m-1$. Suppose $E$ has nonnegative sectional curvature and $|\omega|=1$. Equality in the sharp estimate of Theorem \ref{alg:even-sharp}, namely
    \[
        Q_p=-\frac{p(n-p)}{3n}S,
    \]
    holds if and only if either $E=0$, or there are an orthogonal complex structure $J$ on all of $V$ and $t>0$ such that
    \begin{equation}\label{eqcase:even-model}
        E=tE_J, \qquad \omega=\pm\frac{\Omega_J^a}{a!\sqrt{\binom ma}}.
    \end{equation}
    In particular, the lowest eigenspace of a nonzero equality tensor is one-dimensional.
\end{theorem}

\begin{proof}
    Suppose first that $E\ne0$, so $S>0$. Set $c=n-p$, which is a positive even integer. The proof of Theorem \ref{alg:even-sharp} combines the two nonnegative remainders $Q_p+\frac c3F_p$ and $Q_p+\frac p3(S-F_p)$ with positive weights $p/n$ and $c/n$. Equality in the scalar estimate forces both to be zero. Thus
    \begin{equation}\label{eqcase:even-ricci-value}
        Q_p=-\frac c3F_p, \qquad F_p=\frac pnS>0.
    \end{equation}

    We now give a descending induction. At a degree $j\ge3$ and dimension $N=j+c$, write $E_i=E|_{e_i^\perp}$ and $\alpha_i=\iota_{e_i}\omega$. By Proposition \ref{alg:averages},
    \begin{align}
        \overline{\sum_i\left(Q_{j-1}(E_i,\alpha_i)+\frac c3F_{j-1}(E_i,\alpha_i)\right)} &= \frac{(j-1)(N+1)}{N+2}\left(Q_j+\frac c3F_j\right),\label{eqcase:descent}\\
        \overline{\sum_iF_{j-1}(E_i,\alpha_i)} &= \frac{3Q_j+[(j-1)N+j-4]F_j}{N+2}.\label{eqcase:descent-ricci}
    \end{align}
    By Theorem \ref{alg:Ricci-all}, $Q_{j-1}(E_i,\alpha_i)+(c/3)F_{j-1}(E_i,\alpha_i)\ge0$. If $Q_j+(c/3)F_j=0$, then the right of \eqref{eqcase:descent} is zero, so $Q_{j-1}(E_i,\alpha_i)+(c/3)F_{j-1}(E_i,\alpha_i)=0$ for every frame and every $i$, which does not require $F_j>0$. If in addition $F_j>0$, i.e.,
    \[
        Q_j+\frac c3F_j=0, \qquad F_j>0,
    \]
    then by \eqref{eqcase:descent-ricci},
    \[
        \overline{\sum_iF_{j-1}(E_i,\alpha_i)}=(j-2)F_j>0.
    \]
    In particular, there is an oriented orthonormal frame $\{e_i\}$ and an index $i$ such that $F_{j-1}(E_i,\alpha_i)>0$. Choosing such a frame and such an $i$, we reduce to the case $j-1$:
    \[
        Q_{j-1}(E_i,\alpha_i)+\frac c3F_{j-1}(E_i,\alpha_i)=0, \qquad F_{j-1}(E_i,\alpha_i)>0.
    \]
    Since the base case $j=p$ has been established by \eqref{eqcase:even-ricci-value}, repeating the argument reduces to the case $j=2$. This shows that
    \[
        Q_2(E|_H,\alpha_H)+\frac c3F_2(E|_H,\alpha_H)=0
    \]
    for every ordered orthonormal $(p-2)$-frame $\{e_1,\cdots,e_{p-2}\}$, where
    \[
        H=\bigl(\operatorname{span}\{e_1,\cdots,e_{p-2}\}\bigr)^\perp, \qquad \alpha_H=\iota_{p-2}\cdots\iota_1\omega\in\Lambda^2H^*,
    \]
    and there exists at least one such frame for which $F_2(E|_H,\alpha_H)>0$.

    Fix an orthonormal frame $\{e_i\}$ of $V$. For $i=1,\cdots,p-2$, define maps $\varphi_i:\operatorname{SO}(n)\to V$ by $\varphi_i(g)=ge_i$. Set
    \[
        H_g=(\operatorname{span}\{\varphi_1(g),\cdots,\varphi_{p-2}(g)\})^\perp, \qquad \alpha_g=\iota_{\varphi_{p-2}(g)}\cdots\iota_{\varphi_1(g)}\omega\in\Lambda^2H_g^*.
    \]
    Since $g\in\operatorname{SO}(n)$, we have $H_g=\operatorname{span}\{g e_{p-1},\cdots,g e_n\}$. Define the function $\Phi(g):=F_2(E|_{H_g},\alpha_g)$, which is a polynomial in $g_{ij}$, $E_{ijkl}$, and $\omega_{i_1,\cdots,i_p}$, and hence is real analytic in $g_{ij}$. By the preceding paragraph, $\Phi\not\equiv0$. If $P:=\{g\in\operatorname{SO}(n) \mid \Phi(g)>0\}$ is not dense, then there exists a nonempty open set $O\subset \operatorname{SO}(n)\setminus P=\{g\mid\Phi(g)=0\}$. Thus the real analytic function $\Phi$ vanishes identically on the nonempty open set $O$; by real analytic continuation and the connectedness of $\operatorname{SO}(n)$, $\Phi\equiv0$, a contradiction. Therefore, $\overline P=\operatorname{SO}(n)$.

    Take $H=H_g$ in Lemma \ref{eqcase:two-form-endpoint}. Since
    \[
        \dim H_g=c+2=2\left(\frac c2+1\right) \qquad {\rm and} \qquad \frac c3=\frac{2((\frac c2+1)-1)}3,
    \]
    Lemma \ref{eqcase:two-form-endpoint} applies on $P$ and shows that for all $g\in P$, all Jacobi operators of $E|_{H_g}$ have rank at most one. We first claim that for all $g\in\operatorname{SO}(n)$, all Jacobi operators of $E|_{H_g}$ have rank at most one. Indeed, for $t=(t_{p-1},\cdots,t_n)\in\mathbb R^{c+2}$, let
    \[
        x(g,t)=\sum_{i=p-1}^n t_ige_i\in H_g,
    \]
    and define the matrix
    \[
        A_{ij}(g,t)=E(x(g,t),ge_i,x(g,t),ge_j), \qquad i,j=p-1,\cdots,n.
    \]
    This is the matrix of $J_{x(g,t)}$ on $H_g$. Its rank being at most one is equivalent to the vanishing of all $2\times2$ minors, that is, $A_{ij}A_{kl}-A_{il}A_{jk}=0$. For fixed $t,i,j,k,l$, the left is a continuous function of $g$. It vanishes identically on the dense set $P$, and hence identically on all of $\mathrm{SO}(n)$. This proves the claim. We further claim that all Jacobi operators of $E$ have rank at most one. Otherwise, suppose that there exists $x\in V$ such that the rank of $J_x$ is at least $2$. Since $J_x$ is symmetric, we may choose two orthogonal eigenvectors $y,z$ corresponding to nonzero eigenvalues. Then
    \[
        \begin{pmatrix}
            E(x,y,x,y)&E(x,y,x,z)\\
            E(x,z,x,y)&E(x,z,x,z)
        \end{pmatrix}
        =
        \begin{pmatrix}
            \lambda&0\\
            0&\mu
        \end{pmatrix},
        \qquad \lambda\mu\ne0.
    \]
    Therefore, on any $(c+2)$-dimensional subspace $H$ containing $x,y,z$, the rank of $J_x|_H$ is at least $2$. Since $c+2\ge4\ge3$, a $(c+2)$-dimensional subspace $H$ containing $x,y,z$ indeed exists. But every $(c+2)$-dimensional subspace $H$ is of the form $H_g$ for some $g\in\operatorname{SO}(n)$, which contradicts the preceding claim. This proves the further claim. Lemma \ref{eqcase:rank-one-jacobi} now gives $E=E_A$.

    It remains to analyze equality for this single $E_A$. Choose orthonormal two-planes $P_i$ and orthonormal bases $\{e_{2i-1},e_{2i}\}$ of $P_i$ such that
    \[
        A|_{P_i}=t_iJ_i, \qquad t_i\ge0, \qquad J_ie_{2i-1}=e_{2i}, \quad J_ie_{2i}=-e_{2i-1}.
    \]
    Complexify the exterior algebra and use the orthonormal basis $\{1,u_i^+,u_i^-,\theta_i\}$ of $(\Lambda^*P_i^*)_{\mathbb C}$ as in the proof of Theorem \ref{strong:sharp}. Fix a set $L\subset\{1,\cdots,m\}$ of $s$ two-planes from which we choose one-form basis vectors. For every $j\in L$, choose a sign $\sigma_j\in\{+1,-1\}$, and use the one-form $u_j^{\sigma_j}$. Among the remaining $d=m-s$ two-planes, choose exactly $r$ planes to contribute their full two-form $\theta_i$, where $r$ is determined by the degree condition
    \begin{equation}\label{eqcase:degree-count}
        2a=2r+s.
    \end{equation}
    Let $I\subset\{1,\cdots,m\}\setminus L$ with $|I|=r$. Define $f_I=(\bigwedge_{j\in L}u_j^{\sigma_j})\wedge(\bigwedge_{i\in I}\theta_i)$. The vectors $f_I$ form an orthonormal basis of the subspace
    \[
        \mathcal H_{L,\sigma,r}:=\operatorname{span}_{\mathbb C}\left\{f_I \mid I\subset\{1,\cdots,m\}\setminus L,\ |I|=r\right\} \subset (\Lambda^{2a}V^*)_{\mathbb C}.
    \]
    As in the proof of Theorem \ref{strong:sharp},
    \[
        q_{2a}(E_A)\big|_{\mathcal H_{L,\sigma,r}}=\Bigl[2\sum_{j\in L}t_j^2+\Bigl(\sum_{j\in L}\sigma_jt_j\Bigr)^2\Bigr]\operatorname{Id}-4\mathcal A_r, \qquad \mathcal A_rf_I:=\sum_{i\in I}\sum_{k\notin I\cup L}t_it_kf_{(I\setminus\{i\})\cup\{k\}}.
    \]
    The two Cauchy--Schwarz estimates for $\mathcal A_r$ in the proof of Theorem \ref{strong:sharp} give
    \[
        \langle \mathcal A_r f,f\rangle_{\mathbb C} \le \frac{r(d-r)}d\sum_{i\in\{1,\cdots,m\}\setminus L}t_i^2 \le \frac{r(d-r)}d\sum_{i=1}^m t_i^2.
    \]
    By \eqref{eqcase:degree-count} and $d=m-s$, we have $r=a-s/2\le a$ and $d-r=m-a-s/2\le m-a$. Therefore,
    \[
        \frac{r(d-r)}d=\frac{r(d-r)}{r+(d-r)}\le\frac{a(m-a)}m.
    \]
    Consequently,
    \begin{equation}\label{eqcase:estimate-of-q2a}
        q_{2a}(E_A)\big|_{\mathcal H_{L,\sigma,r}} \ge -\frac{4a(m-a)}m\sum_{i=1}^m t_i^2 = -\frac{2a(m-a)}{3m}\operatorname{Scal}(E_A).
    \end{equation}

    Now let $f\in\Lambda^{2a}V^*\subset(\Lambda^{2a}V^*)_{\mathbb C}$ be a real unit form attaining the equality, then
    \begin{align*}
        &\quad \sum_{\mathcal H}\langle q_{2a}(E_A)f_{\mathcal H},f_{\mathcal H}\rangle_{\mathbb C} = \langle q_{2a}(E_A)f,f\rangle_{\mathbb C} \\
        &= -\frac{2a(m-a)}{3m}\operatorname{Scal}(E_A) = -\frac{2a(m-a)}{3m}\operatorname{Scal}(E_A)\sum_{\mathcal H}|f_{\mathcal H}|^2.
    \end{align*}
    This together with \eqref{eqcase:estimate-of-q2a} shows that for every subspace $\mathcal H$,
    \[
        \langle q_{2a}(E_A)f_{\mathcal H},f_{\mathcal H}\rangle_{\mathbb C}=-\frac{2a(m-a)}{3m}\operatorname{Scal}(E_A)|f_{\mathcal H}|^2.
    \]
    If $s>0$, then both $r=a-s/2<a$ and $d-r=m-a-s/2<m-a$, so $r(d-r)/d<a(m-a)/m$; in this case the strict inequality in \eqref{eqcase:estimate-of-q2a} holds; therefore $f$ cannot have a component in $\mathcal H$ when $s>0$, since $\sum_it_i^2>0$, i.e., $A\ne0$, which is guaranteed by $E\ne0$. If $s=0$, then $r=a$ and $d=m$; in this case the subspace becomes
    \[
        \mathcal H_{\varnothing,\varnothing,a}=\operatorname{span}_{\mathbb C}\left\{\theta_I=\bigwedge_{i\in I}\theta_i \;\middle|\; I\subset\{1,\cdots,m\},\ |I|=a\right\}.
    \]
    Thus we can write $f=\sum_{|I|=a}a_I\theta_I$, where $a_I\in\mathbb R$ with $\sum a_I^2=1$. The two Cauchy--Schwarz estimates used to bound $\mathcal A_a$ were combined with strictly positive weights $a/m$ and $(m-a)/m$. Both must be equalities. Hence for any $K\subset\{1,\cdots,m\}$ with $|K|=a-1$ and any $i,j\notin K$,
    \begin{equation}\label{eqcase:two-cs-equalities}
        t_i a_{K\cup\{i\}}=t_j a_{K\cup\{j\}}, \qquad t_j a_{K\cup\{i\}}=t_i a_{K\cup\{j\}}.
    \end{equation}

    We claim that $t_i>0$ for all $i=1,\cdots,m$. Suppose for contradiction that some $t_{i_0}=0$. Since $A\ne0$, some $t_k>0$; since $f\ne0$, some $a_I\ne0$. Because $a<m$, the set $\{1,\cdots,m\}\setminus I\ne\varnothing$. If $t_i=t_j$ for all $i\in I$ and $j\notin I$, then $t_1=\cdots=t_m$, contradicting $t_k>0=t_{i_0}$. Hence there must exist $i\in I$ and $j\notin I$ such that $t_i\ne t_j$. Let $K=I\setminus\{i\}$, then $i,j\notin K$. By \eqref{eqcase:two-cs-equalities},
    \[
        t_ia_I=t_ja_{K\cup\{j\}}, \qquad t_ja_I=t_ia_{K\cup\{j\}}.
    \]
    Then $t_j^2a_I=t_it_ja_{K\cup\{j\}}=t_i^2a_I$. Since $a_I\ne0$, $t_j^2=t_i^2$, contradicting $t_i\ne t_j$ and $t_i,t_j\ge0$. This proves the claim.

    We next claim that all $a_I\ne0$. Since all $t_i>0$, we have $a_{K\cup\{i\}}=(t_j/t_i)a_{K\cup\{j\}}$ by \eqref{eqcase:two-cs-equalities}; hence for any $I,J\subset\{1,\cdots,m\}$ with $|I|=|J|=a$ and $|I\cap J|=a-1$, we have
    \[
        a_I\ne0 \quad \Longleftrightarrow \quad a_J\ne0.
    \]
    Since $f\ne0$, some $a_I\ne0$. Then for any $J\subset\{1,\cdots,m\}$ with $|J|=a$, it follows by finitely many successive applications of the above equivalence that $a_J\ne0$.

    Now we show that all $t_i$ are equal. Fix $i,j=1,\cdots,m$ with $i\ne j$. Since $a-1\le m-2$, we can choose $K\subset\{1,\cdots,m\}\setminus\{i,j\}$ with $|K|=a-1$. By \eqref{eqcase:two-cs-equalities}, $t_j^2a_{K\cup\{i\}}=t_it_ja_{K\cup\{j\}}=t_i^2a_{K\cup\{i\}}$. Since all $a_I\ne0$, we have $t_i=t_j$.

    It follows that all $a_I$ are also equal. Indeed, by \eqref{eqcase:two-cs-equalities}, we have $a_{K\cup\{i\}}=(t_j/t_i)a_{K\cup\{j\}}=a_{K\cup\{j\}}$; hence for any $I,J\subset\{1,\cdots,m\}$ with $|I|=|J|=a$ and $|I\cap J|=a-1$, we have $a_I=a_J$. Therefore, for any two different $a$-element subsets $I'$ and $J'$, by finitely many successive applications of this identity, $a_{I'}=a_{J'}$.

    Since all $t_i$ are equal, $A=t'J$ for an orthogonal complex structure $J$ on $V$. Let $t=(t')^2$, then $E=E_A=E_{t'J}=(t')^2E_J=tE_J$. Also, $\Omega_J=\sum_{i=1}^m\theta_i$ and $\Omega_J^a=a!\sum_{|I|=a}\theta_I$. Since all $a_I$ are equal and $\sum_{|I|=a}a_I^2=1$, we have
    \[
        f=\sum_{|I|=a}a_I\theta_I=\pm\frac1{\sqrt{\binom ma}}\sum_{|I|=a}\theta_I=\pm\frac{\Omega_J^a}{a!\sqrt{\binom ma}}.
    \]
    Thus the equality line is $\mathbb R\Omega_J^a$. Conversely the model calculation in Theorem \ref{alg:even-sharp} gives $q_{2a}(tE_J)\Omega_J^a=-4a(m-a)t\Omega_J^a$ and $S=6mt$, proving equality.
\end{proof}

\begin{proof}[Proof of Theorem~\ref{intro:main}]
    Theorem~\ref{alg:scalar-all} combines the Ricci-weighted bounds of Theorems~\ref{alg:Ricci-all} and~\ref{alg:three-lift} in complementary degrees to give the coefficient $C_{n,p}$ for $2\le p\le n-2$. It also gives the coefficient zero in degrees $0,1,n-1,n$. In dimension seven, Theorem~\ref{spin:rational} improves the degree-three coefficient to the value in \eqref{intro:Chat}. By Lemma~\ref{lem:HodgeDual}, the Hodge star intertwines $q_3(E)$ and $q_4(E)$ in this dimension, so the same coefficient applies in degree four.

    For even $n$ and even $2\le p\le n-2$, Theorem~\ref{alg:scalar-all} gives the sharp coefficient $p(n-p)/(3n)$, and Theorem~\ref{eqcase:even} identifies the equality pairs in \eqref{eqcase:even-model}. The sharp degree-three coefficients in dimensions four, five and six follow from Corollary~\ref{six:sharp}.
\end{proof}

\subsection{Degree-three equality under strong nonnegativity}

\begin{theorem}\label{eqcase:strong-three}
    Suppose $n\ge5$ and $E$ is strongly nonnegative. Put $M=\lfloor(n-1)/2\rfloor$ and $d_n=2(M-1)/(3M)$. For a unit three-form, $Q_3=-d_nS$ holds if and only if $E=0$, or
    \[
        V=W\oplus F,\quad \dim W=2M,\quad E=tE_J, \quad \omega=\frac1{\sqrt M}\Omega_J\wedge\zeta,
    \]
    where $t>0$, $J$ is an orthogonal complex structure on $W$, extended by zero, and $\zeta$ is a unit covector in $F^*$. For a nonzero equality tensor, the full equality space is $\Omega_J\wedge F^*$, of dimension one for odd $n$ and two for even $n$.
\end{theorem}

\begin{proof}
    We first classify the equality case for the model tensor $E_A\ne0$. Let $m=\lfloor n/2\rfloor$. Choose orthonormal two-planes $P_i$ and orthonormal bases $\{e_{2i-1},e_{2i}\}$ of $P_i$ such that
    \[
        A|_{P_i}=t_iJ_i, \qquad t_i\ge0, \qquad J_ie_{2i-1}=e_{2i}, \quad J_ie_{2i}=-e_{2i-1}.
    \]
    If $n$ is odd, denote the remaining unit covector by $z$. Complexify the exterior algebra and use the orthonormal basis $\{1,u_i^+,u_i^-,\theta_i\}$ of $(\Lambda^*P_i^*)_{\mathbb C}$ as in the proof of Theorem \ref{strong:sharp}. Fix a set $L\subset\{1,\cdots,m\}$ of $s$ two-planes from which we choose one-form basis vectors. For every $j\in L$, choose a sign $\sigma_j\in\{+1,-1\}$, and use the one-form $u_j^{\sigma_j}$. If $n$ is odd, choose whether $z$ occurs. Denote this choice by $\epsilon\in\{0,1\}$. Among the remaining $d=m-s$ two-planes, choose exactly $r$ planes to contribute their full two-form $\theta_i$, where $r$ is determined by the degree condition
    \begin{equation}\label{eqcase:degree-count2}
        3=2r+s+\epsilon.
    \end{equation}
    Let $I\subset\{1,\cdots,m\}\setminus L$ with $|I|=r$. Define
    \[
        f_I=\Bigl(\bigwedge_{j\in L}u_j^{\sigma_j}\Bigr)\wedge\Bigl(\bigwedge_{i\in I}\theta_i\Bigr)\wedge z^\epsilon.
    \]
    As in the proof of Theorem \ref{strong:sharp}, the vectors $f_I$ form an orthonormal basis of the subspace $\mathcal H_{L,\sigma,r,\epsilon}$, and
    \[
        q_3(E_A)\big|_{\mathcal H_{L,\sigma,r,\epsilon}}=\Bigl[2\sum_{j\in L}t_j^2+\Bigl(\sum_{j\in L}\sigma_jt_j\Bigr)^2\Bigr]\operatorname{Id}-4\mathcal A_r, \qquad \mathcal A_rf_I:=\sum_{i\in I}\sum_{k\notin I\cup L}t_it_kf_{(I\setminus\{i\})\cup\{k\}}.
    \]
    Let $f\in\Lambda^3V^*\subset(\Lambda^3V^*)_{\mathbb C}$ be a real unit form attaining the equality, then as in the proof of Theorem \ref{eqcase:even}, the component $f_{\mathcal H}$ of $f$ in each subspace $\mathcal H$ attains the equality since $E_A\ne0$. But when $r=0$, we have $\mathcal A_r=0$ and $q_3(E_A)\big|_{\mathcal H}\ge0$. Therefore $f$ cannot have a component in $\mathcal H$ when $r=0$. By \eqref{eqcase:degree-count2}, $r=1$ and $s+\epsilon=1$. If $s=1$ and $\epsilon=0$, the subspace is $\mathcal H_{\{l\},\sigma,1,0}$; if $s=0$ and $\epsilon=1$, the subspace is $\mathcal H_{\varnothing,\varnothing,1,1}$. Therefore
    \[
        f=\sum_{l=1}^m\sum_{\sigma=\pm1}f_{\mathcal H_{\{l\},\sigma,1,0}}+f_{\mathcal H_{\varnothing,\varnothing,1,1}}.
    \]

    On $\mathcal H_{\{l\},\sigma,1,0}$, $I$ is the singleton set $\{i\}$ with $i\ne l$, and
    \[
        f_i:=f_I=u_l^{\sigma_l}\wedge\theta_i, \qquad \mathcal A_1f_i=t_i\sum_{k\ne i,l}t_kf_k.
    \]
    Therefore
    \[
        q_3(E_A)f_i=\Bigl[2t_l^2+(\sigma_lt_l)^2\Bigr]f_i-4t_i\sum_{k\ne i,l}t_kf_k=3t_l^2f_i-4t_i\sum_{k\ne i,l}t_kf_k.
    \]
    Write $f_{\mathcal H_{\{l\},\sigma,1,0}}=\sum_{i\ne l}a_if_i$. Then
    \[
        \langle q_3(E_A)f_{\mathcal H_{\{l\},\sigma,1,0}},f_{\mathcal H_{\{l\},\sigma,1,0}}\rangle_{\mathbb C} = 3t_l^2|f_{\mathcal H_{\{l\},\sigma,1,0}}|^2-4\Bigl|\sum_{i\ne l}t_ia_i\Bigr|^2+4\sum_{i\ne l}t_i^2|a_i|^2.
    \]
    Applying the Cauchy--Schwarz inequality twice,
    \[
        \Bigl|\sum_{i\ne l}t_ia_i\Bigr|^2-\sum_{i\ne l}t_i^2|a_i|^2 \le \frac{d-1}d\Bigl|\sum_{i\ne l}t_ia_i\Bigr|^2 \le \frac{d-1}d\Bigl(\sum_{i\ne l}t_i^2\Bigr)\sum_{i\ne l}|a_i|^2.
    \]
    Therefore
    \[
        \langle q_3(E_A)f_{\mathcal H_{\{l\},\sigma,1,0}},f_{\mathcal H_{\{l\},\sigma,1,0}}\rangle_{\mathbb C} \ge \Bigl[3t_l^2-\frac{4(m-2)}{m-1}\Bigl(\sum_{i\ne l}t_i^2\Bigr)\Bigr]|f_{\mathcal H_{\{l\},\sigma,1,0}}|^2.
    \]

    If $n$ is odd, then $M=m$ and
    \[
        -d_nS = -\frac{4(m-1)}m\sum_{i=1}^mt_i^2 < 3t_l^2-\frac{4(m-2)}{m-1}\Bigl(\sum_{i\ne l}t_i^2\Bigr),
    \]
    where the strict inequality holds since the $t_i$ are not all zero, which is guaranteed by $E_A\ne0$. Hence $f$ cannot have a component in $\mathcal H_{\{l\},\sigma,1,0}$ when $n$ is odd; in this case $f=f_{\mathcal H_{\varnothing,\varnothing,1,1}}$. On $\mathcal H_{\varnothing,\varnothing,1,1}$, $I$ is the singleton set $\{i\}$ with $i=1,\cdots,m$, and
    \[
        f_i:=f_I=\theta_i\wedge z, \qquad q_3(E_A)f_i=-4\mathcal A_1f_i=-4t_i\sum_{k\ne i}t_kf_k.
    \]
    Writing $f=\sum_{i=1}^ma_if_i$ and applying the Cauchy--Schwarz inequality twice,
    \[
        \langle q_3(E_A)f,f\rangle_{\mathbb C} = -4\Bigl(\Bigl|\sum_{i=1}^{m}t_ia_i\Bigr|^2-\sum_{i=1}^{m}t_i^2|a_i|^2\Bigr) \ge -\frac{4(m-1)}{m}\sum_{i=1}^{m}t_i^2|f|^2.
    \]
    The equality makes all $t_i$ are equal and $a_i=z_0t_i$ for a nonzero constant $z_0$. Hence $A=t'J$ for an orthogonal complex structure $J$ on $W:=\bigoplus_{i=1}^mP_i$. Let $t=(t')^2$, then $E_A=tE_J$. Moreover,
    \[
        f=\pm\frac1{\sqrt m}\sum_{i=1}^mf_i=\pm\frac1{\sqrt m}\Bigl(\sum_{i=1}^m\theta_i\Bigr)\wedge z=\pm\frac1{\sqrt m}\Omega_J\wedge z=\pm\frac1{\sqrt M}\Omega_J\wedge z.
    \]

    If $n$ is even, then $\epsilon=0$ and $\mathcal H_{\varnothing,\varnothing,1,1}=\varnothing$. Moreover, $M=m-1$ and
    \[
        \Bigl[3t_l^2-\frac{4(m-2)}{m-1}\Bigl(\sum_{i\ne l}t_i^2\Bigr)\Bigr] - (-d_nS) = \left(3+\frac{4(m-2)}{m-1}\right)t_l^2 \ge 0.
    \]
    Since $f\ne0$, $f$ has a component in $\mathcal H_{\{l\},\sigma,1,0}$ for some $l\in\{1,\cdots,m\}$, then the equality makes $t_l=0$. The equalities in the two Cauchy--Schwarz inequalities imply that $t_i=t'>0$ for all $i\ne l$ and $a_{i,\sigma}=z_\sigma t_i$ for a nonzero constant $z_\sigma$ depending on $\sigma$. Hence $A=t'J$ for an orthogonal complex structure $J$ on $W:=\bigoplus_{i\ne l}P_i$ and $E_A=tE_J$. Moreover, if we replace $l$ in the above displayed inequality with $l'\ne l$, then the inequality sign becomes a strict inequality sign since $t_{l'}^2=(t')^2>0$. Therefore $f$ cannot have a component in other $\mathcal H_{\{l'\},\sigma,1,0}$. Now
    \begin{align*}
        f&=f_{\mathcal H_{\{l\},+,1,0}}+f_{\mathcal H_{\{l\},-,1,0}}=\sum_{i\ne l}a_{i,+}f_{i,+}+\sum_{i\ne l}a_{i,-}f_{i,-}=t'\sum_{i\ne l}(z_+f_{i,+}+z_-f_{i,-})\\
        &=t'(z_+u_l^++z_-u_l^-)\wedge\Bigl(\sum_{i\ne l}\theta_i\Bigr)=\eta\wedge\Omega_J=\frac1{\sqrt M}\Omega_J\wedge\zeta,
    \end{align*}
    where $\eta:=t'(z_+u_l^++z_-u_l^-)$, $\zeta:=\sqrt M\eta\in F^*$ and $F:=P_l$. Since $f$ is real, $\zeta$ is also real.

    Now suppose that $E\ne0$. Lemma \ref{strong:decomposition} yields a finite decomposition $E=(1/3)\sum_\nu E_{A_\nu}$. From the proof of Lemma \ref{strong:decomposition}, each $E_{A_\nu}$ is strongly nonnegative. Theorem \ref{strong:sharp} therefore gives $Q_3(E_{A_\nu})\ge -d_nS(E_{A_\nu})$ for every $\nu$. By the linearity of both $q_3$ and the scalar curvature, $Q_3=-d_nS$ is equivalent to $\sum_\nu Q_3(E_{A_\nu})=-d_n\sum_\nu S(E_{A_\nu})$. This is also equivalent to $Q_3(E_{A_\nu})=-d_nS(E_{A_\nu})$ for every $\nu$. If $E_{A_\nu}=0$, it contributes nothing and can be ignored. By the above classification of the equality case for the model tensor, for each $E_{A_\nu}\ne0$,
    \[
        V=W_\nu\oplus F_\nu,\quad \dim W_\nu=2M,\quad E_{A_\nu}=t_\nu E_{J_\nu}, \quad \omega=\frac1{\sqrt M}\Omega_{J_\nu}\wedge\zeta_\nu,
    \]
    where $t_\nu>0$, $J_\nu$ is an orthogonal complex structure on $W_\nu$, extended by zero, and $\zeta_\nu$ is a unit covector in $F_\nu^*$.

    To see that these models coincide up to scale, consider any form of the shape $\omega=\Omega_J\wedge\zeta/\sqrt M$, where $W$ is a $2M$-dimensional subspace, $J$ is an orthogonal complex structure on $W$ extended by zero to $F=W^\perp$, and $\zeta\in F^*$ is a unit covector. We claim that $\{\alpha\in V^* \mid \alpha\wedge\omega=0\}=\mathbb R\zeta$. First, $\zeta\wedge\omega=0$. Conversely, decompose $\alpha=a\zeta+\beta+\eta$, where $a\in\mathbb R$, $\beta\in W^*$, and $\eta\in F^*\cap\zeta^\perp$. Then
    \[
        0=\alpha\wedge\omega=\frac1{\sqrt M}\bigl(a\zeta\wedge\Omega_J\wedge\zeta+\beta\wedge\Omega_J\wedge\zeta+\eta\wedge\Omega_J\wedge\zeta\bigr).
    \]
    The first term vanishes because $\zeta\wedge\zeta=0$. The $F$-degrees of the second and third terms are $1$ and $2$, respectively. Their $F$-degrees are different, so they must vanish separately:
    \[
        \beta\wedge\Omega_J\wedge\zeta=0, \qquad \eta\wedge\Omega_J\wedge\zeta=0.
    \]
    The first equation gives $\beta\wedge\Omega_J=0$ since $\beta\wedge\Omega_J\in\Lambda^3W^*$ and $0\ne\zeta\in F^*$. Similarly, the second equation gives $\eta\wedge\zeta=0$. On $W$, since $\Omega_J$ is nondegenerate, $\beta=\iota_x\Omega_J$ for some $x\in W$. Since $n\ge5$, we have $M-1\ge1$, and hence
    \[
        0=\beta\wedge\Omega_J^{M-1}=(\iota_x\Omega_J)\wedge\Omega_J^{M-1}=\frac1M\iota_x(\Omega_J^M).
    \]
    But $\Omega_J^M$ is a nonzero volume form on $W$, so $x=0$, and $\beta=0$. Also, $\eta\wedge\zeta=0$ implies $\eta=0$ since $\eta\in F^*\cap\zeta^\perp$. Hence $\alpha=a\zeta+\beta+\eta=a\zeta$. This proves the claim. Thus for any $\nu$ and $\nu'$,
    \[
        \mathbb R\zeta_\nu=\{\alpha\in V^* \mid \alpha\wedge\omega=0\}=\mathbb R\zeta_{\nu'}.
    \]
    Therefore, there is a covector $\zeta_0$ such that for all $\nu$, $\zeta_\nu=\epsilon_\nu\zeta_0$ for some $\epsilon_\nu\in\{-1,1\}$. Moreover, $\zeta_0\in F_\nu^*$ for every $\nu$. Also,
    \[
        \omega=\frac1{\sqrt M}\Omega_{J_\nu}\wedge(\epsilon_\nu\zeta_0)=\frac1{\sqrt M}(\epsilon_\nu\Omega_{J_\nu})\wedge\zeta_0, \qquad \epsilon_\nu\Omega_{J_\nu}=\iota_{\zeta_0^\sharp}(\epsilon_\nu\Omega_{J_\nu}\wedge\zeta_0)=\sqrt M\,\iota_{\zeta_0^\sharp}\omega.
    \]
    Since $\sqrt M\,\iota_{\zeta_0^\sharp}\omega$ is independent of $\nu$, $\epsilon_\nu\Omega_{J_\nu}=\epsilon_{\nu'}\Omega_{J_{\nu'}}$ for any $\nu$ and $\nu'$. Equivalently, $\epsilon_\nu J_\nu=\epsilon_{\nu'}J_{\nu'}$ for any $\nu$ and $\nu'$. Let $J_0:=\epsilon_\nu J_\nu$, then $J_0^2=J_\nu^2=-\operatorname{Id}$ on $W_\nu$ and $J_0=0$ on $F_\nu$. But $J_0$ is independent of $\nu$, so all $W_\nu$ are equal to some $W_0$, and hence all $F_\nu=W_0^\perp=:F_0$. Therefore $J_0$ is an orthogonal complex structure on $W_0$, extended by zero. Since $E_J$ is quadratic in $J$, $E_{-J}=E_J$; hence $J_\nu=\epsilon_\nu J_0$ gives $E_{J_\nu}=E_{J_0}$. Finally,
    \[
        \omega=\frac1{\sqrt M}\Omega_{J_0}\wedge\zeta_0, \qquad E=\frac13\sum_\nu E_{A_\nu}=\frac13\sum_\nu(t_\nu E_{J_\nu})=\Bigl(\frac13\sum_\nu t_\nu\Bigr)E_{J_0}=t_0E_{J_0}, \quad t_0>0.
    \]

    The formulas $S=6Mt$ and $q_3(E)(\Omega_J\wedge\zeta)=-4(M-1)t\Omega_J\wedge\zeta$ verify the converse.
\end{proof}

\begin{remark}\label{eqcase:scope}
    In dimension four the coefficient is zero and the classification is different: writing $\omega=*v^\flat$, one has $Q_3=\operatorname{Ric}(E)(v,v)$, so equality means precisely that $v$ is a full curvature-null direction. The tensor on $v^\perp$ can then be any three-dimensional sectional-nonnegative tensor.
\end{remark}

\section{Sectional--scalar pinching and cohomology}\label{geom:section}

\subsection{From an algebraic estimate to a pinching theorem}

\begin{theorem}\label{geom:transfer}
    Let $(M^n,g)$ be a closed, connected Riemannian manifold, and let $1\le p\le n-1$. Suppose that $C\ge0$ satisfies
    \[
        q_p(E)\ge-C\operatorname{Scal}(E)\operatorname{Id}
    \]
    for every algebraic curvature tensor $E$ with nonnegative sectional curvature in dimension $n$. Set
    \begin{equation}\label{geom:beta}
        \beta(C)=\frac{Cn(n-1)}{p(n-p)+Cn(n-1)}.
    \end{equation}
    If $K_{\min}\ge\beta(C)S_0$ everywhere, every harmonic form in degrees $p$ and $n-p$ is parallel. For $k\in\{p,n-p\}$ and every $x\in M$,
    \[
        b_k(M)=\dim(\Lambda^k T_x^*M)^{\operatorname{Hol}_x(g)}\le\binom nk,
    \]
    where $\operatorname{Hol}_x(g)$ is the full holonomy group. If, in addition, the pinching inequality is strict at some point, then $H^p(M;\mathbb R)=H^{n-p}(M;\mathbb R)=0$.
\end{theorem}

\begin{proof}
    With $K_{\min}$ and $S_0$ as in \eqref{geom:normalization}, decompose the curvature at each point as $\operatorname{Rm}=K_{\min}I+E$. The tensor $E$ has nonnegative sectional curvature and scalar curvature $\tau=n(n-1)(S_0-K_{\min})$. By linearity and $q_p(I)=p(n-p)\operatorname{Id}$,
    \begin{align}
        q_p(\operatorname{Rm})&=p(n-p)K_{\min}\operatorname{Id}+q_p(E)\nonumber\\
        &\ge\bigl[p(n-p)K_{\min}-C\tau\bigr]\operatorname{Id}\nonumber\\
        &=\bigl[p(n-p)+Cn(n-1)\bigr]\bigl[K_{\min}-\beta(C)S_0\bigr]\operatorname{Id}\nonumber\\
        &\ge\bigl[K_{\min}-\beta(C)S_0\bigr]\operatorname{Id}.\label{geom:bochner-lower}
    \end{align}
    A choice of orientation at a point identifies $q_{n-p}$ with $q_p$ by Hodge conjugacy (Lemma \ref{lem:HodgeDual}). Hence the same bound holds for $q_{n-p}$ at every point. Only the pointwise values of $E$ enter the estimate.

    Let $k\in\{p,n-p\}$ and let $\omega$ be a harmonic $k$-form, so $\Delta_{\mathrm{H}}\omega=0$. Integrating the Weitzenb\"ock formula $\Delta_{\mathrm{H}}\omega=\nabla^*\nabla\omega+q_p(\mathrm{Rm})\omega$ gives
    \begin{equation}\label{geom:bochner-integral}
        0=\int_M\langle\Delta_{\mathrm{H}}\omega,\omega\rangle\,d\operatorname{vol}_g=\int_M\bigl(|\nabla\omega|^2+\langle q_p(\operatorname{Rm})\omega,\omega\rangle\bigr)\,d\operatorname{vol}_g.
    \end{equation}
    Integration uses the Riemannian density and does not require an orientation. Both summands are nonnegative and continuous, so $\omega$ is parallel. If $K_{\min}>\beta(C)S_0$ at $x$, then $q_p(\operatorname{Rm})>0$ at $x$ by \eqref{geom:bochner-lower}, so the second summand forces $\omega(x)=0$. Parallel transport and connectedness then give $\omega\equiv0$, which means there are no nonzero harmonic $k$-forms. By the Hodge theorem, $H^k(M;\mathbb{R})=0$.

    Conversely, every parallel form is closed and coclosed, and hence harmonic. Thus when $k\in\{p,n-p\}$, we have $\{\text{harmonic }k\text{-forms}\}=\{\text{parallel }k\text{-forms}\}$. We claim that for any $k\in\mathbb N$ and $x\in M$,
    \[
        \{\text{parallel }k\text{-forms}\}\cong\bigl(\Lambda^kT_x^*M\bigr)^{\operatorname{Hol}_x(g)}:=\bigl\{\alpha\in\Lambda^kT_x^*M \mid P_\gamma^*\alpha=\alpha\text{ for all closed loop }\gamma\text{ at }x\bigr\},
    \]
    where $P_\gamma^*:\Lambda^kT_x^*M\to\Lambda^kT_x^*M$ is the induced action on exterior forms of $P_\gamma:T_xM\to T_xM$. Assuming the claim, we obtain
    \[
        H^k(M;\mathbb{R})\cong\{\text{harmonic }k\text{-forms}\}=\{\text{parallel }k\text{-forms}\}\cong\bigl(\Lambda^kT_x^*M\bigr)^{\operatorname{Hol}_x(g)}, \qquad k\in\{p,n-p\},
    \]
    and hence
    \[
        b_k(M) := \dim H^k(M;\mathbb{R}) = \dim\bigl(\Lambda^kT_x^*M\bigr)^{\operatorname{Hol}_x(g)} \le \dim\Lambda^kT_x^*M = \binom nk, \qquad k\in\{p,n-p\},
    \]
    which proves the theorem. Thus it remains to prove the claim.

    The evaluation map from the space of parallel $k$-forms to $\Lambda^kT_x^*M$, given by $\omega\mapsto\omega(x)$, is injective since $\omega(x)=0$ implies $\omega\equiv0$. It maps into the holonomy-invariant subspace: if $\gamma$ is a closed loop at $x$, then, since $\omega$ is parallel, $\omega(x)=P_\gamma^*\omega(x)$; hence $\omega(x)\in\bigl(\Lambda^kT_x^*M\bigr)^{\operatorname{Hol}_x(g)}$. Conversely, suppose $\alpha\in(\Lambda^kT_x^*M)^{\operatorname{Hol}_x(g)}$. For any $y\in M$, choose a path $\gamma$ from $x$ to $y$, and define $\omega(y):=P_\gamma^*\alpha$. This is well defined: if $\gamma'$ is another path from $x$ to $y$, then $\gamma'^{-1}\gamma$ is a closed loop at $x$, so $P_{\gamma'^{-1}\gamma}^*\alpha=\alpha$, and hence $P_{\gamma'}^*\alpha=P_{\gamma'}^*\bigl(P_{\gamma'^{-1}\gamma}^*\alpha\bigr)=P_\gamma^*\alpha$. Thus $\omega$ is a globally well-defined parallel form. Therefore the evaluation map is surjective onto $(\Lambda^kT_x^*M)^{\operatorname{Hol}_x(g)}$ and gives the claimed isomorphism.
\end{proof}

\begin{proof}[Proof of Theorem \ref{intro:geometry}]
    Take $C=\widehat C_{n,p}$ in Theorem \ref{geom:transfer}. Theorem \ref{intro:main} gives its algebraic hypothesis, and \eqref{intro:beta} gives $\beta(C)=\widehat\beta_n^{(p)}$. This proves the vanishing assertion.

    For sharpness, let $n=2m$ and $p=2a$, $1\le a\le m-1$. On $\mathbb{CP}^m$ normalize the Fubini--Study metric so that its curvature tensor is $I+E_J$. For orthonormal $X,Y$ one has
    \[
        K(X,Y)=1+3\langle JX,Y\rangle^2, \qquad K_{\min}=1, \qquad S_0=\frac{n+2}{n-1}.
    \]
    Thus $K_{\min}=\widehat\beta_n^{(p)}S_0$ everywhere by \eqref{geom:even-threshold}. Hence every smaller coefficient in the pinching inequality would make this example satisfy the strict condition and would force $H^p(\mathbb{CP}^m;\mathbb{R})=0$. But the parallel K\"ahler form $\Omega$ has a nonzero harmonic power $\Omega^a$ in degree $p$. Thus $H^p(\mathbb{CP}^m;\mathbb{R})\ne0$. A contradiction.
\end{proof}

\begin{proof}[Proof of Theorem \ref{intro:parallel}]
    With the same choice of $C$, the weak assertion of Theorem \ref{geom:transfer} gives parallelism in both degrees and the formula for their Betti numbers. In the flat case write $M=\mathbb R^n/\Gamma$ and let $F\subset O(n)$ be the finite linear holonomy group. A constant form on $\mathbb R^n$ descends precisely when $F$ fixes it. Thus in every degree
    \begin{equation}\label{geom:flat-betti}
        b_k(M)=\dim\bigl(\Lambda^k(\mathbb R^n)^*\bigr)^F.
    \end{equation}
    For a flat torus $F$ is trivial, and the bound $\binom nk$ is attained.
\end{proof}

\subsection{The even-degree rigidity theorem}

\begin{proof}[Proof of Theorem \ref{geom:even-rigidity}]
    Set $p=2a$ and $C=p(n-p)/(3n)$. Since $b_p(M)>0$, there exists a nonzero harmonic $p$-form $\omega$. By Theorem \ref{intro:parallel}, $\omega$ is parallel, so we may normalize it to have unit length everywhere. The Weitzenb\"ock formula gives $q_p(\operatorname{Rm})\omega=0$. Since $(n-1)S_0/(n+2)\le K_{\min}\le S_0$, we have $S_0\ge0$ and $K_{\min}\ge0$. At a nonflat point, $S_0>0$ and $K_{\min}>0$.

    Fix a nonflat point. If $E=0$, then $q_p(\operatorname{Rm})=q_p(K_{\min}I)=p(n-p)K_{\min}\operatorname{Id}>0$, contradicting $q_p(\operatorname{Rm})\omega=0$. Thus $E\ne0$. By \eqref{geom:bochner-lower}, $q_p(\operatorname{Rm})\ge0$. Hence when \eqref{geom:bochner-lower} is applied to $\omega$, all the inequalities in it become equalities, i.e., $q_p(E)\omega=-C\tau\omega$ and $K_{\min}=(n-1)S_0/(n+2)$. Theorem \ref{eqcase:even} gives an orthogonal complex structure $J$ and $t>0$ such that
    \[
        E=tE_J, \qquad \omega=\pm\frac{\Omega_J^a}{a!\sqrt{\binom ma}}.
    \]
    By Lemma \ref{lem:model-powers}, $q_p(tE_J)\omega=-p(n-p)t\omega$, so
    \[
        0=q_p(\operatorname{Rm})\omega=q_p(K_{\min}I+tE_J)\omega=p(n-p)K_{\min}\omega-p(n-p)t\omega,
    \]
    which gives $K_{\min}=t$. Consequently
    \begin{equation}\label{geom:pointwise-fs}
        \operatorname{Rm}=t(I+E_J), \qquad \operatorname{Ric}=(n+2)tg
    \end{equation}
    at every nonflat point. At a flat point the Ricci tensor is zero, so globally $\operatorname{Ric}=f g$, where the smooth function $f=\operatorname{Scal}_g/n$ is nonnegative. The contracted second Bianchi identity gives $df=\frac12d(nf)$, whence $(n-2)df=0$. Since $n\ge4$ and $M$ is connected, $f$ is constant. In the nonflat case it is positive, so there are no flat points and $t=f/(n+2)>0$ is constant.

    We next verify that the pointwise complex structures admit a smooth choice whose resulting complex structure is parallel; this is not automatic from their pointwise existence. The curvature operator $\mathcal R_{\operatorname{Rm}}$ of \eqref{geom:pointwise-fs} has eigenvalues $(2m+2)t$, $2t$, and $0$. Indeed, if $\Omega_J=\sum_{i=1}^me^{2i-1,2i}$, then the eigenspace for $(2m+2)t$ is the K\"ahler line $\mathbb R\Omega_J$, the eigenspace for $2t$ is the space of primitive real $(1,1)$-forms,
    \[
        \operatorname{span}\{e^{2i-1,2i}-e^{2j-1,2j}, e^{2i-1,2j-1}+e^{2i,2j}, e^{2i-1,2j}-e^{2i,2j-1} \mid i<j\},
    \]
    and the eigenspace for $0$ is the space of real $(2,0)+(0,2)$-forms,
    \[
        \operatorname{span}\{e^{2i-1,2j-1}-e^{2i,2j}, e^{2i-1,2j}+e^{2i,2j-1} \mid i<j\}.
    \]
    Its largest eigenvalue $(2m+2)t$ is simple and separated from the others. Thus its eigenline $\mathbb R\Omega_J$ is a smooth line bundle. On a sufficiently small connected neighborhood choose a smooth generator $\Omega$ of norm $\sqrt m$; its associated endomorphism is a smooth orthogonal complex structure $J'$, by the pointwise classification. At every point $x$ in this neighborhood, $\Omega_x$ and $\Omega_{J,x}$ are either equal or differ only by a sign. Since $\omega$ and $\Omega$ are both smooth, the relation $\omega=\pm\Omega^a/(a!\sqrt{\binom ma})$ has locally constant sign, 
    so differentiation gives $\Omega^{a-1}\wedge\nabla_X\Omega=0$ for every $X$.

    Write $L\alpha=\Omega\wedge\alpha$ and $\Lambda=L^*$. The two-form $\eta=\nabla_X\Omega$ is primitive since $\langle\eta,\Omega\rangle=\frac12X|\Omega|^2=0$; thus $\Lambda\eta=0$. We claim that the elementary Lefschetz identities imply
    \[
        |L^r\eta|^2=\frac{r!(m-2)!}{(m-2-r)!}\,|\eta|^2, \qquad 0\le r\le m-2.
    \]
    Applying $[\Lambda,L]|_{\Lambda^k}=(m-k)\operatorname{Id}$ to $\eta$ and $L^{r-1}\eta$ separately, we obtain $\Lambda L\eta=(m-2)\eta$ and $\Lambda L(L^{r-1}\eta)=L\Lambda L^{r-1}\eta+(m-2r)L^{r-1}\eta$. Hence by induction, $\Lambda L^r\eta=r(m-2-r+1)L^{r-1}\eta$ for $r\ge1$. Let $f_r:=|L^r\eta|^2$. Since $\Lambda$ is the adjoint of $L$, we have
    \[
        f_r=\langle\Lambda L^r\eta,L^{r-1}\eta\rangle=\langle r(m-r-1)L^{r-1}\eta,L^{r-1}\eta\rangle=r(m-r-1)f_{r-1}.
    \]
    Since $f_0=|\eta|^2$, the claim follows. But $L^{a-1}\eta=\Omega^{a-1}\wedge\nabla_X\Omega=0$ and $a-1\le m-2$, so the claimed formula shows $\nabla_X\Omega=0$. Hence $J'$ is locally parallel and \eqref{geom:pointwise-fs}, with constant $t$, is a parallel curvature tensor.

    On the simply connected cover the parallel eigenline has a global parallel unit generator: a metric connection on a real line bundle has zero curvature, and its only possible monodromy is a sign, which is trivial on a simply connected space. The cover is therefore a complete simply connected K\"ahler manifold with curvature tensor $t(I+E_J)$, equivalently with constant positive holomorphic sectional curvature $4t$. The classification of complete simply connected complex space forms identifies it, up to scale, with Fubini--Study $\mathbb{CP}^m$ \cite{Besse1987}.

    It remains to examine the deck transformations, because the conclusion about the cover alone does not exclude an antiholomorphic quotient. Since the covering group $\Gamma$ acts properly discontinuously on compact $\mathbb{CP}^m$, $\Gamma$ is finite. Every isometry preserves the intrinsic K\"ahler eigenline of its curvature operator, and hence sends $\Omega$ to either $\Omega$ or $-\Omega$: if $F$ is an isometry from the cover to itself, then $\mathcal R_{\operatorname{Rm}}\circ F^*=F^*\circ\mathcal R_{\operatorname{Rm}}$, so $F^*\Omega$ is an eigenvector of $\mathcal R_{\operatorname{Rm}}$ with eigenvalue $(2m+2)t$; since $F$ is an isometry, $F^*\Omega$ equals $\Omega$ or $-\Omega$. The sign is constant and defines a homomorphism $\epsilon:\Gamma\to\{1,-1\}$ by
    \[
        \epsilon(F)=
        \begin{cases}
            +1,&F^*\Omega=\Omega,\\
            -1,&F^*\Omega=-\Omega.
        \end{cases}
    \]
    Since $F$ is an isometry, $F^*\Omega=\Omega$ is equivalent to $F_*\circ J=J\circ F_*$, so
    \[
        \ker\epsilon=\{\text{holomorphic deck transformations}\}.
    \]
    Every holomorphic isometry of Fubini--Study space is induced by a unitary linear transformation of $\mathbb C^{m+1}$, so it has a fixed point. Since deck transformations act freely, $\ker\epsilon$ is trivial, so $\Gamma\cong\Gamma/\ker\epsilon\cong\operatorname{im}\epsilon\le\{1,-1\}$. Therefore $\Gamma$ is either trivial or of order two, with its nontrivial element antiholomorphic.

    Suppose $\Gamma$ is of order two. Let $\sigma\in\Gamma$ be the antiholomorphic isometry. Then $\sigma([z])=[A\bar z]$ for some unitary $A$. If $\sigma$ is lifted to a map on $\mathbb C^{m+1}\setminus\{0\}$, a natural lift is $T(z)=A\bar z$. Since $\Gamma$ is of order two, $\sigma^2=\operatorname{Id}_{\mathbb{CP}^m}$, which gives $T^2(z)=\lambda z$, $|\lambda|=1$. Also, $T^2(z)=T(A\bar z)=A\overline Az$. Thus $A\overline A=\lambda\operatorname{Id}$ and $\overline A=\lambda A^{-1}$. This implies $A\overline A=\lambda\operatorname{Id}=\overline AA$, so $\lambda=\overline\lambda$, and $\lambda=\pm1$. We claim that $\lambda=-1$. If $\lambda=1$, then $T^2=\operatorname{Id}$ and there is a nonzero fixed vector $w$. Indeed, for any $z\ne0$, $T(z+Tz)=z+Tz$; hence if $z+Tz\ne0$, take $w=z+Tz$; if $z+Tz=0$, then $T(iz)=\bar i\,Tz=-i(-z)=iz$ and take $w=iz$. Thus $\sigma([w])=[w]$, contrary to freeness. Thus the claim follows and $T^2=-\operatorname{Id}$. In this case $\sigma$ has no fixed point; otherwise $Tz=cz$ would imply $-z=T^2z=|c|^2z$. Since $A\overline A=-\operatorname{Id}$, $|\det A|^2=(-1)^{m+1}$, which implies that $m+1$ is even, so $m=2s+1$.

    Since $A$ is unitary, $T$ is antiunitary. An orthonormal pair $u,Tu$ spans a $T$-invariant complex two-plane $\operatorname{span}_{\mathbb C}\{u,Tu\}$. These vectors are orthogonal by antiunitarity and $T^2=-\operatorname{Id}$. Its orthogonal complement is also $T$-invariant. Repeating this construction gives a unitary basis in which
    \[
        T(z_1,z_2,\cdots,z_m,z_{m+1})=(-\bar z_2,\bar z_1,\cdots,-\bar z_{m+1},\bar z_m).
    \]

    Finally $\sigma^*(\Omega^r)=(\sigma^*\Omega)^r=(-\Omega)^r=(-1)^r\Omega^r$, so $\Omega^r$ can descend to the quotient manifold if and only if $(-1)^r=1$, i.e., if and only if $r$ is even. On $\mathbb{CP}^m$ each even-degree harmonic space is $\mathbb R\Omega^r$ and odd-degree cohomology vanishes. Since the covering map is a local isometry, its pullback map commutes with $\Delta_{\mathrm{H}}$. Hence the pullback identifies harmonic forms on a finite Riemannian quotient with invariant harmonic forms on its cover:
    \[
        \mathcal H^k(M)\cong\{\beta\in\mathcal H^k(\mathbb{CP}^m) \mid \sigma^*\beta=\beta\}.
    \]
    This proves the stated Betti numbers and $b_{2a}=1$. Since $m$ is odd, $\sigma^*(\Omega^m)=-\Omega^m$, so the quotient is nonorientable.
\end{proof}

\begin{remark}\label{geom:orientation-warning}
    The quotient $\mathbb{CP}^3/\langle\sigma\rangle$ has real dimension $6$ but $b_4=1\ne0=b_2$. Thus equality of the complementary ordinary Betti numbers must not be inferred from Hodge duality on a nonorientable manifold. The local algebraic conjugacy of the curvature endomorphisms still proves the simultaneous strict vanishing in Theorem \ref{geom:transfer}; that argument does not assert an isomorphism between the two ordinary cohomology groups at a weak endpoint.
\end{remark}

\subsection{Holonomy and the degree-three endpoint}

\begin{proof}[Proof of Theorem \ref{geom:three-holonomy}]
    First, $\beta S_0\le K_{\min}\le S_0$ implies $S_0\ge0$ and $K_{\min}\ge0$. At a nonflat point, $S_0>0$ and $K_{\min}>0$. The universal cover is complete and simply connected. If its restricted holonomy representation were reducible, the de Rham decomposition would split it isometrically as a nontrivial Riemannian product \cite{deRham1952}. Every point of a nontrivial product has a mixed two-plane of sectional curvature zero. This contradicts the fact that $K_{\min}>0$ at a lift of a nonflat point. Therefore the holonomy representation of the cover is irreducible.

    If this irreducible cover is not locally symmetric, the Riemannian holonomy classification gives
    \[
        \operatorname{SO}(n), \quad U(m), \quad SU(m), \quad \operatorname{Sp}(m), \quad \operatorname{Sp}(m)\operatorname{Sp}(1), \quad G_2, \quad \operatorname{Spin}(7),
    \]
    with their standard dimension restrictions. We use the complete modern form of Berger's classification, as stated in \cite[Theorem~2.3]{Joyce2004}, together with \cite{Berger1955}. The groups $SU(m)$, $\operatorname{Sp}(m)$, $G_2$, and $\operatorname{Spin}(7)$ force Ricci-flatness. Nonnegative sectional curvature and $\operatorname{Ric}=0$ then force $\operatorname{Rm}=0$ by Lemma \ref{eqcase:nullity}, contradicting the nonflat case. The tangent representations of $U(m)$ and $\operatorname{Sp}(m)\operatorname{Sp}(1)$ contain $-\operatorname{Id}$. Its pullback on $k$-forms is $(-\operatorname{Id})^*|_{\Lambda^kT_x^*M}=(-1)^k\operatorname{Id}_{\Lambda^kT_x^*M}$ since
    \[
        \bigl((-\operatorname{Id})^*\omega\bigr)(v_1,\cdots,v_k)=\omega(-v_1,\cdots,-v_k)=(-1)^k\omega(v_1,\cdots,v_k).
    \]
    In particular $(-\operatorname{Id})^*|_{\Lambda^3T_x^*M}=-\operatorname{Id}_{\Lambda^3T_x^*M}$, so these groups have no nonzero invariant three-form. Finally, we claim that $\operatorname{SO}(n)$ has no nonzero invariant $k$-form when $1\le k\le n-1$. Indeed, let $\omega\in\Lambda^kT_x^*M$ satisfy $h^*\omega=\omega$ for all $h\in\operatorname{SO}(n)$. Fix $k$ distinct indices $i_1,\cdots,i_k$. Since $k\le n-1$, choose $i_0\notin\{i_1,\cdots,i_k\}$. Define $h\in\operatorname{SO}(n)$ by
    \[
        he_{i_0}=-e_{i_0}, \qquad he_{i_1}=-e_{i_1}, \qquad he_j=e_j \quad (j\ne i_0,i_1).
    \]
    Then $\det h=(-1)^2=1$, so $h\in\operatorname{SO}(n)$. Using $h^*\omega=\omega$, we obtain
    \begin{align*}
        \omega_{i_1\cdots i_k}&=\omega(e_{i_1},\cdots,e_{i_k})=(h^*\omega)(e_{i_1},\cdots,e_{i_k})\\
        &=\omega(he_{i_1},\cdots,he_{i_k})=\omega(-e_{i_1},e_{i_2},\cdots,e_{i_k})=-\omega_{i_1\cdots i_k}.
    \end{align*}
    Hence $\omega_{i_1\cdots i_k}=0$ and $\omega=0$. This proves the claim. Since $3\le n-1$, the cover has no nonzero parallel three-form.

    If the irreducible cover is locally symmetric, it is a complete simply connected irreducible symmetric space. Positive curvature at the chosen point forces compact type and rank one, since a space of rank at least two has a totally geodesic flat two-plane through each point, again contrary to $K_{\min}>0$ at the chosen point. Thus it is a compact rank-one symmetric space: a sphere, a complex projective space, a quaternionic projective space, or the Cayley plane; see \cite{Besse1987}. A sphere of dimension at least four has holonomy $\operatorname{SO}(n)$ and no nonzero parallel three-form. The isotropy representations of the other three families contain $-\operatorname{Id}$, so they have no nonzero invariant odd-degree forms. For the Cayley plane this is the central element $-1$ of $\operatorname{Spin}(9)$ acting as $-\operatorname{Id}$ on its real $16$-dimensional spin representation. The complex and quaternionic cases have the same central tangent action already used above. Thus the cover has no nonzero parallel three-form in either case.

    Next suppose $n$ is odd. In the non-locally-symmetric case, the groups $U(m)$ and $\operatorname{Sp}(m)\operatorname{Sp}(1)$ have even dimension, while $SU(m)$, $\operatorname{Sp}(m)$, $G_2$, and $\operatorname{Spin}(7)$ have already been excluded by Ricci-flatness, hence only $\operatorname{SO}(n)$ remains. In the locally-symmetric case, the complex projective, quaternionic projective, and Cayley planes have even dimension, so the only odd-dimensional one is the sphere, whose holonomy is $\operatorname{SO}(n)$.

    Finally, if $q_3(\operatorname{Rm})\ge0$, then \eqref{geom:bochner-integral} makes every harmonic three-form parallel. A parallel form on $M$ pulls back to a parallel form on its cover, but the cover has no nonzero parallel three-form in the nonflat case, so $H^3(M;\mathbb R)=0$. The simply connected cover is oriented; a parallel $(n-3)$-form there would have a parallel Hodge-dual three-form. The local Hodge conjugacy also makes every harmonic $(n-3)$-form on $M$ parallel. Pulling it back then proves $H^{n-3}(M;\mathbb R)=0$, without an orientability assumption on $M$.
\end{proof}

\begin{proof}[Proof of Theorem \ref{geom:three-rigidity}]
    For $n\ge5$, the coefficient $\widehat C_{n,3}$ is positive and finite, so $0<\widehat\beta_n^{(3)}<1$. Under the weak pinching hypothesis, \eqref{geom:bochner-lower} with $C=\widehat C_{n,3}$ gives $q_3(\operatorname{Rm})\ge0$. Theorem \ref{geom:three-holonomy} gives the conclusion. For $n=5,6,7,8$ the coefficients $\widehat C_{n,3}$ are, respectively, $1/3$, $1/3$, $5168/11253$, and $175/344$. Substitution in \eqref{intro:beta} gives the stated thresholds.
\end{proof}

In the flat case \eqref{geom:flat-betti} gives $b_3>0$ precisely when the linear holonomy fixes a nonzero three-form. Connectedness matters for the nonflat conclusion. A disjoint union of a round sphere and flat tori satisfies the weak inequality and has arbitrarily large $b_3$ as the number of tori increases.

\subsection{Rational homology spheres in odd dimensions}

\begin{proof}[Proof of Corollary \ref{geom:odd-rational-sphere}]
    For $n\ge5$ we have $0<B_n<1$. Since $K_{\min}\le S_0$, the weak hypothesis gives $S_0\ge0$ and $K_{\min}\ge0$. By complementary-degree symmetry of \eqref{intro:beta} and $\widehat\beta_n^{(1)}=0$, we have $\widehat\beta_n^{(p)}\le B_n$ for each $1\le p\le n-1$. Thus $K_{\min}\ge B_nS_0\ge\widehat\beta_n^{(p)}S_0$. Theorem \ref{geom:transfer} then makes every harmonic form of degree $1,\cdots,n-1$ parallel. Since $M$ is nonflat and odd-dimensional, Theorem \ref{geom:three-holonomy} gives restricted holonomy $\operatorname{SO}(n)$, which has no nonzero invariant $p$-form from the proof of Theorem \ref{geom:three-holonomy}. Thus $b_1=\cdots=b_{n-1}=0$. For $n=3$, the strict condition is $K_{\min}>B_3S_0=0$. Theorem \ref{geom:transfer} with $p=1$, $C=0$ applies, since $q_1(E)=\operatorname{Ric}(E)^\sharp$ is nonnegative when $\sec(E)\ge0$. It gives $b_1=b_2=0$.

    The product $S^2\times S^1$ with its standard product metric has $K_{\min}=0$ and $b_1=1$, which shows why the nonflat weak assertion does not extend to $B_3=0$.

    Every closed odd-dimensional manifold has Euler characteristic zero by Poincar\'e duality over $\mathbb F_2$. Computing this same characteristic over $\mathbb R$ gives $0=1-b_n$ since the connectedness of $M$ gives $b_0=1$. Extension of scalars from $\mathbb Q$ to $\mathbb R$ now gives $H_*(M;\mathbb Q)\cong H_*(S^n;\mathbb Q)$.
\end{proof}

\section{Optimal constants and the remaining odd-degree problem}
\label{sec:open}

We finish by specifying the scope of the sharpness statements. For $0\le p\le n$, let
\[
    \gamma_{n,p}^*=\inf\{\gamma\ge0 \mid q_p(E)\ge-\gamma\operatorname{Scal}(E)\operatorname{Id}\text{ for every }\sec(E)\ge0\}.
\]
The scalar normalization used in this definition is harmless. Indeed, if $\sec(E)\ge0$ and $\operatorname{Scal}(E)=0$, then $\sec(E)=0$, hence $E=0$. The normalized set
\[
    \mathcal K_n=\{E \mid \sec(E)\ge0,\ \operatorname{Scal}(E)=1\}
\]
is compact. Indeed, in a fixed frame each $K_{ij}\in[0,1/2]$ since $\sum_{i<j}K_{ij}=1/2$. Jacobi--Cauchy--Schwarz \eqref{alg:Jacobi-CS} bounds all curvature components with three distinct indices, and \eqref{alg:NW-equation} bounds those with four. There are only finitely many components, so all components of $E$ are uniformly bounded. Thus $\mathcal K_n$ is bounded, and it is closed by continuity of sectional and scalar curvature. Thus the continuous energy $Q_p(E,\omega)$ attains a minimum on $\mathcal K_n\times\{\omega \mid |\omega|=1\}$. For $2\le p\le n-2$, the model calculation below shows that the minimum is negative; its negative is then the optimal coefficient.

The proved general comparison is
\begin{equation}\label{open:interval}
    \frac{2ab}{3(a+b)}\le\gamma_{n,p}^*\le\widehat C_{n,p}, \qquad a=\lfloor p/2\rfloor,\quad b=\lfloor(n-p)/2\rfloor,
\end{equation}
When $a+b=0$ the left side means zero. The lower bound follows from the sharp examples for Theorem~\ref{strong:sharp}, which also have nonnegative sectional curvature.

\begin{conjecture}\label{open:conjecture}
    For every $n$ and $p$, the optimal nonnegative coefficient on the full sectional-nonnegative cone is
    \[
        \gamma_{n,p}^*=\frac{2\lfloor p/2\rfloor\lfloor(n-p)/2\rfloor}{3\bigl(\lfloor p/2\rfloor+\lfloor(n-p)/2\rfloor\bigr)},
    \]
    with the zero-denominator case interpreted as zero.
\end{conjecture}

This is a conjecture only in the branches not proved above. The two-form case and its Hodge dual, every even-dimensional even-degree case, and degree three in dimensions at most six are already theorems. On the strongly nonnegative cone it holds in every degree. The appendix supplies another class on which the degree-three value is sharp without strengthening the curvature hypothesis.

\appendix
\section{A sharp estimate for a family with a three-dimensional remainder}
\label{sec:block-family}

The estimates above apply to arbitrary forms. We now isolate a special family of three-forms on which the conjectured degree-three constant follows directly from the NW inequality \eqref{alg:NW-equation} and two scalar Cauchy--Schwarz inequalities. The covectors in the three-dimensional remainder need not be parallel.

\begin{theorem}\label{block:sharp}
    Let $m\ge2$ and let
    \[
        H^{2m+1}=P_1^2\oplus\cdots\oplus P_{m-1}^2\oplus U^3
    \]
    be an orthogonal decomposition. Choose unit area forms $\eta_a$ on $P_a$ and an orientation of $U$. Suppose
    \begin{equation}\label{block:form}
        \omega=\sum_{a=1}^{m-1}\eta_a\wedge u_a+c\operatorname{vol}_U, \qquad u_a\in U^*, \quad c\in\mathbb R.
    \end{equation}
    If $\sec(E)\ge0$, then
    \begin{equation}\label{block:bound}
        Q_3(E,\omega)\ge-\frac{2(m-1)}{3m}\operatorname{Scal}(E)|\omega|^2.
    \end{equation}
    The same assertion holds on $V^{2m+2}=H\oplus\mathbb Rz$ for forms of the form \eqref{block:form}, extended by zero to $V$. The constant is sharp in both dimensions. If $E\ne0$ and $\omega\ne0$, equality holds precisely when, for some $2m$-plane $W\subset V$, orthogonal complex structure $J$ on $W$, unit $\zeta\in(W^\perp)^*$, and numbers $t>0$, $\lambda\ne0$,
    \[
        E=tE_J, \qquad \omega=\lambda\Omega_J\wedge\zeta,
    \]
    where $E_J$ is extended by zero from $W$ and the resulting form belongs to the family \eqref{block:form}.
\end{theorem}

\begin{proof}
    We first work in $H$. Choose oriented orthonormal bases $(e_{a1},e_{a2})$ of $P_a$. Put
    \begin{align*}
        k_a&=E(e_{a1},e_{a2},e_{a1},e_{a2}),& r_{ab}&=\langle\mathcal R_E\eta_a,\eta_b\rangle,\\
        A_{ab}&=\sum_{\alpha,\beta=1}^2E(e_{a\alpha},e_{b\beta},e_{a\alpha},e_{b\beta}),& \mathsf P_a(x,y)&=\frac12\sum_{\alpha=1}^2E(e_{a\alpha},x,e_{a\alpha},y).
    \end{align*}
    Here $\mathsf P_a$ is a symmetric form on $U$. Use the metric to identify covectors and vectors, and define
    \[
        D(x,y)=\langle\mathcal R_E(*_Ux^\flat),*_Uy^\flat\rangle, \quad d=\operatorname{tr} D, \quad \langle h_a,x\rangle=\langle\mathcal R_E\eta_a,*_Ux^\flat\rangle.
    \]
    We have $k_a,A_{ab}\ge0$ and $\mathsf P_a\ge0$. Since $\dim U=3$, every two-form on $U$ is decomposable, so sectional nonnegativity also gives $D\ge0$; thus $D\le d\operatorname{Id}_U$. This argument does not require $\mathcal R_E\ge0$ on the full exterior square. Counting all coordinate sectional curvatures once gives
    \begin{equation}\label{block:scalar}
        S=2\Bigl(\sum_a k_a+\sum_{a<b}A_{ab}+2\sum_a\operatorname{tr}\mathsf P_a+d\Bigr).
    \end{equation}

    Write $a_a=|u_a|$ and $N=|\omega|^2=\sum_a a_a^2+c^2$. Since $\eta_a=e^{a1}\wedge e^{a2}$, we have $\iota_{e_{a1}}\omega=e^{a2}\wedge u_a$ and $\iota_{e_{a2}}\omega=-e^{a1}\wedge u_a$. For an oriented orthonormal basis $(f_i)_{i=1}^3$ of $U$,
    \[
        \iota_{f_i}\omega=\sum_a u_a(f_i)\eta_a+c*_Uf^i.
    \]
    It follows, with every mixed entry retained, that
    \begin{align}
        \sum_j\langle\mathcal R_E\iota_j\omega,\iota_j\omega\rangle &= \sum_a\sum_{\alpha=1}^2\langle\mathcal R_E\iota_{e_{a\alpha}}\omega,\iota_{e_{a\alpha}}\omega\rangle+\sum_{i=1}^3\langle\mathcal R_E\iota_{f_i}\omega,\iota_{f_i}\omega\rangle \notag\\
        &= 2\sum_a\mathsf P_a(u_a,u_a)+\sum_a k_aa_a^2+c^2d+2\sum_{a<b}r_{ab}\langle u_a,u_b\rangle+2c\sum_a\langle h_a,u_a\rangle. \label{block:pair-energy}
    \end{align}
    The contraction Gram matrix $T$ is $a_a^2\operatorname{Id}$ on $P_a$ and $\sum_a u_a\otimes u_a+c^2\operatorname{Id}$ on $U$, with all other blocks zero. On the curvature side,
    \[
        \operatorname{tr}_{P_a}\operatorname{Ric}=2k_a+\sum_{b\ne a}A_{ab}+2\operatorname{tr}\mathsf P_a, \qquad \operatorname{Ric}|_U=2\sum_a\mathsf P_a+d\operatorname{Id}-D.
    \]
    Substituting these expressions and \eqref{block:pair-energy} in \eqref{eq:three-energy} gives the exact identity
    \begin{align}
        Q_3 &= \sum_{a<b}\left[A_{ab}(a_a^2+a_b^2)-4r_{ab}\langle u_a,u_b\rangle\right] \notag\\
        &\quad +2\sum_a\Bigl[(a_a^2+c^2)\operatorname{tr}\mathsf P_a-\mathsf P_a(u_a,u_a)+\sum_{b\ne a}\mathsf P_a(u_b,u_b)\Bigr] \notag\\
        &\quad +d\sum_a a_a^2-\sum_aD(u_a,u_a)-4c\sum_a\langle h_a,u_a\rangle. \label{block:exact}
    \end{align}
    In particular, arbitrary directions of the $u_a$ introduce no unaccounted Ricci terms.

    First, $\mathsf P_a\ge0$ gives $\sum_{b\ne a}\mathsf P_a(u_b,u_b)\ge0$. Second, $D\le d\operatorname{Id}$ makes $da_a^2-D(u_a,u_a)\ge0$. Lemma \ref{alg:nw} gives $6|r_{ab}|\le A_{ab}+4\sqrt{k_ak_b}$ for each pair $a<b$, so $-4r_{ab}\langle u_a,u_b\rangle\ge-4|r_{ab}||\langle u_a,u_b\rangle|\ge-4(\frac16A_{ab}+\frac23\sqrt{k_ak_b})|\langle u_a,u_b\rangle|$. Substituting this into the first line of \eqref{block:exact}, the coefficient of $A_{ab}$ is $a_a^2+a_b^2-\frac23|\langle u_a,u_b\rangle|\ge0$. The remaining negative term is $-\frac83\sqrt{k_ak_b}|\langle u_a,u_b\rangle|$. Hence
    \begin{equation}\label{block:estimate}
        Q_3 \ge -\frac83\sum_{a<b}\sqrt{k_ak_b}|\langle u_a,u_b\rangle|+2\sum_a\Bigl[(a_a^2+c^2)\operatorname{tr}\mathsf P_a-\mathsf P_a(u_a,u_a)-2|c||\langle h_a,u_a\rangle|\Bigr].
    \end{equation}

    For $a_a>0$, apply the same lemma to the four-plane $P_a\oplus(u_a^\perp\cap U)$. The sectional curvature of the two-plane $u_a^\perp\cap U$ is $D(u_a,u_a)/a_a^2$, and the sum of four transverse sectional curvatures of $P_a\oplus(u_a^\perp\cap U)$ is $2\operatorname{tr}\mathsf P_a-2\mathsf P_a(u_a,u_a)/a_a^2$. Therefore Lemma \ref{alg:nw} gives
    \[
        6|\langle h_a,u_a\rangle| \le 2a_a\operatorname{tr}\mathsf P_a-\frac{2}{a_a}\mathsf P_a(u_a,u_a)+4\sqrt{k_aD(u_a,u_a)}.
    \]
    Using this to replace $-2|c||\langle h_a,u_a\rangle|$ in \eqref{block:estimate} yields the additional term $-\frac43|c|\sqrt{k_aD(u_a,u_a)}$, while the remaining $\mathsf P_a$-terms are
    \begin{align*}
        &\quad\left(a_a^2+c^2-\frac{2|c|a_a}{3}\right)\operatorname{tr}\mathsf P_a+\left(\frac{2|c|}{3a_a}-1\right)\mathsf P_a(u_a,u_a)\\
        &=\frac{c^2}{a_a^2}\mathsf P_a(u_a,u_a)+\left(a_a^2+c^2-\frac23|c|a_a\right)\left(\operatorname{tr}\mathsf P_a-\frac1{a_a^2}\mathsf P_a(u_a,u_a)\right)\ge0.
    \end{align*}
    If $a_a=0$, the brackets in \eqref{block:estimate} is nonnegative, while $\sqrt{k_aD(u_a,u_a)}=0$. We finally obtain
    \begin{equation}\label{block:negative-budget}
        Q_3\ge-\frac83\Bigl\{\sum_{a<b}\sqrt{k_ak_b}|\langle u_a,u_b\rangle|+|c|\sum_a\sqrt{k_aD(u_a,u_a)}\Bigr\}.
    \end{equation}

    Set $x_a=a_a\sqrt{k_a}$ for $a<m$ and $x_m=|c|\sqrt d$. The expression in braces is bounded above by $\sum_{i<j}x_ix_j$. Two applications of Cauchy--Schwarz give
    \begin{align}
        \sum_{i<j}x_ix_j &= \frac12\Bigl[\Bigl(\sum_i x_i\Bigr)^2-\sum_i x_i^2\Bigr]\le\frac{m-1}{2m}\Bigl(\sum_i x_i\Bigr)^2 \notag\\
        &\le \frac{m-1}{2m}\Bigl(\sum_a a_a^2+c^2\Bigr)\Bigl(\sum_a k_a+d\Bigr)\le\frac{m-1}{4m}NS. \label{block:cs}
    \end{align}
    The last inequality uses \eqref{block:scalar}. Equations \eqref{block:negative-budget} and \eqref{block:cs} prove \eqref{block:bound}.

    We next keep track of equality, assuming $E\ne0$ and $\omega\ne0$. Then $S,N>0$. Every inequality in \eqref{block:cs} must be an equality, and its first one forces $x_1=\cdots=x_m>0$; its second one forces $(a_1,\cdots,a_{m-1},|c|)$ proportional to $(\sqrt{k_1},\cdots,\sqrt{k_{m-1}},\sqrt d)$; its last one forces $A_{ab}=0$ and $\mathsf P_a=0$. Combining the first two requirements yields numbers $a_0,k_0>0$ such that
    \[
        a_a=|c|=a_0, \quad k_a=d=k_0, \qquad A_{ab}=0, \quad \mathsf P_a=0.
    \]
    Equality in $D(u_a,u_a)\le d a_a^2$ implies that $D$ has an eigenvalue $d$, but its trace is $d$, so all its remaining eigenvalues are zero since $D\ge0$. Thus there is a unit $\zeta\in U^*$ such that
    \[
        D=k_0\zeta\otimes\zeta, \qquad u_a=\epsilon_a a_0\zeta, \quad \epsilon_a\in\{1,-1\}.
    \]
    Let $P_m=\zeta^\perp\cap U$ and $\eta_m=*_U\zeta$. The form is
    \[
        \omega=a_0\Bigl(\sum_{a<m}\epsilon_a\eta_a+\operatorname{sgn}(c)\eta_m\Bigr)\wedge\zeta.
    \]
    The expression in parentheses is the K\"ahler form of an orthogonal complex structure $J$ on $W=P_1\oplus\cdots\oplus P_m$. Moreover, since $\operatorname{Ric}|_U=2\sum_a\mathsf P_a+d\operatorname{Id}-D$, we have $\operatorname{Ric}(E)(\zeta,\zeta)=d-D(\zeta,\zeta)=0$, so Lemma \ref{eqcase:ricci-kernel} makes $\zeta$ a full curvature null direction.

    In the resulting block frame of $W$, the only nonzero coordinate sectional curvatures are $k_0$ on each $P_i$. As in the proof of Lemma \ref{eqcase:two-form-endpoint}, every three-index component is zero, and every four-index component whose four indices do not make two complete blocks is zero. For two complete blocks, combining equality in the appropriate use of NW, the sign of the mixed energy, $A_{ab}=0$ and $\mathsf P_a=0$, gives $E(\eta_i,\eta_j)=\frac23 k_0\epsilon_i\epsilon_j$ and $\epsilon_m=\operatorname{sgn}(c)$. As in the proof of Lemma \ref{eqcase:two-form-endpoint}, using Lemma \ref{eqcase:four-dimensional} we have $E=(k_0/3)E_J$. Conversely Lemma \ref{lem:model-powers}, with $r=1$, shows that these tensors and forms attain equality, and proves sharpness.

    For the even-dimensional assertion, write $J_z(x,y)=E(z,x,z,y)$ on $H$. For $\omega\in\Lambda^3H^*$, the definition of $q_3$ gives
    \[
        Q_3^V(E,\omega)=Q_3^H(E|_H,\omega)+\langle J_z^{[3]}\omega,\omega\rangle, \qquad S_V=S_H+2\operatorname{tr} J_z.
    \]
    Both $\langle J_z^{[3]}\omega,\omega\rangle$ and $\operatorname{tr} J_z$ are nonnegative, so $Q_3^H\ge-2(m-1)/(3m)S_H|\omega|^2$ implies $Q_3^V\ge-2(m-1)/(3m)S_V|\omega|^2$. Equality forces $\operatorname{Ric}(E)(z,z)=\operatorname{tr}J_z=0$. Lemma \ref{eqcase:ricci-kernel} now eliminates all components involving $z$, and the odd-dimensional equality classification applies. There are therefore exactly two curvature-null directions, $\zeta$ and $z$, in the nonzero even-dimensional model, and exactly one, $\zeta$, in the odd-dimensional model.
\end{proof}

\end{document}